\documentclass[11pt,reqno]{amsart}

\usepackage{amssymb, amsmath, amsfonts, latexsym}
\usepackage{enumerate}
\usepackage{color}
\usepackage{threeparttable}
\usepackage{dcolumn}
\usepackage{multirow,tabularx}
\usepackage{rotate,graphics,epsfig}
\usepackage{booktabs}
\usepackage{float}
\usepackage{bm}
\usepackage{epstopdf}

\newcommand{\bd}{\bold}

\newtheorem{thm}{Theorem}[]
\newtheorem{lem}{Lemma}[section]
\newtheorem{cor}{Corollary}[section]
\newtheorem{prop}{Proposition}[section]

\theoremstyle{definition}
\newtheorem{defn}{Definition}[section]

\numberwithin{equation}{section} \theoremstyle{remark}

\title[Random Orthogonal Matrices and Independent Normals]{\bf Distances between Random Orthogonal Matrices and Independent Normals }

\author{Tiefeng JIANG}
\address{Tiefeng JIANG \\ School of Statistics, University of Minnesota, 224 Church Street SE, Minneapolis, MN 55455 USA.}
\thanks{The research of Tiefeng Jiang was supported in part by NSF Grant DMS-1209166
and DMS-1406279.}
\email{jiang040@umn.edu}

\author{Yutao Ma}
\address{Yutao MA\\ School of Mathematical Sciences $\&$ Laboratory  of Mathematics and Complex Systems of Ministry of Education, Beijing Normal University, 100875 Beijing, China.}
\thanks{The research of Yutao Ma was supported in part by NSFC 11431014, 11371283, 11571043 and 985 Projects.}
\email{mayt@bnu.edu.cn}

\newcommand{\ee}{\mathbb{E}}

\def\ml{\mathcal}
\def\f{\mathbf}

\def\lbl{\label}

\def\<{\langle}
\def\>{\rangle}

\def\beaa{\begin{eqnarray*}}
\def\eeaa{\end{eqnarray*}}
\def\bea{\begin{eqnarray}}
\def\eea{\end{eqnarray}}
\def\lbl{\label}

\def\be{\begin{equation}}
\def\ee{\end{equation}}

\def\bdef{\begin{defn}}
\def\ndef{\end{defn}}

\def\bthm{\begin{thm}}
\def\nthm{\end{thm}}

\def\bprop{\begin{prop}}
\def\nprop{\end{prop}}

\def\brmk{\begin{remarks}}
\def\nrmk{\end{remarks}}

\def\bexa{\begin{exa}}
\def\nexa{\end{exa}}

\def\blem{\begin{lem}}
\def\nlem{\end{lem}}

\def\bcor{\begin{cor}}
\def\ncor{\end{cor}}

\def\bexe{\begin{exe}}
\def\nexe{\end{exe}}

\def\bprf{\begin{proof}}
\def\nprf{\end{proof}}

\def\bdes{\begin{description}}
\def\ndes{\end{description}}

\begin{document}
\maketitle

\begin{abstract}
Let $\bold{\Gamma}_n$ be an $n\times n$ Haar-invariant orthogonal matrix.  Let $\f Z_n$ be the $p\times q$
upper-left submatrix of
$\bold{\Gamma}_n,$ where $p=p_n$ and $q=q_n$ are two positive integers.
Let $\f G_n$ be a $p\times q$ matrix whose $pq$ entries are independent standard normals.
In this paper we consider the distance between $\sqrt{n}\f Z_n$ and $\f G_n$ in terms of the total
variation distance, the Kullback-Leibler distance, the Hellinger distance and the Euclidean distance.
We prove that each of the first three distances  goes to zero as long as $pq/n$ goes to zero, and not so
if  $(p, q)$ sits on the curve $pq=\sigma n$, where $\sigma$ is a constant. However, it is different for the Euclidean
distance, which goes to zero provided $pq^2/n$ goes to zero, and not so if  $(p,q)$ sits on the curve
$pq^2=\sigma n.$ A previous work by Jiang \cite{Jiang06}
shows that the total variation distance goes to
zero if both $p/\sqrt{n}$ and $q/\sqrt{n}$ go to zero, and it is not true provided
$p=c\sqrt{n}$ and $q=d\sqrt{n}$ with $c$ and $d$ being constants. One of the above results confirms a conjecture that the total variation distance goes to zero as long as $pq/n\to 0$ and the distance does not go to zero if $pq=\sigma n$ for some constant $\sigma$.
\end{abstract}

\bigskip\bigskip\bigskip\bigskip\bigskip



\newpage

\section{Introduction}\label{chap:intro}

Let $O(n)$ be the orthogonal group consisting of all $n\times n$ orthogonal matrices.
Let $\bold{\Gamma}_n=(\gamma_{ij})_{n\times n}$ be a random orthogonal matrix which is
uniformly distributed on the orthogonal group $O(n)$, or equivalently, $\bold{\Gamma}_n$ follows the Haar-invariant probability measure on $O(n)$. We sometimes also say that $\bold{\Gamma}_n$ is an Haar-invariant orthogonal matrix. Let $\f Z_n$ be the $p\times q$
upper-left submatrix of
$\bold{\Gamma}_n,$ where $p=p_n$ and $q=q_n$ are two positive integers.
Let $\f G_n$ be a $p\times q$ matrix from which the $pq$ entries are independent standard normals.
In this paper we will study the distance between  $\sqrt{n}\f Z_n$ and $\f G_n$ in terms of the
total variation distance, the Hellinger distance, the Kullback-Leibler distance and the
Euclidean distance (or equivalently, the trace norm). Throughout this paper, we will frequently encounter the notations $p_n, q_n.$ For simplicity, we will use  $p$ and $q$ rather than $p_n$ and $q_n$, respectively, if there is no confusion.

It has long been observed that the entries of $\bold{\Gamma}_n$ are roughly independent random variables with distribution $N(0, \frac{1}{n}).$ Historically, authors show that the distance between $\sqrt{n}\bd{Z}_n$ and $\f G_n$, say, $d(\sqrt{n}\bd{Z}_n, \bd{G}_n)$ goes to zero under condition $(p, q)=(1,1)$, $(p, q)=(\sqrt{n},1)$, $(p, q)=(o(n),1)$ or $(p, q)=(n^{1/3}, n^{1/3})$. Readers are referred to, for instance, Maxwell \cite{max75, max78}, Poincar\'{e} \cite{poicare}, Stam \cite{stam}, Diaconis {\it et al}. \cite{DLE} and Collins \cite{Collins}. A more detailed recounts can be seen from Diaconis and Freedman \cite{DF87} and Jiang \cite{Jiang06}.

Obviously, with more research being done, it is known that the values of $p$ and $q$ become larger and
larger such that $d(\sqrt{n}\bd{Z}_n, \bd{G}_n)$ goes to zero.
Diaconis \cite{persi03} then asks the largest values of $p_n$ and $q_n$ such that the distance
between $\sqrt{n}\f Z_n$ and $\f G_n$ goes to zero. Jiang \cite{Jiang06}  settles the problem by showing that
$p=o(n^{1/2})$ and $q=o(n^{1/2})$ are the largest orders to make the total variation distance go to
zero. If the distance is the weak distance,
or equivalently, the maximum norm,  Jiang \cite{Jiang06} further proves  that the largest order of
$q$ is $\frac{n}{\log n}$ with $p=n$. Based on this work some applications are obtained, for example,  for
the properties of eigenvalues  of the Jacobi ensemble in the random matrix theory \cite{Jiang09},
the wireless communications \cite{Li2016, Li2014,  Li2016b} and data storage from Big Data \cite{ChenKe}.

However, even with the affirmative answer by Jiang \cite{Jiang06}, a conjecture [(1) below] and
a question [(2) below] still remain.

(1) If $pq/n\to 0$, and $p$ and $q$ do not have to be in the same scale, does the total variation distance
still go to zero?

(2) What if the total variation distance and weak norm are replaced by other popular distances,
say, the Hellinger distance, the Kullback-Leibler distance or the Euclidean distance?

Conjecture (1) is natural because it is shown by Diaconis and Freedman \cite{DF87} that the total variation distance goes to zero if $p=o(n)$ and $q=1$. The work by Jiang \cite{Jiang06} proves that the same holds if $p=o(n^{1/2})$ and $q=o(n^{1/2})$. In both occasions, $(p,q)$ satisfies that $pq=o(n)$.

In this paper we will answer conjecture (1) and question (2). For conjecture (1), we show that the total variation
distance between $\sqrt{n}\bd{Z}_n$ and $\f G_n$ goes to zero as long as
$p\geq 1, q\geq 1$ and $\frac{pq}{n}\to 0$, and the orders are sharp in the sense
that the distance does not go to zero if $\frac{pq}{n}\to \sigma>0$, where $\sigma$ is a constant.

For question (2), we prove that the same answer as that for (1) is also true for the Hellinger
distance and the Kullback-Leibler distance. However, it is different for the Euclidean distance.
We prove that the Euclidean distance between them goes to zero as long as $\frac{pq^2}{n}\to 0$, and the conclusion no longer holds for any $p\geq 1$ and $q\geq 1$ satisfying $\frac{pq^2}{n}\to \sigma>0$. In order to compare
these results clearly, we make Table \ref{table1} for some special cases. One may like to read
the table through its caption and the statements of Theorems
\begin{table}[H]\setlength{\tabcolsep}{10pt}\label{table1}
 \centering
  \begin{threeparttable}[b]
   \begin{tabular}{c|c}
     distance $d$  & order of $(p, q)$ \\
     \midrule
     total variation    & $(\sqrt{n}, \sqrt{n})$    \\
     Hellinger    & $(\sqrt{n}, \sqrt{n})$    \\
     Kullback-Leibler     & $(\sqrt{n}, \sqrt{n})$   \\
     Euclidean    & $(\sqrt[3]{n}, \sqrt[3]{n} )$   \\
     weak     & $(n, \frac{n}{\log n})$   \\
   \end{tabular}
  \end{threeparttable}
 \caption{\small\it Largest orders of $p$ and $q$  such that $d(\sqrt{n}\bd{Z}_n, \bd{G}_n) \to 0$, where $\bd{Z}_n$ and $\bd{G}_n$ are the $p\times q$ upper-left submatrix of an $n\times n$ Haar-invariant orthogonal matrix and a $p\times q$ matrix whose entries are i.i.d. $N(0,1)$, respectively.}
 \end{table}
\noindent
\ref{main1} and \ref{ttt} below.

Before stating our main results, let us review rigorously the distances aforementioned.
Let $\mu$ and $\nu$ be two probability measures on $(\mathbb{R}^m, \ml{B}),$
where $\mathbb{R}^m$ is the $m$-dimensional Euclidean space and $\ml{B}$ is the Borel $\sigma$-algebra.
Recall the total variation distance between $\mu$ and $\nu,$ denoted by $\|\mu-\nu\|_{\rm TV},$
is defined by
\bea\lbl{norm}
\|\mu-\nu\|_{\rm TV}=2\cdot \sup_{A\in \ml{B}}|\mu(A)-\nu(A)|=\int_{\mathbb{R}^m}|f(x)-g(x)|\, dx,
\eea
provided $\mu$ and $\nu$ have density functions $f$ and $g$ with respect to the Lebesgue measure, respectively. The Hellinger distance  $H(\mu, \nu)$ between $\nu$ and $\mu$ is defined by
$$H^2(\mu, \nu)=\frac12\int_{\mathbb{R}^m} |\sqrt{f(x)}-\sqrt{g(x)}\,|^2 dx. $$
The Kullback-Leibler distance between  $\mu$ and $\nu$ is defined by
\beaa\lbl{Kull}
D_{\rm KL}(\mu||\nu)=\int_{\mathbb{R}^m} \frac{d\mu}{d\nu}\log\frac{d\mu}{d\nu} d\nu.
\eeaa
The three distances have the following relationships:
\bea
& & 2 H^2(\mu, \nu)\le \|\mu-\nu\|_{\rm TV} \le 2\sqrt{2}H(\mu, \nu);\lbl{H_TV}\\
& & \|\mu-\nu\|_{\rm TV}^2\le 2 D_{\rm KL}(\mu||\nu).\lbl{TV_KL}
\eea
Readers are referred to, for example,  \cite{Kraft} and \cite{Csiszar} for (\ref{H_TV})
and (\ref{TV_KL}), respectively. In particular, the assertion in (\ref{TV_KL}) is called the Pinsker inequality.

\bthm\label{main1}
Suppose $p=p_n$ and $q=q_n$ satisfy $1\leq p, q\leq n$. For each $n\geq 1$, let $\bd{Z}_n$ and $\bd{G}_n$ be the $p\times q$ submatrices aforementioned.
Let $d(\sqrt{n}\bd{Z}_n, \bd{G}_n)$ be the total variation distance, the Hellinger distance or the Kullback-Leibler  distance
 between the probability distributions of $\sqrt{n}\bd{Z}_n$ and $\bd{G}_n$. Then
\begin{itemize}
\item[(i)] $\lim_{n\to\infty}d(\sqrt{n}\bd{Z}_n, \bd{G}_n)=0$ for any $p\geq 1$ and $q\geq 1$ with $\lim_{n\to\infty}\frac{pq}{n}=0$;

\item[(ii)] $\liminf_{n\to\infty}d(\sqrt{n}\bd{Z}_n, \bd{G}_n) > 0$
if  $\lim_{n\to\infty}\frac{pq}{n}=\sigma\in(0, \infty).$
 \end{itemize}
\nthm
%
When $d(\cdot, \cdot)$ is the total variation distance, Jiang \cite{Jiang06} obtains (i) with $p=o(\sqrt{n})$ and
$q=o(\sqrt{n})$ and (ii) with $p=[x\sqrt{n}\,]$ and $q=[y\sqrt{n}\,]$ where $x>0$ and $y>0$ are constants. Theorem \ref{main1} confirms a conjecture by the first author.

Now we study the approximation in terms of the Euclidean distance. Let $\bd{Y}_n=(\bd{y}_1, \cdots, \bd{y}_n)=(y_{ij})_{n\times n}$ be an $n\times n$ matrix, where $y_{ij}$'s are i.i.d. random variables with distribution $N(0,1)$. Perform the Gram-Schmidt algorithm on the column vectors $\bd{y}_1, \cdots, \bd{y}_n$ as follows.
\bea
&& \bd{w}_1=\bd{y}_1,\ \ \bm{\gamma}_1=\frac{\bd{w}_1}{\|\bd{w}_1\|}; \nonumber\\
&& \bd{w}_k=\bd{y}_k-\sum_{i=1}^{k-1}\langle\bd{y}_k, \bm{\gamma}_i\rangle \bm{\gamma}_i,\ \ \bm{\gamma}_k=\frac{\bd{w}_k}{\|\bd{w}_k\|} \lbl{sea}
\eea
for $k=2,\cdots, n$, where $\langle\bd{y}_k, \bm{\gamma}_i\rangle$ is the inner product of the two vectors.
Then $\bm{\Gamma}_n=(\bm{\gamma}_1, \cdots, \bm{\gamma}_n)=(\gamma_{ij})$ is an $n\times n$
Haar-invariant orthogonal matrix. Set $\bd{\Gamma}_{p\times q}=(\gamma_{ij})_{1\leq i\leq p, 1\leq j \leq q}$
and $\bd{Y}_{p\times q}=(y_{ij})_{1\leq i\leq p, 1\leq j \leq q}$ for $1\leq p, q\leq n.$
We consider the Euclidean distance between $\sqrt{n}\bd{\Gamma}_{p\times q}$ and $\bd{Y}_{p\times q}$, that is,
the Hilbert-Schmidt norm defined by
\bea\lbl{Hiha}
\|\sqrt{n}\bd{\Gamma}_{p\times q}-\bd{Y}_{p\times q}\|^2_{\rm HS}=\sum_{i=1}^p\sum_{j=1}^q(\sqrt{n}\gamma_{ij}-y_{ij})^2.
\eea

Throughout the paper the notation $\xi_n\overset{p}{\to} \xi$ indicates that random variable $\xi_n\to \xi$ in probability as $n\to\infty$.

\bthm\label{ttt}
Let the notation $\bd{\Gamma}_{p\times q}$ and $\bd{Y}_{p\times q}$ be as in the above. If $p=p_n, q=q_n$ satisfy $1\leq p, q\leq n$ and $\lim_{n\to\infty}\frac{pq^2}{n}= 0$, then
$\|\sqrt{n}\bd{\Gamma}_{p\times q}-\bd{Y}_{p\times q}\|_{\rm HS} \overset{p}{\to} 0$ as $n\to\infty$.
Further, if $1\leq p, q\leq n$ satisfy $\lim_{n\to\infty}\frac{pq^2}{n}= \sigma \in (0, \infty)$, then
\bea\lbl{indication_phone}
\liminf_{n\to\infty}P(\|\sqrt{n}\bd{\Gamma}_{p\times q}-\bd{Y}_{p\times q}\|_{\rm HS}\geq \epsilon)>0
\eea
for every $\epsilon \in (0, \sqrt{\sigma/2})$.
\nthm
We also obtain an upper bound in Proposition \ref{lala}:  $\mathbb{E}
\|\sqrt{n}\bd{\Gamma}_{p\times q}-\bd{Y}_{p\times q}\|^2_{\rm HS} \leq \frac{24pq^2}{n}$
for any $n\ge 2$ and $1\leq p, q\leq n$.  Further, we obtain cleaner results than (\ref{indication_phone}) for two special cases. It is proved in Lemma \ref{piano_black} that
$\|\sqrt{n}\bd{\Gamma}_{p\times 1}-\bd{Y}_{p\times 1}\|_{\rm HS} \to \sqrt{\frac{c}{2}}\cdot |N(0,1)|$ weakly provided  $p/n\to c\in (0,1]$. In the proof of Theorem \ref{ttt}, we show that
$\|\sqrt{n}\bd{\Gamma}_{p\times q}-\bd{Y}_{p\times q}\|_{\rm HS} \overset{p}{\to} \sqrt{\sigma/2}$ if $q\to \infty$ and $(pq^2)/n\to  \sigma>0$.

In order to compare the orders for all different norms, we make Table \ref{table1} for
the special case that $p$ and $q$ are of the same scale except for the weak norm.
The weak norm is defined by $\|\bd{A}-\bd{B}\|_{\rm max}=\max_{1\leq i \leq p, 1\leq j \leq q}
|a_{ij}-b_{ij}|$ for $\bd{A}=(a_{ij})_{p\times q}$ and $\bd{B}=(b_{ij})_{p\times q}$.
The distance $\|\sqrt{n}\bd{Z}_n-\bd{G}_n\|_{\rm max}$ for the case $p=n$
is studied in \cite{Jiang06}.

\medskip

\noindent\textbf{Remarks and future questions} 

{\bf A}. Compared to the techniques employed in \cite{Jiang06}, the proofs of the results
in this paper use the following new elements:

\begin{enumerate}
\item  Tricks of calculating the means of monomials
of the entries from $\bd{\Gamma}_n$ are used in Lemmas \ref{what} and \ref{great}.

\item A subsequence argument is applied to the proofs of both theorems. In particular,  the proof of Theorem \ref{main1} is reduced to the case $q/p\to 0$ and the case  $q\equiv 1$.

\item  A central limit theorem (CLT) on $\mbox{tr}[(\bd{G}_n'\bd{G}_n)^2]$  for the case
$q/p\to 0$  is established in Lemma \ref{awesome_pen}. The CLT for the case $q/p\to c>0$
is well known; see, for example, \cite{Bai_Jack} or \cite{Jonsson}.

\item  Some properties of the largest and
the smallest eigenvalues of $\bd{G}_n'\bd{G}_n$ for the case $q/p\to 0$ is proved in
 Lemma \ref{fly_ground}. This is a direct consequence of a recent result by Jiang
and Li \cite{JiangLi15}. The situation for $q/p\to c>0$ is well known; see,
for example, \cite{Bai_Jack}.

\item  Connections in (\ref{H_TV}) and (\ref{TV_KL}) among distances
provide an efficient way to use known properties of Wishart matrices and Haar-invariant orthogonal matrices. The Wishart matrices appear in ``Proof of (ii) of Theorem \ref{main1}" and the Haar-invariant orthogonal matrices occur in ``Proof of (i) of Theorem \ref{main1}."
\end{enumerate}

{\bf B}. In this paper we approximate the Haar-invariant orthogonal matrices by independent normals with
various probability measures. It can be proved that similar results also hold for Haar-invariant unitary and symplectic matrices without difficulty.
This can be done by the method employed here together with those from \cite{Jiang09, Jiang10}.

{\bf C}. As mentioned earlier, the work \cite{Jiang06} has been applied to other random matrix problems \cite{Jiang09},
the wireless communications \cite{Li2016, Li2014,  Li2016b} and a problem from Big Data \cite{ChenKe}. In this paper
we consider other three probability metrics:  the Hellinger distance, the Kullback-Leibler distance
and the Euclidean distance. We expect more applications. In particular, since
Hellinger distance and Kullback-Leibler distance are popular in Statistics and Information Theory, respectively,
we foresee some applications in the two areas.

{\bf D}. In Theorem \ref{ttt}, the Haar-invariant orthogonal matrices are obtained
by the Gram-Schmidt algorithm. The approximation by independent normals via the Hilbert-Schmidt norm is valid if $pq^2=o(n)$.
There are other ways to generate Haar-invariant orthogonal matrices; see, for example, \cite{Mezzadri}.
It will be interesting to see the cut-off orders of $p$ and $q$ such that (\ref{indication_phone}) holds under the new couplings.

{\bf E}. So far five popular probability metrics are applied to study the distance
between $\sqrt{n}\bd{Z}_n$ and
$\bd{G}_n$. They are the total variation distance, the Hellinger distance,
the Kullback-Leibler distance,
the Euclidean distance in this paper and the weak norm in \cite{Jiang06}. Their corresponding conclusions show different features.  There are many other distances of
probability measures which include the Prohorov distance, the Wasserstein distance and the Kantorovich
transport distance. It will be interesting to see the largest orders of $p$ and $q$ such that these distances go to zero.
Of course, applications of the results along this line are welcomed.

\medskip

Finally, the structure of the rest paper is organized as follows.

\noindent\textbf{Section \ref{sky_phone}: Proof of Theorem \ref{main1}}

Section \ref{sun_leave}: Preliminary Results.

Section \ref{Proof_Main1}: The Proof of Theorem \ref{main1}.

\noindent\textbf{Section \ref{yes_ttt}: Proof of Theorem \ref{ttt}}

Section \ref{papapa}: Auxiliary Results.

Section \ref{big_mouse}: The Proof of Theorem \ref{ttt}.

\noindent \textbf{Section \ref{Good_Tue}: Appendix}.

\section{Proof of Theorem \ref{main1}}\lbl{sky_phone}

\subsection{Preliminary Results}\lbl{sun_leave} Throughout the paper we
will adopt the following notation.

\noindent\textbf{Notation}. (a) $X\sim \chi^2(k)$ means that random variable $X$ follows the
chi-square distribution with degree of freedom $k$;

(b) $N_p(\bm{\mu}, \bd{\Sigma})$ stands for the $p$-dimensional normal distribution of mean vector $\bm{\mu}$ and covariance matrix $\bd{\Sigma}.$ We write $\bd{X} \sim N_p(\bm{\mu}, \bd{\Sigma})$ if random vector $\bd{X}$ has the distribution  $N_p(\bm{\mu}, \bd{\Sigma})$. In particular, we write $\bd{X} \sim N_p(\bm{0}, \bd{I})$ if  the $p$ coordinates of $\bd{X}$ are independent $N(0, 1)$-distributed random variables.

(c) For two sequences of numbers $\{a_n;\, n\geq 1\}$ and $\{b_n;\, n\geq 1\}$, the notation $a_n=O(b_n)$ as $n\to\infty$ means that $\limsup_{n\to\infty}|a_n/b_n|<\infty.$ The notation $a_n=o(b_n)$ as $n\to\infty$ means that $\lim_{n\to\infty}a_n/b_n=0$, and the symbol $a_n \sim b_n$ stands for $\lim_{n\to\infty}a_n/b_n=1$.

(d) $X_n=o_p(a_n)$ means $\frac{X_n}{a_n}\to 0$ in probability as $n\to\infty$.
The symbol $X_n=O_p(a_n)$ means that $\{\frac{X_n}{a_n};\, n\geq 1\}$
 are stochastically bounded, that is,
$\sup_{n\geq 1}P(|X_n|\geq b a_n)\to 0$ as $b\to \infty.$

Before proving Theorem \ref{main1}, we need some preliminary results. They appear in a series of lemmas.

The following is taken from Proposition 2.1 by Diaconis, Eaton and Lauritzen \cite{DLE} or Proposition 7.3 by Eaton \cite{ME2}.

\blem\lbl{del} Let $\f \Gamma_n$ be an $n\times n$ random matrix which is uniformly
distributed on the orthogonal group $O(n)$ and let $\f Z_n$ be the upper-left $ p\times q$  submatrix   of $\f \Gamma_n.$ If $p+q\leq n$ and $q\leq p$ then the joint density function of entries of $\f Z_n$ is
\bea\lbl{heavy}
f(z)=(\sqrt{2\pi})^{-pq}\frac{\omega(n-p, q)}{\omega(n, q)}\left\{\mbox{det}(I_{q}-z'z)^{(n-p-q-1)/2}\right\}I_0(z'z)
\eea
where $I_0(z'z)$ is the indicator function of the set that all $q$ eigenvalues of $z'z$ are in $(0,1),$  and $\omega(\cdot, \cdot)$ is the Wishart constant defined by
\bea\lbl{fire_you}
\frac{1}{\omega(s, t)}=\pi^{t(t-1)/4}2^{st/2}\prod_{j=1}^t\Gamma\left(\frac{s-j+1}{2}\right).
\eea
Here $t$ is a positive integer and $s$ is a real number, $s>t-1.$ When $p< q,$ the density of $\f Z_n$ is obtained by interchanging $p$ and $q$ in the above Wishart constant.
\nlem

The following result is taken from \cite{Jiang2009}. For any integer $a\geq 1$, set $(2a-1)!!=1\cdot 3\cdots (2a-1)$ and $(-1)!!=1$ by convention.
\blem\lbl{Jiang2009} Suppose $m\geq 2$ and $\xi_1, \cdots, \xi_m$ are i.i.d. random variables with $\xi_1 \sim N(0,1).$ Define $U_i=\frac{\xi_i^2}{\xi_1^2 + \cdots + \xi_m^2}$ for $1\leq i \leq m$. Let $a_1, \cdots, a_m$ be non-negative integers and $a=\sum_{i=1}^m a_i$. Then
\beaa
E\big(U_1^{a_1}\cdots U_m^{a_m}\big) = \frac{\prod_{i=1}^m(2a_i-1)!!}{\prod_{i=1}^a(m+2i-2)}.
\eeaa
\nlem

\medskip

The expectations of some monomials of the entries of Haar-orthogonal matrices will be  computed next. Recall $\bold{\Gamma}_n=(\gamma_{ij})_{n\times n}$ is an Haar-invariant orthogonal matrix. The  following facts will be repeatedly used later. They follow from the property of the Haar invariance.
\begin{itemize}
\item[{\bf F1})]\, The vector $(\gamma_{11},\cdots, \gamma_{n1})'$ and $\frac{1}{\sqrt{\xi_1^2+\cdots +\xi_n^2}}(\xi_1, \cdots, \xi_n)'$ have the same probability distribution, where $\xi_1, \cdots, \xi_n$ are i.i.d. $N(0,1)$-distributed random variables.
\item[{\bf F2})]\, By the orthogonal invariance, for any $1\leq k \leq n$, any $k$  different rows/columns of $\bold{\Gamma}_n$ have the same joint distribution as that of the first $k$  rows of $\bold{\Gamma}_n$.
\end{itemize}
\blem\lbl{what} Let $\bd{\Gamma}_n=(\bm{\gamma}_1, \cdots, \bm{\gamma}_n)=(\gamma_{ij})_{n\times n}$ be an Haar-invariant orthogonal matrix. Then

(a)\, $\mathbb{E}(\gamma_{11}^2)=\frac{1}{n}$ and $\mathbb{E}(\gamma_{11}^4)=\frac{3}{n(n+2)}$;

(b)\, $\mathbb{E}(\gamma_{11}^2\gamma_{12}^2)=\frac{1}{n(n+2)}\ \ \mbox{and} \ \ \mathbb{E}(\gamma^2_{11}\gamma^2_{22})=\frac{n+1}{n(n-1)(n+2)}$;

(c)\, $\mathbb{E}(\gamma_{11}\gamma_{12}\gamma_{21}\gamma_{22})=-\frac{1}{n(n-1)(n+2)}$.
\nlem
\noindent\textbf{Proof}.
By Property {\bf F1}), picking $m=n$, $a_1=1, a_2=\cdots =a_n=0$ from Lemma \ref{Jiang2009}, we see $\mathbb{E}(\gamma_{11}^2)=\frac{1}{n}$. Choosing $a_1=2, a_2=\cdots =a_n=0$, we obtain $\mathbb{E}(\gamma_{11}^4)=\frac{3}{n(n+2)}$. Selecting $a_1=a_2=1, a_3=\cdots =a_n=0$, we see $\mathbb{E}(\gamma_{11}^2\gamma_{12}^2)=\frac{1}{n(n+2)}$.

Now, since $\|\bm{\gamma}_1\|=\|\bm{\gamma}_2\|=1$, by {\bf F2})
\beaa
1=\mathbb{E}\big(\|\bm{\gamma}_1\|^2\|\bm{\gamma}_2\|^2\big)
&=&\mathbb{E}\bigg(\sum_{i=1}^n \gamma_{i1}^2\gamma_{i2}^2+\sum_{1\le i\neq j \le n}\gamma_{i1}^2\gamma_{j2}^2\bigg)\\
& = & n\mathbb{E}(\gamma_{11}^2\gamma_{12}^2) + n(n-1)\mathbb{E}(\gamma_{11}^2 \gamma_{22}^2)\\
&=&\frac{1}{n+2}+n(n-1)\mathbb{E}(\gamma_{11}^2 \gamma_{22}^2).
\eeaa
The second conclusion of (b) is yielded.
Now we work on conclusion (c). In fact, since the first two columns of $\f {\Gamma}_n$ are orthogonal, we know
\bea\lbl{what_green}
0=\big(\sum_{i=1}^n\gamma_{i1}\gamma_{i2}\big)^2=\sum_{i=1}^n\gamma_{i1}^2\gamma_{i2}^2+\sum_{1\leq i\ne j \leq n}\gamma_{i1}\gamma_{i2}\gamma_{j1}\gamma_{j2}.
\eea
By Property {\bf F2)} again,
$$\mathbb{E}(\gamma_{i1}\gamma_{i2}\gamma_{j1}\gamma_{j2})=\mathbb{E}(\gamma_{11}\gamma_{12}\gamma_{21}\gamma_{22})$$ for any $i \ne j$. Hence, take expectations of both sides of (\ref{what_green}) to see
\beaa\lbl{sun_land}
\mathbb{E}(\gamma_{11}\gamma_{12}\gamma_{21}\gamma_{22})=-\frac{n}{n(n-1)}\mathbb{E}(\gamma_{11}^2\gamma_{12}^2)=-\frac{1}{(n-1)n(n+2)}.
\eeaa
\hfill$\blacksquare$

In order to understand the trace of the third power of an Haar-invariant orthogonal matrix, we need the following expectations of monomials of the matrix elements.
\blem\lbl{great} Let $\f \Gamma_n=(\gamma_{ij})_{n\times n}$ be a  random matrix with the uniform
distribution on the orthogonal group $O(n)$, $n\geq 3.$ The following holds:

(a) $\mathbb{E}(\gamma_{11}^2\gamma_{21}^2\gamma_{31}^2)=\frac{1}{n(n+2)(n+4)}$.

(b) $\mathbb{E}(\gamma_{11}\gamma_{12}\gamma_{21}\gamma_{22}\gamma_{23}^2)=
-\frac{1}{(n-1)n(n+2)(n+4)}.$

(c) $\mathbb{E}(\gamma_{11}^2\gamma_{21}^2\gamma_{22}^2)=
\frac{1}{(n-1)n(n+2)}-\frac{3}{(n-1)n(n+2)(n+4)}.$

(d) $\mathbb{E}(\gamma_{11}\gamma_{12}\gamma_{21}\gamma_{22}^3)=-\frac{3}{(n-1)n(n+2)(n+4)}$.

(e) $\mathbb{E}(\gamma_{11}\gamma_{12}\gamma_{22}\gamma_{23}\gamma_{31}\gamma_{33})
=\frac{2}{(n-2)(n-1)n(n+2)(n+4)}.$
\nlem
Obviously, Lemma \ref{great} is more complex than Lemma \ref{what}.
We postpone its proof in Appendix from Section \ref{Good_Tue}.

Based on Lemma \ref{what}, we now present two identities that will be used later.

\blem\label{KL-Key}  Let $\lambda_1, \cdots, \lambda_q$ be the eigenvalues of $\f {Z}_n'\f {Z}_n$, where $\f Z_n$ is defined as in Lemma \ref{del}. Then
$$\aligned &\mathbb{E}\sum_{i=1}^q\lambda_i=\frac{pq}{n}; \\
 &\mathbb{E}\sum_{i=1}^q\lambda_i^2=\frac{pq}{n(n+2)}\Big[p+q+1-\frac{(p-1)(q-1)}{n-1}\Big].
\endaligned $$
\nlem
\noindent\textbf{Proof}.  
The first equality is trivial since
\beaa\label{mean}
\mathbb{E}\sum_{i=1}^q \lambda_i=\mathbb{E}\mbox{tr}(\f Z_n'\f Z_n)=\mathbb{E}\sum_{i=1}^{p}\sum_{j=1}^q \gamma_{ij}^2=\frac{pq}{n}
\eeaa
since $E(\gamma_{ij}^2)=E(\gamma_{11}^2)=\frac{1}{n}$ for any $i,j$ by (a) of Lemma \ref{what}. For the second equality, first
\bea
\sum_{i=1}^q\lambda_i^2=\mbox{tr}(\f {Z}_n'\f{Z}_n\f{Z}_n'\f{Z}_n)
&=& \sum_{1\leq j,l\leq p; 1\leq i, k\leq q}\gamma_{ji}\gamma_{jk}\gamma_{lk}\gamma_{li}\nonumber\\
& =: & A + B +C,\lbl{additional_care}
\eea
where $A$ corresponds to that $j=l,\ i=k$; $B$ corresponds to that $j=l,\, i \ne k$ or $j \ne l,\, i=k$; $C$ corresponds to that $j \ne l,\, i \ne k.$ It is then easy to see that
\beaa
& & A=\sum_{1\leq j \leq p,\, 1\leq i \leq q}\gamma_{ji}^4;\ \ \ B=\sum_{1\leq j \leq p,\, 1\leq i \ne k\leq q}\gamma_{ji}^2\gamma_{jk}^2 + \sum_{1\leq j\ne l \leq p,\, 1\leq i \leq q}\gamma_{ji}^2\gamma_{li}^2;\\
& & C=\sum_{1\leq j\ne l\leq p;\, 1\leq i\ne k\leq q}\gamma_{ji}\gamma_{jk}\gamma_{lk}\gamma_{li}.
\eeaa
\medskip
By Properties {\bf F1}) and {\bf F2}) and Lemma \ref{what}, we see
\beaa
\mathbb{E}A&=&pq \cdot E(\gamma_{11}^4)=
\frac{3pq}{n(n+2)};\lbl{thousand}\\
\mathbb{E}B &=& [pq(q-1)+pq(p-1)]\cdot \mathbb{E}(\gamma_{11}^2\gamma_{12}^2)
=\frac{pq(p+q-2)}{n(n+2)};\nonumber\\
\mathbb{E}C &= & pq(p-1)(q-1)\cdot \mathbb{E}(\gamma_{11}\gamma_{12}\gamma_{21}\gamma_{22})=-\frac{pq(p-1)(q-1)}{(n-1)n(n+2)}.\nonumber
\eeaa
%
Consequently,
\beaa
\mathbb{E}\sum_{i=1}^q\lambda_i^2
&=&\mathbb{E}A+\mathbb{E}B + \mathbb{E}C\\
&=& \frac{3pq}{n(n+2)} + \frac{pq(p+q-2)}{n(n+2)}-\frac{pq(p-1)(q-1)}{(n-1)n(n+2)}\\
& = & \frac{pq}{n(n+2)}\Big[p+q+1-\frac{(p-1)(q-1)}{n-1}\Big].
\eeaa
The proof is completed. \hfill$\blacksquare$

\medskip

With Lemma \ref{great}, we are ready to compute the following quantity.

 \blem\label{tri}  Let $\lambda_1, \cdots, \lambda_q$ be the eigenvalues of $\f {Z}_n'\f {Z}_n.$ Then,
\beaa
\mathbb{E}\sum_{i=1}^q \lambda_i^3&=&\frac{pq}{n(n+2)(n+4)}\big[p^2+q^2+3pq+3(p+q)+4\big]\\
&&\quad\quad +\frac{pq(p-1)(q-1)}{(n-1)n(n+2)(n+4)}\bigg[\frac{2(p-2)(q-2)}{n-2}-3(p+q)\bigg].
\eeaa
\nlem
\noindent\textbf{Proof}. By definition,
\bea\lbl{sumsum}
\sum_{i=1}^q\lambda_i^3=\mbox{tr}(\f {Z}_n'\f{Z}_n\f{Z}_n'\f{Z}_n\f{Z}_n'\f{Z}_n)
&= & \sum_{1\leq i, j,k\leq p; 1\leq l, s, t\leq q}\gamma_{il}\gamma_{is}\gamma_{js}\gamma_{jt}\gamma_{kt}\gamma_{kl}\nonumber\\
& = & A_1+A_2 + A_3,
\eea
where $A_1$ corresponds to the sum over $i=j=k$, $A_2$ corresponds to the sum that only two of $\{i, j, k\}$ are identical, and $A_3$ corresponds to the sum $i\ne j\ne k$. We next compute each term in detail.

{\it Case 1: $i=j=k$}. Each term in the sum has the expression $\mathbb{E}(\gamma_{il}^2\gamma_{is}^2\gamma_{it}^2)$. The corresponding sum then becomes
\beaa
A_1&=& \sum_{i=1}^p\sum_{1\leq l,s,t\leq q}\mathbb{E}(\gamma_{il}^2\gamma_{is}^2\gamma_{it}^2)\\
& = & pq\mathbb{E}(\gamma_{11}^6)+ 3pq(q-1)\mathbb{E}(\gamma_{11}^4\gamma_{12}^2)+  pq(q-1)(q-2)\mathbb{E}(\gamma_{11}^2\gamma_{12}^2\gamma_{13}^2)\\
& = & \frac{15pq}{n(n+2)(n+4)} + \frac{9pq(q-1)}{n(n+2)(n+4)} + \frac{pq(q-1)(q-2)}{n(n+2)(n+4)}
\eeaa
by {\bf F2)}, Lemmas \ref{Jiang2009}, \ref{what} and  \ref{great}.

{\it Case 2: only two of $\{i, j, k\}$ are identical.} The corresponding sum is
\beaa
A_2 & =& 3\sum_{1\leq i\ne k\leq p}\sum_{1\leq l,s,t\leq q}\mathbb{E}(\gamma_{il}\gamma_{is}^2\gamma_{it}\gamma_{kt}\gamma_{kl})\\
&=& 3p(p-1) \sum_{1\leq l,s,t\leq q}\mathbb{E}(\gamma_{1l}\gamma_{1s}^2\gamma_{1t}\gamma_{2t}\gamma_{2l}).
\eeaa
 By symmetry and {\bf F2)}
\beaa
\frac{A_2}{3p(p-1)}&=&\sum_{l=1}^q\mathbb{E}(\gamma_{1l}^4\gamma_{2l}^2) + \sum_{1\leq l\ne s\leq q}\mathbb{E}(\gamma_{1l}\gamma_{1s}^3\gamma_{2s}\gamma_{2l})+ \sum_{1\leq l\ne s\leq q}\mathbb{E}(\gamma_{1l}^2\gamma_{1s}^2\gamma_{2l}^2)\\
 &&\quad \quad\quad\quad + \sum_{1\leq l \ne s\leq q}\mathbb{E}(\gamma_{1l}^3\gamma_{1t}\gamma_{2t}\gamma_{2l})
 + \sum_{1\leq l\ne s\ne t\leq q}\mathbb{E}(\gamma_{1l}\gamma_{1s}^2\gamma_{1t}\gamma_{2t}\gamma_{2l})\\
 & = & q\cdot\mathbb{E}(\gamma_{11}^4\gamma_{21}^2) + q(q-1)\cdot\mathbb{E}(\gamma_{11}\gamma_{12}\gamma_{21}\gamma_{22}^3) + q(q-1)\cdot\mathbb{E}(\gamma_{11}^2\gamma_{21}^2\gamma_{22}^2)\\
& &   + q(q-1)\cdot\mathbb{E}(\gamma_{11}\gamma_{12}\gamma_{21}\gamma_{22}^3) + q(q-1)(q-2)\cdot\mathbb{E}(\gamma_{11}\gamma_{12}\gamma_{21}\gamma_{22}\gamma_{23}^2)
\eeaa
where the sums in the first equality appearing in order correspond to $l=s=t$,  $l\ne s=t$, $s\ne l= t$, $t\ne l= s$ and  $l \ne s \ne t$,  respectively. By Lemmas \ref{Jiang2009} and \ref{great},
\beaa
& & A_2\\
 & =& 3p(p-1)\Big[\frac{3q}{n(n+2)(n+4)} - \frac{6q(q-1)}{(n-1)n(n+2)(n+4)}\\
& & +q(q-1)\Big(\frac{1}{(n-1)n(n+2)}-\frac{3}{(n-1)n(n+2)(n+4)}\Big)-
\frac{q(q-1)(q-2)}{(n-1)n(n+2)(n+4)}\Big)\Big].
\eeaa

{\it Case 3:  $i\ne j\ne k$.} The corresponding sum becomes
\beaa
A_3=p(p-1)(p-2)\sum_{1\leq l, s, t\leq q}\mathbb{E}\gamma_{1l}\gamma_{1s}\gamma_{2s}\gamma_{2t}\gamma_{3l}\gamma_{3t}.
\eeaa
By symmetry and the same classification as that in {\it Case 2},
\beaa
& & \frac{A_3}{p(p-1)(p-2)}\\
& = & q\cdot \mathbb{E}(\gamma_{11}^2\gamma_{21}^2\gamma_{31}^2) + 3q(q-1)\cdot \mathbb{E}(\gamma_{11}\gamma_{12}\gamma_{21}\gamma_{22}\gamma_{23}^2)\\
& & \quad\quad\quad\quad\quad\quad +\,  q(q-1)(q-2)\cdot \mathbb{E}(\gamma_{11}\gamma_{12}\gamma_{22}\gamma_{23}\gamma_{31}\gamma_{33})\\
&= & \frac{q}{n(n+2)(n+4)}-\frac{3q(q-1)}{(n-1)n(n+2)(n+4)} + \frac{2q(q-1)(q-2)}{(n-2)(n-1)n(n+2)(n+4)}.
\eeaa
Combing (\ref{sumsum}) and the formulas on $A_1, A_2$ and $A_3$, we see
\beaa
\mathbb{E}\sum_{i=1}^q \lambda_i^3&=&\frac{pq}{n(n+2)(n+4)}\bigg[p^2+q^2+6(p+q)+1\bigg]+\frac{3pq(p-1)(q-1)}{(n-1)n(n+2)}\\
&&-\frac{pq(p-1)(q-1)}{(n-1)n(n+2)(n+4)}\bigg[15+3(p+q)-\frac{2(p-2)(q-2)}{n-2}\bigg].
\eeaa
Now write
\beaa
\frac{3pq(p-1)(q-1)}{(n-1)n(n+2)}=\frac{3pq(p-1)(q-1)}{n(n+2)(n+4)}+
\frac{15pq(p-1)(q-1)}{(n-1)n(n+2)(n+4)}.
\eeaa
By making a substitution, we obtain the desired formula.
\hfill$\blacksquare$

\medskip

The normalizing constant from (\ref{heavy}) needs to be understood. It is given below.
\blem\lbl{KO} For $1\leq q\le p<n$, define
\be\lbl{eq1} K_n:=\left(\frac{2}{n}\right)^{pq/2}\prod_{j=0}^{q-1}\frac{\Gamma((n-j)/2)}{\Gamma((n-p-j)/2)}.\ee  If $p=p_n\to \infty$, $\limsup_{n\to\infty}\frac{p}{n}<1$ and  $pq=O(n)$, then
\be\lbl{KmO}\log K_n=-\frac{pq}{2}+\frac{q(q+1)}{4}\log \big(1+\frac p{n-p}\big)-\frac{pq^3}{12n^2}-c_n q\log\big(1-\frac pn\big)+o(1)\ee
as $n\to\infty$, where $c_n:=\frac{1}{2}(n-p-q-1).$
\nlem
\noindent\textbf{Proof}. Recalling the Stirling formula (see, e.g., p. 204 from \cite{Ahlfors} or p. 368 from \cite{Gamelin}),
$$\log\Gamma(x)=\big(x-\frac{1}{2}\big)\log x-x+\log\sqrt{2\pi}+\frac1{12x}+O(\frac{1}{x^3})$$
as $x\to+\infty.$
 Then, we have from the fact $q=o(n)$ that
\beaa\lbl{Kmgeneral}
\log K_n &=& -\frac{pq}{2}\log \frac{n}{2}+\sum_{j=0}^{q-1}\Big[\log\Gamma\big(\frac{n-j}{2}\big)-
\log\Gamma\big(\frac{n-p-j}{2}\big)\Big]\\
& = & -\frac{pq}{2}\log \frac{n}{2}+\sum_{j=0}^{q-1}\Big[\frac{n-j-1}{2}\log \frac{n-j}{2}-\frac{n-p-j-1}{2}\log \frac{n-p-j}{2}
-\frac{p}{2}\Big]\\
& & \quad \quad \quad \quad \quad \quad\quad \quad \quad\quad \quad \quad \quad \quad \quad\quad \quad \quad \quad \quad \quad \quad \quad \quad\quad +o(1).
\eeaa
Now, writing
$\frac{n-j-1}{2}=\frac{n-p-j-1}{2}+\frac{p}{2}$ and putting term ``$-\frac{pq}{2}\log \frac{n}{2}$" into ``$\sum_{j=0}^{q-1}$", we see
\bea\lbl{Knmiddle}
\log K_n
= -\frac{pq}{2}+\sum_{j=0}^{q-1}\frac{n-p-j-1}{2}\log\frac{n-j}{n-p-j}+\frac p2\sum_{j=0}^{q-1}\log\frac{n-j}{n}+o(1).
\eea
It is easy to check
\bea\lbl{substi}
\log\frac{n-j}{n-p-j}=-\log\big(1-\frac{p}{n}\big)+\log\Big[1+\frac{pj}{n(n-p-j)}\Big].
 \eea
Putting \eqref{substi} back into the expression \eqref{Knmiddle}, we have
\beaa\lbl{semi-final}
\log K_n &=& -\frac{pq}{2}+\frac p2\sum_{j=0}^{q-1}\log(1-\frac{j}{n})\\
&&+\sum_{j=0}^{q-1}\frac{n-p-j-1}{2}\Big[-\log(1-\frac{p}{n})+
\log\Big(1+\frac{pj}{n(n-p-j)}\Big)\Big]+o(1)\\
&=& -\frac{pq}{2}-\log(1-\frac{p}{n})\sum_{j=0}^{q-1}\frac{n-p-j-1}{2}\\
&& +\sum_{j=0}^{q-1}\Big[\frac{p}{2}\log\big(1-\frac{j}{n}\big)+
\frac{n-p-j-1}{2}\log\Big(1+\frac{pj}{n(n-p-j)}\Big)\Big]+o(1)\\
&=& -\frac {pq} 2-\bigg[\frac{(n-p)q}{2}-\frac{q(q+1)}{4}\bigg]\log(1-\frac{p}{n})\\
&& +\sum_{j=0}^{q-1}\Big[\frac{p}{2}\log\Big(1-\frac{j}{n}\big)+
\frac{n-p-j-1}{2}\log\Big[1+\frac{pj}{n(n-p-j)}\Big]+o(1).
\eeaa
Since $\log (1+x)=x-\frac{1}{2}x^2+o(x^3)$ as $x\to 0$, we have
\beaa &&\frac{p}{2}\log\big(1-\frac{j}{n}\big)+\frac{n-p-j-1}{2}\log\Big[1+\frac{pj}{n(n-p-j)}\Big]\\
&=&-\frac{pj}{2n}-\frac{pj^2}{4n^2}+
\frac{n-p-j-1}{2}\cdot\frac{pj}{n(n-p-j)}
+O\Big(\frac{pq^3}{n^3}+ \frac{1}{n}\Big)\\
&=&-\frac{pj}{2n}-\frac{pj^2}{4n^2}+\frac{pj}{2n}-\frac{pj}{2n(n-p-j)}+
O\Big(\frac{1}{n}\Big)\\
&=&-\frac{pj^2}{4n^2}+O\Big(\frac{1}{n}\Big),
\eeaa
uniformly for all $1\leq j \leq q$, where we use the fact $\max_{1\leq j \leq q}\frac{j}{n}=\frac{q}{n}$, $\max_{1\leq j \leq q}\frac{pj}{n(n-p-j)}\leq \frac{pq}{n(n-p-q)}=O(\frac{1}{n})$ and $\frac{pq^3}{n^3}=O(\frac{1}{n})$ by the condition $p\to \infty$, $q\leq p$ and $pq=O(n)$ in the calculation. Combining the last two assertions, we conclude
$$\aligned &\log K_n=-\frac {pq} 2-\bigg[\frac{(n-p)q}{2}-\frac{q(q+1)}{4}\bigg]\log\Big(1-\frac{p}{n}\Big)-
\frac{pq^3}{12n^2}+O\Big(\frac{q}{n}\Big)\\
&=-\frac{pq}{2}-\frac{q(q+1)}{4}\log\Big(1-\frac pn\Big)-\frac{pq^3}{12n^2}-c_n q\log\Big(1-\frac pn\Big)+o(1)
\endaligned$$
with $c_n=\frac{1}{2}(n-p-q-1).$
\hfill$\blacksquare$
\medskip

Now we present some properties of the chi-square distribution.

\blem\lbl{party_love} Given integer $m\ge 1,$ review the random variable $\chi^2_{m}$ has density function $$
f(x)=\frac{x^{\frac m2-1} e^{-\frac x2}}{2^{\frac{m}{2}}\Gamma(\frac m2)}$$ for any $x>0$.  Then
$$ \mathbb{E}(\chi^2_{m})^k=\prod_{l=0}^{k-1}(m+2l)
$$
for any positive integer $k.$
In particular,  we have
\beaa
& & {\rm Var}(\chi^2_{m})=2m, \quad  {\rm Var}\big((\chi^2_m-m)^2\big)=8m(m+6); \\
& & {\rm Var}\big((\chi^2_m)^2\big)=8m(m+2)(m+3); \\
& & \mathbb{E}((\chi^2_m-m)^3)=8m, \quad \mathbb{E}((\chi^2_m-m)^4)=12m(m+4).
\eeaa
\nlem

\noindent {\bf Proof.}
Note that
\be\lbl{chik}\aligned \mathbb{E}(\chi^2_{m})^k&=\frac1{2^{\frac{m}{2}}\Gamma(\frac m2)}\int_{0}^{\infty} x^{\frac {m+2k}2-1} e^{-\frac x2} dx\\
&=\frac{2^k\Gamma(\frac m2+k)}{\Gamma(\frac m2)}=\prod_{i=0}^{k-1}(m+2i)
\endaligned
\ee
for any $k\ge 1$.
Here for the last equality we use the property of the Gamma function that $\Gamma(l+1)=l\Gamma(l)$ for any $l>0.$
By \eqref{chik},
 it is easy to check that
\beaa
 \mathbb{E}(\chi^2_{m}-m)^2 &=& \mathbb{E}[(\chi^2_m)^2]-2m\mathbb{E}(\chi_m^2)+m^2=2m; \\
{\rm Var}\big((\chi^2_m)^2\big)&=& \mathbb{E}[(\chi^2_m)^4]-\big[\mathbb{E}(\chi^2_m)^2\big]^2\\
&=& m(m+2)[(m+4)(m+6)-m(m+2)]\\
&=& 8m(m+2)(m+3)
\eeaa
and
$$\aligned
{\rm Var}\big((\chi^2_m-m)^2\big)&={\rm Var}\big((\chi^2_m)^2\big)+4m^2{\rm Var}\big(\chi^2_m\big)-4m\cdot{\rm Cov}\big((\chi^2_m)^2, \chi^2_m\big)\\
&=8m(m+2)(m+3)+8m^3-4m\big[m(m+2)(m+4)-m^2(m+2)\big]\\
&=8m(m+6),
\endaligned $$ where we use the formula ${\rm Cov}\big((\chi^2_m)^2, \chi^2_m\big)=\mathbb{E}[(\chi^2_m)^3]-\mathbb{E}[(\chi^2_m)^2]\cdot \mathbb{E}(\chi^2_m).$
Similarly by the binomial formula, we have
$$\aligned
\mathbb{E}\big((\chi^2_m-m)^3\big)&=\mathbb{E}\big[(\chi^2_m)^3-3m(\chi^2_m)^2+3m^2(\chi^2_m)-m^3\big]\\
&=m(m+2)(m+4)-3m^2(m+2)+3m^3-m^3\\
&=8m\\
\endaligned $$ and
$$\aligned
\mathbb{E}\big((\chi^2_m-m)^4\big)&=\mathbb{E}\big[(\chi^2_m)^4-4m(\chi^2_m)^3+6m^2(\chi^2_m)^2-4m^3(\chi^2_m)+m^4\big]\\
&=m(m+2)(m+4)(m+6)-4m^2(m+2)(m+4)+6m^3(m+2)\\
&\ \ \ \ \ \ \ \ \ \ \ \ \ \ \ \ \ \ \ \ \ \ \ \ \ \ \ \ \ \ \ \ \ \ \ \ \ -4m^4+m^4\\
&=12m(m+4).
\endaligned $$ The proof is completed.
\hfill$\blacksquare$

\medskip

\medskip

The next result is on Wishart matrices. A Wishart matrix is determined by parameters $p$ and $q$ if it is generated by a random sample from $N_p(\bd{0}, \bd{I}_p)$ with sample size $q$. Let $p=p_n$ and $q=q_n$. Most popular work on this matrix has been taken under the condition $\lim_{n\to\infty}q_n/p_n= c\in (0, \infty).$ For instance, the Marchenko-Pastur distribution \cite{Marchenko}, the central limit theorem (e.g., \cite{Bai_Silverstein}) and the large deviations of its eigenvalues (e.g., \cite{Hiai}) are obtained. The following conclusion is based on the extreme case that $q_n/p_n\to 0$. It is one of the key ingredients in the proof of Theorem \ref{main1}.

\begin{lem}\lbl{fly_ground} Let $\bd{g}_1, \cdots, \bd{g}_q$ be i.i.d. random vectors with distribution $N_p(\bd{0}, \bd{I}_p)$.  Set $\bd{X}_n=(\bd{g}_1, \cdots, \bd{g}_q).$ Let $\lambda_1, \cdots, \lambda_q$ be the eigenvalues of $\bd{X}_n'\bd{X}_n$. Let $p=p_n$ and $q=q_n$ satisfy $p\to\infty,\ q\to \infty,\ \frac{q}{p}\to 0$, then  $\max_{1\leq i \leq q}|\frac{\lambda_i}{p}-1|\stackrel{p}{\to} 0$  as $n\to\infty$.
\end{lem}
\noindent\textbf{Proof}. Review (1.2) from \cite{JiangLi15}. Take $\beta=1$ and treating  $n$ as our ``$q$"  in Theorems 2 and 3 from \cite{JiangLi15}. The rate function $I$ satisfies $I(1)=0$ in both Theorems. By  the large deviations in the two Theorems, we see
\beaa
\frac{1}{p}\max_{1\leq i \leq q}\lambda_i \stackrel{p}{\to} 1\ \ \mbox{and}\ \ \frac{1}{p}\min_{1\leq i \leq q}\lambda_i \stackrel{p}{\to} 1
\eeaa
as $n\to\infty$. The conclusion then follows from the inequality
\beaa
\max_{1\leq i \leq q}\Big|\frac{\lambda_i}{p}-1\Big| \leq \Big|\frac{1}{p}\max_{1\leq i \leq q}\lambda_i-1\Big| + \Big|\frac{1}{p}\min_{1\leq i \leq q}\lambda_i-1\Big|.
\eeaa
The proof is completed. \hfill$\blacksquare$

\begin{lem}\lbl{awesome_pen} Let $\bd{g}_1, \cdots, \bd{g}_q$ be i.i.d. random vectors
with distribution $N_p(\bd{0}, \bd{I}_p)$. Assume $p=p_n\to\infty, q=q_n\to\infty$
 and $\frac{q}{p}\to 0$.
Then, as $n\to\infty$,
\beaa
\frac{1}{pq}\sum_{1\leq i\ne j\leq q}\big[(\bd{g}_i'\bd{g}_j)^2-p\big]
\ \mbox{converges weakly to}\   N(0, 4).
\eeaa
\end{lem}

The proof of Lemma \ref{awesome_pen} is based on a central limit theorem on martingales. Due to its length, we put it as an appendix in Section \ref{Good_Tue}. Figure \ref{fig_Xinmei}, which will be presented later, simulates the densities of $W:=\frac{1}{2pq}\sum_{1\leq i\ne j\leq q}\big[(\bd{g}_i'\bd{g}_j)^2-p\big]$ for various values of $(p, q).$ They indicate that the density of $W$ is closer to the density of $N(0,1)$ as both $p$ and $q$ are larger, and $\frac{q}{p}$ are smaller.

We would like to make a remark on Lemma \ref{awesome_pen} here. Assume $p=1$ instead of the condition $p\to\infty$ in Lemma \ref{awesome_pen}, the conclusion is no longer true. In fact, realizing that $\bd{g}_i$'s are real-valued random variables as $p=1$, we see
\beaa
& & \sum_{1\leq i\ne j\leq q}\big[(\bd{g}_i'\bd{g}_j)^2-1\big]\\
&=& (\bd{g}_1^2+\cdots +\bd{g}_q^2-q)(\bd{g}_1^2+\cdots +\bd{g}_q^2+q)+q-\sum_{i=1}^q\bd{g}_i^4.
\eeaa
By the Slutsky lemma, it is readily seen that $1/(pq^{3/2})\sum_{1\leq i\ne j\leq q}\big[(\bd{g}_i'\bd{g}_j)^2-1\big]$ converges weakly to $N(0, 8)$ as $q\to\infty$. The scaling ``$pq^{3/2}$" here is obviously different from ``$pq$".

\medskip
We will use the following result to prove Lemma \ref{rectangle_performance}.
\blem\label{var-e} Let $\bd{X}=(g_{ij})_{p\times q}$ where $g_{ij}$'s  are independent standard normals. Then

(i) 
${\rm Var}({\rm tr}[({\bf X'X})^2]) =4p^2q^2+8pq(p+q)^2+20pq(p+q+1);$

(ii)  ${\rm Cov}\big({\rm tr}({\bf X}^{\prime}{\bf X}), {\rm tr}[({\bf X}^{\prime}{\bf X})^2]\big)=4pq(p+q+1).$
\nlem
The assertion (i) corrects an error  appeared in (i) of Lemma 2.4 from \cite{Jiang06}, the correct coefficient of the term $p^2q^2$ is ``$4$." However, this does not affect the the main conclusions from \cite{Jiang06}.  The proof of Lemma \ref{var-e} is postponed in Appendix.

In the proof of Theorem \ref{main1}, we will need a slightly more general version of a result from \cite{Jiang06} as follows.
\begin{lem}\lbl{rectangle_performance} Let $\bd{Z}_n$ and $\bd{G}_n$ be as in Theorem \ref{main1}. If $\lim_{n\to\infty}\frac{p_n}{\sqrt{n}}=x\in (0, \infty)$ and $\lim_{n\to\infty}\frac{q_n}{\sqrt{n}}=y\in (0, \infty)$,  then
 $$\liminf_{n\to\infty}\|\ml{L}(\sqrt{n}\bd{Z}_n)-\ml{L}(\bd{G}_n)\|_{\rm TV}\ge\mathbb{E}|e^{\xi}-1|>0,$$
 where $\xi\sim N(-\frac{x^2y^2}{8}, \frac{x^2y^2}{4}).$
\end{lem}
\noindent{\bf Proof.}  By inspecting the proof of Theorem 2 from \cite{Jiang06}, the variable $\xi$ is the limit of random variable $W_n-\frac{x^2y^2}{8}$ with $W_n$ defined in (2.16) of \cite{Jiang06}. Recall
$$
W_n:=\frac{p+q+1}{2n}h_1-\frac{n-p-q-1}{4n^2}h_2$$
where $h_i={\rm tr}({\bf X'X})^i-\mathbb{E}\,{\rm tr}({\bf X'X})^i.$ It is proved in \cite{Jiang06} that $W_n$ converges weakly to a normal random variable with zero mean. What we need to do is to calculate the limit of
${\rm Var}(W_n).$ In fact,

\beaa
{\rm Var}(W_n)&= &\frac{(p+q+1)^2}{4n^2}{\rm Var}\left({\rm tr}({\bf X'X})\right) + \frac{(n-p-q-1)^2}{16n^4}{\rm Var}\left({\rm tr}\left(({\bf X'X})^2\right)\right)\\
& & - \frac{(p+q+1)(n-p-q-1)}{4n^3}\cdot {\rm Cov}\left({\rm tr}({\bf X'X}), {\rm tr}\left(({\bf X'X})^2\right)\right).
\eeaa
Since ${\rm Var}({\rm tr}({\bf X'X}))=2pq$, by Lemma \ref{var-e} we have
$$
{\rm Var}(W_n)=\frac{p^2q^2}{4n^2}+o(1)\to \frac{x^2y^2}{4}
$$
as $n \to \infty.$ Therefore $W_n\to N(0, \frac{x^2y^2}{4}).$ The rest proof is exactly the same as  the proof of Theorem 2 from \cite{Jiang06}.  \hfill$\blacksquare$

\medskip

Let $p=p_n$ and $q=q_n$. We often need the following setting later:
\bea\lbl{service}
q\to \infty,\ \frac{q}{p}\to 0\ \ \mbox{and}\ \ \frac{pq}{n}\to \sigma\in (0,\infty)
\eea
as $n\to \infty$. The next result reveals a subtle property of the eigenvalue part in the density from (\ref{heavy}) under the ``rectangular" case $\frac{q}{p}\to 0$. It is also one of the building blocks in the proof of Theorem \ref{main1}.

\begin{lem}\lbl{sun_half} Let $p=p_n$ and $q=q_n$ satisfy (\ref{service}).
Suppose $\lambda_1, \cdots, \lambda_q$ are the eigenvalues of $\f X_n'\f X_n$ where  $\bd{X}_n=(g_{ij})_{p\times q}$ and $g_{ij}$'s  are independent standard normals. Define
\beaa L_n'=\Big(1-\frac{p}{n}\Big)^{-\frac{1}{2}(n-p-q-1)q}
\Big\{\prod_{i=1}^q\Big(1-\frac{\lambda_i}{n}\Big)\Big\}^{\frac{n-p-q-1}{2}}
\exp\Big(\frac{1}{2}\sum_{i=1}^q\lambda_i\Big).
 \eeaa
Then, as $n\to\infty$,
\beaa
\log L_n'-\frac{pq}{2} + \frac{pq(q+1)}{4(n-p)}\ \mbox{converges weakly to}\  N\big(0, \frac{\sigma^2}{4}\big).
\eeaa
\end{lem}
\noindent\textbf{Proof}. Write
\bea\lbl{half_moon}
\log L_n' &=& \frac{1}{2}\sum_{i=1}^q\lambda_i+ \frac{n-p-q-1}{2}\sum_{i=1}^q\log \frac{n-\lambda_i}{n-p} \nonumber\\
& = & \frac{1}{2}\sum_{i=1}^q\lambda_i+ \frac{n-p-q-1}{2}\sum_{i=1}^q\log \Big(1+ \frac{p-\lambda_i}{n-p}\Big).
\eea
Let function $h(x)$ be such that $\log (1+x)=x-\frac{x^2}{2}+x^3h(x)$ for all $x>-1$. We are able to further write
\bea
& & \sum_{i=1}^q\log \Big(1+ \frac{p-\lambda_i}{n-p}\Big)\nonumber\\
&=&\frac{1}{n-p}\sum_{i=1}^q(p-\lambda_i) -\frac{1}{2(n-p)^2}\sum_{i=1}^q(p-\lambda_i)^2+ \sum_{i=1}^q\Big(\frac{p-\lambda_i}{n-p}\Big)^3h\Big(\frac{p-\lambda_i}{n-p}\Big). \lbl{lantern_red}
\eea
Notice that
\beaa
& & \frac{n-p-q-1}{2}\cdot \frac{1}{n-p}\sum_{i=1}^q(p-\lambda_i)\\
& = & -\frac{1}{2}\sum_{i=1}^q\lambda_i +\frac{(n-p-q-1)pq}{2(n-p)} + \frac{q+1}{2(n-p)}\sum_{i=1}^q\lambda_i.
\eeaa
This, (\ref{half_moon}) and (\ref{lantern_red}) say that
\bea
& & \log L_n' \nonumber\\
&=&\frac{(n-p-q-1)pq}{2(n-p)} + \frac{q+1}{2(n-p)}\sum_{i=1}^q\lambda_i -\frac{n-p-q-1}{4(n-p)^2}\sum_{i=1}^q(p-\lambda_i)^2  \nonumber\\
& & \quad \quad \quad \quad\quad \quad\quad \ \ +\frac{n-p-q-1}{2}\sum_{i=1}^q\Big(\frac{p-\lambda_i}{n-p}\Big)^3h
\Big(\frac{p-\lambda_i}{n-p}\Big). \lbl{lord}
\eea
We now inspect each term one by one. Since $\lambda_1, \lambda_2, \cdots, \lambda_q$ are the eigenvalues of $\f X_n'\f X_n$ and  $\bd{X}_n=(g_{ij})_{p\times q}$, we have
\beaa
\frac{q+1}{2(n-p)}\sum_{i=1}^q\lambda_i
&=&\frac{pq(q+1)}{2(n-p)} + \frac{q+1}{2(n-p)}\sum_{i=1}^p\sum_{j=1}^q (g_{ij}^2-1)\\
& = & \frac{pq(q+1)}{2(n-p)} +\frac{q+1}{2(n-p)}\sqrt{pq}\cdot O_p(1)\\
& = & \frac{pq(q+1)}{2(n-p)} + o_p(1)
\eeaa
by the central limit theorem on i.i.d. random variables. This together with (\ref{lord}) gives
\bea
& & \log L_n' \nonumber\\
&=&\frac{pq}{2} -\frac{n-p-q-1}{4(n-p)^2}\sum_{i=1}^q(\lambda_i-p)^2  \nonumber\\
& & \quad \quad \quad \quad\quad \quad\quad \ \ +\frac{n-p-q-1}{2}\sum_{i=1}^q\Big(\frac{p-\lambda_i}{n-p}\Big)^3h
\Big(\frac{p-\lambda_i}{n-p}\Big)+o_p(1)\nonumber\\
& & \lbl{bed_time}
\eea
as $n\to\infty.$ Now we study $\sum_{i=1}^q(\lambda_i-p)^2$. To do so, set $\bd{X}_n=(\bd{g}_1, \cdots, \bd{g}_q)$. Then $\bd{g}_1, \cdots, \bd{g}_q$ are i.i.d. with distribution $N_p(\bd{0}, \bd{I}_p)$. So $\lambda_1-p, \cdots, \lambda_q-p$ are the eigenvalues of the $q\times q$ symmetric matrix $\bd{X}_n'\bd{X}_n-p\bd{I}_q=(\bd{g}_i'\bd{g}_j)_{q\times q}-p\bd{I}_q.$ Consequently,
\beaa
\sum_{i=1}^q(\lambda_i-p)^2=\sum_{1\leq i\ne j\leq q}(\bd{g}_i'\bd{g}_j)^2 + \sum_{i=1}^q(\|\bd{g}_i\|^2-p)^2.
\eeaa
Now, for $\|\bd{g}_1\|^2 \sim \chi^2(p),$ by Lemma \ref{party_love} we see
$$\aligned
&\mathbb{E}\sum_{i=1}^q(\|\bd{g}_i\|^2-p)^2=q\cdot\mathbb{E}(\|\bd{g}_1\|^2-p)^2=2pq;\\
&\mbox{Var}\,\Big[\sum_{i=1}^q(\|\bd{g}_i\|^2-p)^2\Big] =q\cdot \mbox{Var}\,\Big[(\|\bd{g}_1\|^2-p)^2\Big]
=8pq(p+6).
\endaligned $$
By the Chebyshev inequality,
\beaa
& & \frac{n-p-q-1}{4(n-p)^2}\sum_{i=1}^q(\|\bd{g}_i\|^2-p)^2\\
& = & \frac{(n-p-q-1)pq}{2(n-p)^2}+ \frac{n-p-q-1}{4(n-p)^2}\Big[-2pq+\sum_{i=1}^q(\|\bd{g}_i\|^2-p)^2\Big]\\
& = & \frac{(n-p-q-1)pq}{2(n-p)^2}+ o_p(1)
\eeaa
by noting $\frac{n-p-q-1}{(n-p)^2}\sim \frac{1}{n}$ and $\frac{p^2q}{n^2}\to 0$ as $n\to\infty$. This concludes
\beaa
& & \frac{n-p-q-1}{4(n-p)^2}\sum_{i=1}^q(\lambda_i-p)^2\\ &=&\frac{(n-p-q-1)pq}{2(n-p)^2}+\frac{n-p-q-1}{4(n-p)^2}\sum_{1\leq i\ne j\leq q}(\bd{g}_i'\bd{g}_j)^2+ o_p(1)\\
& = & \frac{(n-p-q-1)pq(q+1)}{4(n-p)^2} + \frac{n-p-q-1}{4(n-p)^2}\sum_{1\leq i\ne j\leq q}\big[(\bd{g}_i'\bd{g}_j)^2-p\big]+ o_p(1).
\eeaa
By splitting $(n-p-q-1)pq(q+1)=(n-p)pq(q+1)-pq(q+1)^2$ and using the fact $\frac{pq^3}{n^2}\to 0$, we see
\bea\lbl{liberation}
\frac{n-p-q-1}{4(n-p)^2}\sum_{i=1}^q(\lambda_i-p)^2-\frac{pq(q+1)}{4(n-p)}\to N\big(0, \frac{\sigma^2}{4}\big)
\eea
weakly, where Lemma \ref{awesome_pen} and  the assertion $\frac{n-p-q-1}{4(n-p)^2}\sim \frac{1}{4n}\sim \frac{\sigma}{4}\cdot \frac{1}{pq}$ are used. Recalling (\ref{bed_time}), to finish our proof, it is enough to show
\bea\lbl{last_last}
\delta_n:=\frac{n-p-q-1}{2}\sum_{i=1}^q\Big(\frac{p-\lambda_i}{n-p}\Big)^3h
\Big(\frac{p-\lambda_i}{n-p}\Big) \stackrel{p}{\to} 0.
\eea

Review $\log (1+x)=x-\frac{x^2}{2}+x^3h(x)$ for all $x>-1$. Then, $\tau:=\sup_{|x|\leq 1/2}|h(x)|< \infty.$ Hence, by the fact $\frac{p}{n}\to 0$ from (\ref{service}),
\bea
P(|\delta_n|>\epsilon)&=& P\Big(|\delta_n|>\epsilon,\ \max_{1\leq i \leq q}|\frac{p-\lambda_i}{n-p}|\leq \frac{1}{2}\Big) + P\Big(\max_{1\leq i \leq q}|\frac{p-\lambda_i}{n-p}|>\frac{1}{2}\Big)\nonumber\\
& \leq & P\Big(|\delta_n|>\epsilon,\ \max_{1\leq i \leq q}|\frac{p-\lambda_i}{n-p}|\leq \frac{1}{2}\Big) + P\Big(\max_{1\leq i \leq q}|\frac{\lambda_i}{p}-1|>\frac{1}{4}\Big) \nonumber\\
& & \lbl{cloud_pee}
\eea
as $n$ is sufficiently large. Under $\max_{1\leq i \leq q}|\frac{p-\lambda_i}{n-p}|\leq \frac{1}{2}$,
\beaa
|\delta_n| &\leq & (2\tau)\cdot \max_{1\leq i \leq q}|\frac{p-\lambda_i}{n-p}|\cdot \frac{n-p-q-1}{4(n-p)^2}\sum_{i=1}^q(\lambda_i-p)^2\\
& = & (2\tau)\cdot \max_{1\leq i \leq q}|\frac{\lambda_i}{p}-1|\cdot \frac{p}{n-p}\cdot \frac{n-p-q-1}{4(n-p)^2}\sum_{i=1}^q(\lambda_i-p)^2
\eeaa
which goes to zero in probability by Lemma \ref{fly_ground}, (\ref{liberation}) and the fact $\frac{p}{n-p}\cdot \frac{pq(q+1)}{4(n-p)}=O(1)$ from the assumption $pq=O(n)$. This, (\ref{cloud_pee}) and Lemma \ref{fly_ground} again conclude (\ref{last_last}).
\hfill$\blacksquare$

\blem\lbl{sunny_time} Let $p_n$ satisfy $p_n/n\to c$ for some $c\in (0,1)$ and $q_n\equiv 1$. Let $\bd{Z}_n$ and $\bd{G}_n$ be as in the first paragraph in Section \ref{chap:intro}. Then $\liminf_{n\to\infty}\|\sqrt{n}\bd{Z}_n - \bd{G}_n\|_{\rm TV} > 0$.
\end{lem}
\noindent\textbf{Proof}. The argument is similar to that of Lemma \ref{sun_half}. By Lemma \ref{del}, the density function of $\sqrt{n}\f Z_n$ is given by
$$
f_n(z):=(\sqrt{2\pi})^{-p}\Big(\frac{2}{n}\Big)^{p/2}
\frac{\Gamma(\frac{n}{2})}{\Gamma(\frac{n-p}{2})}
\Big(1-\frac{|z|^2}{n}\Big)^{(n-p-2)/2}I(|z|<\sqrt{n}),
$$
where $z\in \mathbb{R}^p$ and $p=p_n$. By Lemma \ref{KO},
\beaa
\log \Big[\Big(\frac{2}{n}\Big)^{p/2}
\frac{\Gamma(\frac{n}{2})}{\Gamma(\frac{n-p}{2})}\Big]=-\frac{1}{2}\log (1-c)-\frac{p}{2}-c_n \log\big(1-\frac pn\big)+o(1)\eeaa
as $n\to\infty$, where $c_n:=\frac{1}{2}(n-p-2).$ The density function of $\bd{G}_n$ is $g_n(z)=(\sqrt{2\pi})^{-p}e^{-|z|^2/2}$ for all $z\in \mathbb{R}^p.$ By a measure transformation,
\be\lbl{discovery_find}
\|\ml{L}(\sqrt{n}\bd{Z}_n)-\ml{L}(\bd{G}_n)\|_{\rm TV}=\int_{\mathbb{R}^{pq}}\Big|\frac{f_n(z)}{g_n(z)}-1\Big|g_n(z)\,dz =\mathbb{E}\Big|\frac{f_n(\bd{G}_n)}{g_n(\bd{G}_n)}-1\Big|,
\ee
where the expectation is taken with respect to  random vector $\bd{G}_n$. It is easy to see
\beaa
\log \frac{f_n(z)}{g_n(z)}=-\frac{1}{2}\log (1-c)+c_n\log \frac{n-|z|^2}{n-p}-\frac{p}{2}+\frac{1}{2}|z|^2
\eeaa
if $|z|<\sqrt{n}$, and it is defined to be $-\infty$ if $|z|\geq \sqrt{n}.$  Define function $h(x)$ such that $\log (1+x)=x-\frac{x^2}{2}+x^3h(x)$ for all $x>-1$ and $h(x)=-\infty$, otherwise. Write $\frac{n-|z|^2}{n-p}=1+\frac{p-|z|^2}{n-p}=1+\eta_n(z).$ For convenience, write $\eta_n=\eta_n(z).$ It follows that
\beaa
& & c_n\log \frac{n-|z|^2}{n-p}\\
&=& \frac{1}{2}(n-p-2)\Big[\frac{p-|z|^2}{n-p}-
\frac{1}{2}\Big(\frac{p-|z|^2}{n-p}\Big)^2+\Big(\frac{p-|z|^2}{n-p}\Big)^3
h\Big(\frac{p-|z|^2}{n-p}\Big)\Big]\\
& = & \frac{p}{2}-\frac{|z|^2}{2}-\frac{1}{4}\Big(\frac{p-|z|^2}{\sqrt{n-p}}\Big)^2 + \frac{1}{2(n-p)^{1/2}}\cdot \Big(\frac{p-|z|^2}{\sqrt{n-p}}\Big)^{3}h\Big(\frac{p-|z|^2}{n-p}\Big)\\
& & \quad\quad\quad\quad\quad\quad\quad\quad\quad\quad\quad -\eta_n+\frac{1}{2}\eta_n^2-\eta_n^3h(\eta_n)
\eeaa
for every $|z|<\sqrt{n}$ by using $\frac{1}{2}(n-p-2)=\frac{1}{2}(n-p)-1.$ The last two assertions imply
\bea
\log \frac{f_n(z)}{g_n(z)} &=&-\frac{1}{2}\log (1-c)-\frac{1}{4}\Big(\frac{p-|z|^2}{\sqrt{n-p}}\Big)^2 + \frac{1}{2(n-p)^{1/2}}\cdot \Big(\frac{p-|z|^2}{\sqrt{n-p}}\Big)^{3}h\Big(\frac{p-|z|^2}{n-p}\Big)\nonumber\\
& & \quad\quad\quad\quad\quad\quad\quad\quad\quad\quad\quad\quad\  -\eta_n+\frac{1}{2}\eta_n^2-\eta_n^3h(\eta_n) \lbl{Tea_houses}
\eea
for every $|z|<\sqrt{n}$, and it is identical to $-\infty$ otherwise. Since $\bd{G}_n\sim N_p(\bd{0}, \bd{I}_p)$, we see $\|\bd{G}_n\|^2\sim \chi^2(p)$, $\frac{p-\|\bd{G}_n\|^2}{\sqrt{n-p}}\to N(0, 2c(1-c)^{-1})$ weakly and $\eta_n(\bd{G}_n)\to 0$ in probability. In particular, this implies $h\big(\frac{p-|\bd{G}_n|^2}{n-p}\big)\to 0$ in probability. Finally, by the law of large numbers, $P(\|\bd{G}_n\|<\sqrt{n})\to 1.$ Consequently, from (\ref{Tea_houses}) we conclude
\beaa
\frac{f_n(\bd{G}_n)}{g_n(\bd{G}_n)}\to \frac{1}{\sqrt{1-c}}\cdot\exp\Big\{-\frac{c}{2(1-c)}\chi^2(1)\Big\}
\eeaa
weakly as $n\to\infty$. This and (\ref{discovery_find}) yield the desired conclusion by the Fatou lemma. \hfill$\blacksquare$

\medskip

For a sequence of real numbers $\{a_n\}_{n=1}^{\infty}$ and for a set $I\subset \mathbb{R}$, the notation $\lim_{n\to\infty}a_n\in I$  represents that $\{a_n\}$ has a limit and the limit is in $I$. The next result reveals the strategy about the proof of (ii) of Theorem \ref{main1}.
\begin{lem}\lbl{High_teacher}
For each $n\geq 1$, let $f_n(p,q):\{1, 2,\cdots, n\}^2\to [0, \infty)$ satisfy that  $f_n(p, q)$ is non-decreasing in  $p\in \{1, 2,\cdots, n\}$ and $q\in \{1, 2,\cdots, n\}$, respectively. Suppose
\bea\lbl{red}
\liminf_{n\to\infty}f_n(p_n, q_n)>0
\eea
for any sequence $\{(p_n, q_n);\, 1\leq q_n\leq p_n\leq n\}_{n=1}^{\infty}$
if any of the following  conditions holds:

(i) $q_n\equiv 1$ and  $\lim_{n\to\infty}p_n/n\in (0,1)$;

(ii)  $q_n\to\infty$, $\lim_{n\to\infty} q_n/p_n=0$ and $\lim_{n\to\infty}(p_n q_n)/n\in (0, \infty)$;

(iii) $\lim_{n\to\infty}p_n/\sqrt{n}\in (0,\infty)$ and $\lim_{n\to\infty}q_n/\sqrt{n}\in (0,\infty)$.

\noindent Then (\ref{red}) holds for any sequence $\{(p_n, q_n);\, 1\leq q_n\leq p_n\leq n\}_{n=1}^{\infty}$ satisfying that  $\lim_{n\to\infty}(p_nq_n)/n\in (0, \infty)$.
\end{lem}
\noindent\textbf{Proof}. Suppose the conclusion is not true, that is, $\liminf_{n\to\infty}f_n(p_n, q_n)=0$ for some sequence $\{(p_n, q_n);\, 1\leq q_n\leq p_n\leq n\}_{n=1}^{\infty}$ with $\lim_{n\to\infty}\frac{p_nq_n}{n}=\alpha$, where $\alpha\in (0, \infty)$ is a constant.
Then there exists a subsequence $\{n_k; \, k\geq 1\}$ satisfying  $1\leq q_{n_k}\leq p_{n_k}\leq n_k$ for all $k\geq 1$,  $\lim_{k\to\infty}(p_{n_k}q_{n_k})/n_k=\alpha$ and
\bea\lbl{shawn}
\lim_{k\to\infty}f_{n_k}(p_{n_k}, q_{n_k})=0.
\eea
There are two possibilities: $\liminf_{k\to\infty}q_{n_k}<\infty$ and $\liminf_{k\to\infty}q_{n_k}=\infty$. Let us discuss the two cases separately.

{\it (a)}. Assume $\liminf_{k\to\infty}q_{n_k}<\infty$. Then there exists a further subsequence $\{n_{k_j}\}_{j=1}^{\infty}$ such that $q_{n_{k_j}}\equiv m\geq 1$. For convenience of notation, write $\bar{n}_j=n_{k_j}$ for all $j\geq 1.$ The condition $\lim_{n\to\infty}\frac{p_nq_n}{n}=\alpha$ implies that $\lim_{j\to\infty}\frac{p_{\bar{n}_j}}{\bar{n}_j}=\frac{\alpha}{m}\in (0, 1].$ By (\ref{shawn}) and the monotonocity,
\bea\lbl{privacy_it}
\lim_{j\to\infty}f_{\bar{n}_j}(p_{\bar{n}_j}, 1)=\lim_{j\to\infty}f_{\bar{n}_j}(p_{\bar{n}_j}, q_{\bar{n}_j})=0.
\eea
Define $\tilde{p}_{\bar{n}_j}=[p_{\bar{n}_j}/2]+1$ for all $j\geq 1.$ Then, $\lim_{j\to\infty}\frac{\tilde{p}_{\bar{n}_j}}{\bar{n}_j}=c:=\frac{\alpha}{2m}\in (0, \frac{1}{2}].$ Construct a new sequence such that
\beaa
\tilde{p}_r=
\begin{cases}
\tilde{p}_{\bar{n}_j}, & \text{if $r=\bar{n}_j$ for some $j\geq 1$};\\
[cr] +1, & \text{if not}
\end{cases}
\eeaa
and $\tilde{q}_r=1$ for $r=1,2, \cdots.$ It is easy to check $1\leq \tilde{q}_r\leq \tilde{p}_r\leq r$ for all $r\geq 1$ and $\lim_{r\to\infty}\tilde{p}_r/r=c\in (0,1/2)$. Moreover, $\tilde{p}_{\bar{n}_j}\leq p_{\bar{n}_j}$ for each $j\geq 1.$  So $\{(\tilde{p}_r, \tilde{q}_r);\, r\geq 1\}$ satisfies condition (i), and hence, $\liminf_{r\to\infty}f_r(\tilde{p}_r, \tilde{q}_r)>0$ by (\ref{red}). This contradicts (\ref{privacy_it}) since $f_r(\tilde{p}_r, \tilde{q}_r)=f_{\bar{n}_j}(\tilde{p}_{\bar{n}_j}, 1)\leq f_{\bar{n}_j}(p_{\bar{n}_j}, 1)$ if  $r=\bar{n}_j$ for some $j\geq 1$ by monotonocity.

{\it (b).} Assume $\liminf_{k\to\infty}q_{n_k}=\infty$.  Since $\{q_{n_k}/p_{n_k};\, k\geq 1\} \subset [0,1]$, there is a further subsequence $\{n_{k_j}\}_{j=1}^{\infty}$ such that $q_{n_{k_j}}/p_{n_{k_j}}\to c\in [0, 1]$ as $j\to\infty$. To ease notation, write $\bar{n}_j=n_{k_j}$ for all $j\geq 1.$ Then, $\lim_{j\to\infty}q_{\bar{n}_j}=\infty,  \lim_{j\to\infty}q_{\bar{n}_j}/p_{\bar{n}_j}=c\in [0, 1]$ and $\lim_{j\to\infty}p_{\bar{n}_j}q_{\bar{n}_j}/\bar{n}_j=\alpha\in(0,\infty).$ There are two situations: $c=0$ and $c\in (0,1]$. Let us discuss these cases, respectively.

{\it (b1). $c=0$.} Define
$$
\aligned \tilde{p}_r&=
\begin{cases}
p_{\bar{n}_j}, & \text{if $r=\bar{n}_j$ for some $j\geq 1$};\\
[r^{2/3}], & \text{if not};
\end{cases}
\\
\tilde{q}_r&=
\begin{cases}
q_{\bar{n}_j}, & \text{if $r=\bar{n}_j$ for some $j\geq 1$};\\
([\alpha r^{1/3}\,]+1)\wedge \tilde{p}_r, & \text{if not}.
\end{cases}
\endaligned
$$
 Trivially, $1\leq \tilde{q}_r\leq \tilde{p}_r\leq r$ for all $r\geq 1$ and condition (ii) holds. Moreover, $\tilde{p}_{\bar{n}_j}=p_{\bar{n}_j}$ and $\tilde{q}_{\bar{n}_j}=q_{\bar{n}_j}$ for all $j\geq 1.$ By assumption,
\beaa
\liminf_{j\to\infty}f_{\bar{n}_j}(p_{\bar{n}_j}, q_{\bar{n}_j})\geq \liminf_{n\to\infty}f_n(\tilde{p}_n, \tilde{q}_n)>0.
\eeaa
This contradicts the second equality in (\ref{privacy_it}).

{\it (b2). $c\in (0, 1]$.} In this scenario, $q_{\bar{n}_{j}}/p_{\bar{n}_{j}}\to c\in (0, 1]$. The argument here is similar to {\it (b1)}. Define
\beaa
\tilde{p}_r=
\begin{cases}
p_{\bar{n}_j}, & \text{if $r=\bar{n}_j$ for some $j\geq 1$};\\
([\sqrt{\frac{\alpha}{c} r}]+1)\wedge r, & \text{if not}
\end{cases}
\eeaa
and
\beaa
\tilde{q}_r=
\begin{cases}
q_{\bar{n}_j}, & \text{if $r=\bar{n}_j$ for some $j\geq 1$};\\
[c\tilde{p}_r] \vee 1, & \text{if not}.
\end{cases}
\eeaa
Obviously $1\le \tilde{q}_r\le \tilde{p}_r\le r.$ Since when $r$ is large enough, $\tilde{p}_r\sim\sqrt{\frac{\alpha}{c}}\sqrt{r}$ and $\tilde{q}_r\sim \sqrt{\alpha c}\sqrt{r},$ which means $(\tilde{p}_r, \tilde{q}_r)$ satisfies condition (iii).  We will also get a contradiction by using the same discussion as that of {\it (b1)}.

In conclusion, any of the cases that $\liminf_{k\to\infty}q_{n_k}<\infty$ and $\liminf_{k\to\infty}q_{n_k}=\infty$ results with a contradiction. So our desired conclusion holds true.
\hfill$\blacksquare$

\subsection{The Proof of Theorem \ref{main1}}\lbl{Proof_Main1} The argument is relatively lengthy. We will prove (i) and (ii) separately.

\medskip

\noindent\textbf{Proof of (i) of Theorem \ref{main1}}.
For simplicity, we will use later $p, q$ to replace $p_n, q_n$, respectively,
if there is no confusion. By (\ref{H_TV}) and (\ref{TV_KL}), it is enough to show
\bea\lbl{Comrade_blue}
\lim_{n\to\infty}D_{\rm KL}\big(\ml{L} (\sqrt{n}\f Z_{n})||\f G_{n}\big)=0
\eea
where $\ml{L} (\sqrt{n}\f Z_{n})$ is the probability distribution of $\sqrt{n}\f Z_{n}$.

We can always take two subsequences of $\{n\}$, one of which is such that $q_n\leq p_n$ and the second is $q_n>p_n$. By the symmetry of $p$ and $q$, we only need to prove one of them. So, without loss of generality, we assume $q\leq p$ in the rest of the proof. From the assumption $\lim_{n\to\infty}\frac{pq}{n}=0$, without loss of generality, we assume $p+q<n.$ By Lemma \ref{del}, the density function of $\sqrt{n}\f Z_n$ is
\bea\lbl{miserable_world}
f_n(z):=(\sqrt{2\pi})^{-pq}n^{-pq/2}\frac{\omega(n-p, q)}{\omega(n, q)}\left\{det\left(I_q-\frac{z'z}{n}\right)^{(n-p-q-1)/2}\right\}I_0(z'z/n)
\eea
where $I_0(z'z/n)$ is the indicator function of the set that all
$q$
eigenvalues of $z'z/n$ are in $(0, 1),$ and $\omega(s, t)$ is as in (\ref{fire_you}). Obviously, $g_n(z):= (\sqrt{2\pi})^{-pq} e^{-\mbox{tr}(z'z)/2}$ is the density function of $\bd{G}_n.$

Let $\lambda_1, \cdots, \lambda_q$ be the eigenvalues of $z'z$. Then, $\mbox{det}(I_q-\frac{z'z}{n})=\prod_{i=1}^q(1-\frac{\lambda_i}{n})$ and $\mbox{tr}(z'z)=\sum_{i=1}^q\lambda_i.$
Define \be\lbl{Lequ2}
 L_n=\left\{\prod_{i=1}^q\left(1-\frac{\lambda_i}{n}\right)\right\}^{c_n}\exp\left(\frac{1}{2}\sum_{i=1}^q\lambda_i\right)
\ee
if all $\lambda_i$'s are in $(0,n),$ and $L_n$ is zero otherwise, where $c_n=\frac{1}{2}(n-p-q-1)$. Then one has
\be\lbl{Radon}
\frac{f_n(z)}{g_n(z)}=K_n\cdot L_n
\ee
where $K_n$ is defined as in  \eqref{eq1}. The condition $pq=o(n)$ implies that $\frac{pq^3}{n^2}\to 0$. From Lemma \ref{KO},
\be\lbl{Kmo}\log K_n=-\frac{pq}{2}+\frac{q(q+1)}{4}\log(1+\frac p{n-p})-c_n q\log(1-\frac pn)+o(1) \ee
as $n\to\infty.$ By definition,
\bea\lbl{KL-dis}
D_{\rm KL}\big(\ml{L} (\sqrt{n}\f Z_{n})||\f G_{n}\big)
&=& \int_{\mathbb{R}^{pq}} \Big[\frac{f_n(z)}{g_n(z)}\log \frac{f_n(z)}{g_n(z)}\Big] g_n(z)\, dz\nonumber\\
&=&\mathbb{E} \log\frac{f_n(\sqrt{n}\f Z_n)}{g_n(\sqrt{n}\f Z_n)}\nonumber\\
&=&\mathbb{E}\log[K_n\cdot L_n],
\eea
where $\lambda_1, \cdots, \lambda_q$ are the eigenvalues
of $n\f Z_n'\f Z_n$ since $f_n(z)$ is the density function of $\sqrt{n}\f Z_n$. We also define
$\log 0=0$ since random variable $\frac{f_n(\sqrt{n}\f Z_n)}{g_n(\sqrt{n}\f Z_n)}>0$ a.s. The definition
 of $I_0(z'z/n)$ from (\ref{miserable_world}) ensures that $I_0(z'z/n)=0$ if  $\max_{1\leq i \leq q}\lambda_i\geq n$ a.s. By Lemma \ref{KL-Key}, $\mathbb{E}\sum_{i=1}^q\lambda_i=pq$. This and (\ref{Kmo}) imply that
the expectation in (\ref{KL-dis}) is further equal to
\bea
&&\log K_n+\frac12\mathbb{E}\sum_{i=1}^q \lambda_i+c_n \mathbb{E}\sum_{i=1}^{q}
\log\big(1-\frac{\lambda_i}{n}\big)\nonumber\\
&=&\frac{q(q+1)}{4}\log\big(1+\frac p{n-p}\big)+c_n\mathbb{E}\sum_{i=1}^q
\log\big(1+\frac{p-\lambda_i}{n-p}\big)+o(1)\lbl{portland}\\
&\le&\frac{pq^2}{4(n-p)}+c_n\mathbb{E}\sum_{i=1}^q
\Big[\frac{p-\lambda_i}{n-p}-\frac{(\lambda_i-p)^2}{2(n-p)^2}+\frac{(p-\lambda_i)^3}{3(n-p)^3} \Big]+o(1),\nonumber
\eea
where we combine the term ``$-c_n q\log(1-\frac pn)$" from (\ref{Kmo}) with ``$\log\big(1-\frac{\lambda_i}{n}\big)$" to get the sum in (\ref{portland}),
and the last step is due to  the elementary inequality
$$\log(1+x)\le x-\frac{x^2}{2}+\frac{x^3}{3}$$
for any $x>-1$. Based on Lemmas \ref{KL-Key} and \ref{tri}, we know that, under the condition  $pq=o(n),$
$$\aligned
&\mathbb{E}\sum_{i=1}^q \lambda_i=pq, \quad \mathbb{E}\sum_{i=1}^q\lambda_i^2=pq(p+q)+O(pq), \\
& \mathbb{E}\sum_{i=1}^q\lambda_i^3=pq(p^2+q^2+3pq)+O(p^2q)
\endaligned $$
(the ``$\lambda_i$" here is $n$ times the ``$\lambda_i$" from Lemmas \ref{KL-Key} and \ref{tri}). These imply that
\bea\lbl{main1-esti}\aligned
\mathbb{E}\sum_{i=1}^q(p-\lambda_i)&=pq-\mathbb{E}\sum_{i=1}^q \lambda_i=0; \\
 \mathbb{E}\sum_{i=1}^q(p-\lambda_i)^2&=p^2q-2p\mathbb{E}
\sum_{i=1}^q \lambda_i+\mathbb{E}\sum_{i=1}^q \lambda_i^2=pq^2+O(pq);\\
 \mathbb{E}\sum_{i=1}^q(p-\lambda_i)^3&=p^3q-3p^2\mathbb{E}\sum_{i=1}^q \lambda_i+3p\mathbb{E}\sum_{i=1}^q \lambda_i^2-\mathbb{E}\sum_{i=1}^q \lambda_i^3\\
&=-pq^3+O(p^2q).
\endaligned \eea
Recall that $c_n=\frac{1}{2}(n-p)-\frac{1}{2}(q+1).$ Plugging \eqref{main1-esti} into \eqref{portland}, we get from (\ref{KL-dis}) that
\beaa\label{klestimate2}
 D_{\rm KL}\bigg(\ml{L} (\sqrt{n}\f Z_{n})||\f G_{n}\bigg)
 &\leq &\frac{pq^2}{4(n-p)}+c_n\bigg[-\frac{pq^2+O(pq) }{2(n-p)^2}-
 \frac{pq^3+O(p^2q)}{3(n-p)^3}\bigg]\\
&= & \frac{pq^2}{4(n-p)}-\frac{pq^2}{4(n-p)}+o(1)\to 0,
\eeaa
where we use the following two limits:
\beaa
& &  \frac{q+1}{2}\cdot
\frac{pq^2+O(pq) }{2(n-p)^2}=O\Big(\frac{pq^3+pq^2}{n^2}\Big)\to 0;\\
& & c_n\cdot \frac{pq^3+O(p^2q)}{(n-p)^3}=O\Big(\frac{pq^3 + p^2q}{n^2}\Big)\to 0\\
\eeaa
 by \eqref{service}. This gives (\ref{Comrade_blue}).
\hfill$\blacksquare$

\medskip

Let $(U_1, V_1)'\in \mathbb{R}^m$ and $(U_2, V_2)'\in \mathbb{R}^m$ be two random vectors with $U_1 \in \mathbb{R}^s, U_2 \in \mathbb{R}^s$ and $V_1 \in \mathbb{R}^t, V_2 \in \mathbb{R}^t$ where $s\geq 1,\, t\geq 1$ and $s+t=m$. It is easy to see from the first identity of (\ref{norm}) that
\bea\lbl{Sunday_Kowalski}
\|\mathcal{L}(U_1, V_1)-\mathcal{L}(U_2, V_2)\|_{\rm TV} \geq \|\mathcal{L}(U_1)-\mathcal{L}(U_2)\|_{\rm TV}
\eea
by taking (special) rectangular sets in the supremum.

\medskip

\noindent\textbf{Proof of (ii) of Theorem \ref{main1}}. Remember that our assumption is
$\lim_{n\to\infty}\frac{p_nq_n}{n}=\sigma\in(0, \infty).$ By the argument at the beginning
of the proof of (i) of Theorem \ref{main1}, without loss of generality,
we assume $q_n\leq p_n$ for all $n\geq 3.$ By (\ref{H_TV}) and (\ref{TV_KL}),
it suffices to show
\bea\lbl{culture_mouse}
\liminf_{n\to\infty}\|\ml{L}(\sqrt{n}\bd{Z}_n)-\ml{L}(\bd{G}_n)\|_{\rm TV}>0.
\eea

Define $f_n(p,q)=\|\ml{L}(\sqrt{n}\bd{Z}_n)-\ml{L}(\bd{G}_n)\|_{\rm TV}$ for $1\leq q\leq p \leq n.$ Here we slightly abuse the  notation: $\bd{Z}_n$ and $\bd{G}_n$ are $p\times q$ matrices with $p$ and $q$ being arbitrary instead of  fixed sizes $p_n$ and $q_n.$ From (\ref{Sunday_Kowalski}) it is immediate that $f_n(p,q)$ is non-decreasing in $p\in\{1,\cdots, n\}$ and $q\in \{1,\cdots, n\}$, respectively. Then,  by Lemmas \ref{rectangle_performance}, \ref{sunny_time} and \ref{High_teacher}, it is enough to prove \eqref{culture_mouse} under assumption
(\ref{service}). For simplicity, from now on  we will write  $p$ for $p_n$ and $q$ for $q_n$, respectively.  Remember the joint density function of entries of $\f Z_n$ is the function $f_n(z)$ defined in \eqref{heavy} and $g_n(z):= (\sqrt{2\pi})^{-pq} e^{-\mbox{tr}(z'z)/2}$ is the density function of $\bd{G}_n.$
Set
 \beaa
& &  L_n'=\left(1-\frac{p}{n}\right)^{-\frac{1}{2}(n-p-q-1)q}
 L_n;\\
& & K_n'= \left(1-\frac{p}{n}\right)^{\frac{1}{2}(n-p-q-1)q}K_n,
 \eeaa
 where $K_n$ and $L_n$ are defined by \eqref{eq1} and \eqref{Lequ2}, respectively.
Evidently,
\beaa
K_n'\cdot L_n'=K_n\cdot L_n.
\eeaa
By the expression \eqref{Radon},
we have
\be\lbl{111}
\frac{f_n(z)}{g_n(z)}=K_n\cdot L_n=K_n'\cdot L_n'.
\ee
Then by definition,
\be\lbl{222}
\|\ml{L}(\sqrt{n}\bd{Z}_n)-\ml{L}(\bd{G}_n)\|_{\rm TV}=\int_{\mathbb{R}^{pq}}\Big|\frac{f_n(z)}{g_n(z)}-1\Big|g_n(z)\,dz
=
\mathbb{E}\Big|\frac{f_n(\bd{G}_n)}{g_n(\bd{G}_n)}-1\Big|,
\ee
where the expectation is taken over random matrix $\bd{G}_n$.
From (\ref{111}) and (\ref{222}), we have
\bea\lbl{333}
\|\ml{L}(\sqrt{n}\bd{Z}_n)-\ml{L}(\bd{G}_n)\|_{\rm TV}= \mathbb{E}|K_n'L_n'-1|,
\eea
where $\lambda_1, \cdots, \lambda_q$ are the eigenvalues of the Wishart matrix $\bd{G}_n'\bd{G}_n$.

First, we know $\frac{pq^3}{n^2}\to 0$ by (\ref{service}). Use (\ref{KmO}) to see
\beaa\lbl{free_hot}
\log K_n'=-\frac{pq}{2}+\frac{q(q+1)}{4}\log \big(1+\frac p{n-p}\big)+o(1)
\eeaa
as $n\to\infty$. By Taylor's expansion,
\beaa
\log \big(1+\frac p{n-p}\big)&=&\frac p{n-p}-\frac {p^2}{2(n-p)^2} +O\big(\frac{p^3}{n^3}\big).
\eeaa
We then get
\beaa
\log K_n'=-\frac{pq}{2}+\frac{pq(q+1)}{4(n-p)}-\frac{\sigma^2}{8} +o(1)
\eeaa
by (\ref{service}). This and Lemma \ref{sun_half} yield
\beaa
\log (K_n'L_n') \to N\big(-\frac{\sigma^2}{8}, \frac{\sigma^2}{4}\big)
\eeaa
weakly as $n\to\infty.$  This implies that $K_n'L_n'$ converges weakly to $e^{\xi}$, where $\xi\sim N\big(-\frac{\sigma^2}{8}, \frac{\sigma^2}{4}\big)$. By (\ref{333}) and the Fatou lemma
\beaa
\liminf_{n\to\infty}\|\ml{L}(\sqrt{n}\bd{Z}_n)-\ml{L}(\bd{G}_n)\|_{\rm TV}\geq \mathbb{E}|e^{\xi}-1|>0.
\eeaa
The proof is completed. \hfill$\blacksquare$

\section{Proof of Theorem \ref{ttt}}\lbl{yes_ttt}

There are two parts in this section. We first need a preparation  and then prove  Theorem \ref{ttt}.

\subsection{Auxiliary Results}\lbl{papapa}
Review the Hilbert-Schmidt norm defined in  (\ref{Hiha}). A limit theorem on  the norm appeared in Theorem \ref{ttt} is given for a special case.

\blem\lbl{piano_black} Let $p=p_n$ satisfy $p_n/n\to c$ for some $c\in (0,1]$.
Let $\bd{\Gamma}_{p\times 1}$ and $\bd{Y}_{p\times 1}$ be as in Theorem \ref{ttt}. Then $\|\sqrt{n}\bd{\Gamma}_{p\times 1}-\bd{Y}_{p\times 1}\|_{\rm HS} \to \sqrt{\frac{c}{2}}\cdot |N(0,1)|$ weakly as $n\to\infty$.
\nlem
\noindent\textbf{Proof}. Let $\bd{y}=(\xi_1, \cdots, \xi_n)'\sim N_n(\bd{0}, \bd{I}_n).$ By the Gram-Schmidt algorithm, $(\bd{y}, \frac{\bd{y}}{\|\bd{y}\|})$ and $(\bd{Y}_{n\times 1}, \bd{\Gamma}_{n\times 1})$ have the same distribution. By the definition of the Hilbert-Schmidt norm, it is enough to show
\beaa\lbl{coffee_tea}
\|\sqrt{n}\bd{\Gamma}_{p\times 1}-\bd{Y}_{p\times 1}\|_{\rm HS}^2\overset{d}{=}\Big(\frac{\sqrt{n}}{\|\bd{y}\|}-1\Big)^2\cdot \sum_{i=1}^p\xi_i^2\to \frac{c}{2}\cdot N(0,1)^2
\eeaa
as $n\to\infty$, where $\|\bd{y}\|^2=\xi_1^2+\cdots + \xi_n^2$. In fact, the middle term of the above is equal to
\beaa
& & \frac{(\|\bd{y}\|^2-n)^2}{\|\bd{y}\|^2(\|\bd{y}\|+\sqrt{n})^2}\sum_{i=1}^p\xi_i^2\\
&= & \frac{p}{n}\cdot \Big(\frac{\|\bd{y}\|^2-n}{\sqrt{n}}\Big)^2\cdot \Big(\frac{\|\bd{y}\|^2}{n}\Big)^{-1}\cdot \Big[1+\Big(\frac{\|\bd{y}\|^2}{n}\Big)^{1/2}\Big]^{-2} \cdot \frac{1}{p}\sum_{i=1}^p\xi_i^2.
\eeaa
By the classical law of large numbers and the central limit theorem, $\frac{\|\bd{y}\|^2}{n}\to 1$ in probability and $\frac{\|\bd{y}\|^2-n}{\sqrt{n}}\to N(0,2)$ weakly as $n\to\infty$. By the Slutsky lemma, the above converges weakly to $\frac{c}{2}\cdot N(0,1)^2$. \hfill$\blacksquare$

\medskip

Review the notation before the statement of Theorem \ref{ttt}. Set
\bea\lbl{notationwise}
\bd{\Sigma}_k:=(\bm{\gamma}_1, \cdots, \bm{\gamma}_k)(\bm{\gamma}_1, \cdots, \bm{\gamma}_k)'=\sum_{i=1}^k\bm{\gamma}_i\bm{\gamma}_i'
\eea
for $1\leq k \leq n$. Easily, $\bd{\Sigma}_k$ has rank $k$ almost surely and it is an idempotent matrix, that is, $\bd{\Sigma}_k^2 =\bd{\Sigma}_k$.
It is easy to check that  $\bd{w}_k=(\bd{I}-\bd{\Sigma}_{k-1})\bd{y}_k$ for $2\leq k \leq n.$

\blem\lbl{party}
Let $1\leq k\leq n$ be given. Then, $\|\bd{w}_k\|^2\sim \chi^2(n-k+1)$ and  $\|\bd{\Sigma}_{k-1}\bd{y}_k\|^2 \sim \chi^2(k-1)$. Further, given $\bd{y}_1, \cdots, \bd{y}_{k-1}$, the two conclusions still hold, and
$\|\bd{w}_k\|^2$ and $\|\bd{\Sigma}_{k-1}\bd{y}_k\|^2$ are conditionally independent.
\nlem
\noindent\textbf{Proof}. First, let us review the following fact.
Suppose $\bd{y}\sim N_n(\bd{0}, \bd{I}_n)$ and $\bd{A}$ is an $n\times n$ symmetric matrix with eigenvalues
$\lambda_1, \cdots, \lambda_n.$  Then
\bea\lbl{normal_fact}
\bd{y}'\bd{A}\bd{y}\ \ \mbox{and}\ \ \sum_{i=1}^n\lambda_i\xi_i^2\ \ \mbox{have the same distribution}.
\eea
In particular,
\bea\lbl{variance2}
\mbox{Var}\, (\bd{y}'\bd{A}\bd{y})=\sum_{i=1}^n\lambda_i^2\,\mbox{Var}\,(\xi_i^2)=2\,\mbox{tr}\,(\bd{A}^2).
\eea
If $\bd{A}$ is an idempotent matrix with rank $r$,
then all of the nonzero eigenvalues of $\bd{A}$ are $1$ with $r$-fold. Thus, $\|\bd{A}\bd{y}\|^2=\bd{y}'\bd{A}\bd{y}\sim \chi^2(r)$. Moreover,
the distribution of $\|\bd{A}\bd{y}\|^2$ depends only on the rank of $\bd{A}.$
Therefore, all conclusions follow except the one on conditional independence.
Now we prove it.

Given $\bd{y}_1, \cdots, \bd{y}_{k-1}$, we see that $\bd{w}_k=(\bd{I}-\bd{\Sigma}_{k-1})\bd{y}_k$ and $\bd{\Sigma}_{k-1}\bd{y}_k$ are two Gaussian random vectors. By using the fact that $\bd{y}_1, \cdots, \bd{y}_k$ are i.i.d. random vectors, we  see that the conditional covariance matrix
\beaa
& & \mathbb{E}[\bd{w}_k(\bd{\Sigma}_{k-1}\bd{y}_k)'\big|\bd{y}_1, \cdots, \bd{y}_{k-1}]\\
& = &  \mathbb{E}[(\bd{I}-\bd{\Sigma}_{k-1})\bd{y}_k\bd{y}_k'\bd{\Sigma}_{k-1}\big|\bd{y}_1, \cdots, \bd{y}_{k-1}]\\
& =& (\bd{I}-\bd{\Sigma}_{k-1})\mathbb{E}(\bd{y}_k\bd{y}_k')\bd{\Sigma}_{k-1}  =0
\eeaa
since $\mathbb{E}(\bd{y}_k\bd{y}_k')=\bd{I}_n$. This implies that $\bd{w}_k$ and $\bd{\Sigma}_{k-1}\bd{y}_k$ are conditionally independent. \hfill$\blacksquare$

\medskip

We next expand the trace of a target matrix in terms of its entries. Then the expectation of the trace can be computed explicitly via Lemma \ref{what}.
\blem\lbl{Togo}  Let $\bd{\Gamma}_n=(\bm{\gamma}_1, \cdots, \bm{\gamma}_n)=(\gamma_{ij})$  be an $n\times n$ matrix. Set $\bd{\Sigma}_k=\sum_{i=1}^k\bm{\gamma}_i\bm{\gamma}_i'$ for $1\leq k \leq n.$ Given $1\leq p \leq n$, denote the upper-left $p\times p$ submatrix of $\bd{\Sigma}_{k}$ by $(\bd{\Sigma}_{k})_{p}$. Then
\beaa
{\rm tr}[(\bd{\Sigma}_{k})_p^2]=\sum_{i=1}^k\sum_{r=1}^p\gamma_{ri}^4 +\sum_{i=1}^k\sum_{1\leq r\ne  s\leq p}\gamma_{ri}^2\gamma_{si}^2& +& \sum_{r=1}^p\sum_{1\leq i \ne j \leq k}\gamma_{ri}^2\gamma_{rj}^2\\
&+ & \sum_{1\leq i \ne j \leq k}\sum_{1\leq r\ne  s\leq p}\gamma_{ri}\gamma_{si}\gamma_{rj}\gamma_{sj}.
\eeaa
\nlem
\noindent\textbf{Proof}. The argument is similar to that of (\ref{additional_care}). However, the following care has to be taken additionally. Observe the $(r,s)$-element of $(\bm{\gamma}_i\bm{\gamma}_i')$ is $\gamma_{ri}\gamma_{si}$.
Since $(\bd{\Sigma}_{k})_{p}=\sum_{i=1}^k(\bm{\gamma}_i\bm{\gamma}_i')_{p}$, we know the
$(r,s)$-element of the symmetric matrix $(\bd{\Sigma}_{k})_{p}$ is $\sum_{i=1}^k\gamma_{ri}\gamma_{si}$ for $1\leq r, s \leq p.$ Note that ${\rm tr}(\bd{U}^2)=\sum_{1\leq r, s\leq p}u_{rs}^2$ for any symmetric matrix $\bd{U}=(u_{ij})_{p\times p}$. We have
\beaa
{\rm tr}[(\bd{\Sigma}_{k})_p^2]=\sum_{1\leq r, s\leq p}\Big(\sum_{i=1}^k\gamma_{ri}\gamma_{si}\Big)^2
=\sum_{1\leq i, j \leq k}\sum_{1\leq r, s\leq p}\gamma_{ri}\gamma_{si}\gamma_{rj}\gamma_{sj}.
\eeaa
Divide  the first sum into two sums corresponding to that $i=j$ and that $i\ne j$, respectively. Similarly, for the second sum, consider the case $r=s$ and the case $r \ne s$, respectively. The conclusion then follows. \hfill$\blacksquare$

\medskip

The study of the trace norm appearing in Theorem \ref{ttt} is essentially reduced to a sum; see the first statement next. It discloses the behavior of the sum on the ``boundary" case.
\blem\lbl{trofSigma}  Let $\bd{\Gamma}_n=(\bm{\gamma}_1, \cdots, \bm{\gamma}_n)=(\gamma_{ij})$ be the
 $n\times n$ Haar-invariant orthogonal matrix generated from $\bd{Y}_n=(\bd{y}_1, \cdots, \bd{y}_n)$ as in  (\ref{sea}). Let
 $\bd{\Sigma}_k$ be as in (\ref{notationwise}) and $(\bd{\Sigma}_{k})_{p}$
 denote the upper-left $p\times p$ submatrix of $\bd{\Sigma}_{k}.$ The following hold.
\begin{itemize}
\item[1)] For $\sigma>0$, assume $p\geq 1$, $q\to\infty$ with $\frac{pq^2}{n}\to\sigma.$ Then, as  $n\to\infty,$
$$\sum_{j=2}^q{\rm tr}[(\bd{\Sigma}_{j-1})_p]\stackrel{p}{\longrightarrow}\frac{\sigma}2.$$
\item[2)] For any $q\geq 2$,
 $$\sum_{j=2}^{q}\mathbb{E}\,{\rm tr}[(\bd{\Sigma}_{j-1})_p^2]=\frac{pq(q-1)(p+2)}{2n(n+2)}+\frac{pq(q-1)(q-2)(n-p)}{3n(n-1)(n+2)}.
 $$
\end{itemize}
\nlem
\noindent\textbf{Proof}. To prove 1), it is enough to show that
\bea\lbl{monkey}
\mathbb{E}\sum_{j=2}^q{\rm tr}[(\bd{\Sigma}_{j-1})_p]=\frac{pq(q-1)}{2n}\ \ \mbox{and}\  \quad {\rm Var}\Big(\sum_{j=2}^q{\rm tr}[(\bd{\Sigma}_{j-1})_p]\Big)\to 0
\eea
as $n\to\infty.$
 Recall
$
\bd{\Sigma}_{j-1}=
\sum_{k=1}^{j-1}\bm{\gamma}_k\bm{\gamma}_k'
$ and the $(r,s)$-element of $(\bm{\gamma}_k\bm{\gamma}_k')$ is $\gamma_{rk}\gamma_{sk}$. For convenience, define
$u_k=\sum_{s=1}^p
\gamma_{sk}^2$ for any $1\le k\le q.$
 Then,
\beaa
{\rm tr}[(\bd{\Sigma}_{j-1})_p]=\sum_{k=1}^{j-1}{\rm tr}[(\bm{\gamma}_k\bm{\gamma}_k')_p]=\sum_{k=1}^{j-1}\sum_{s=1}^p \gamma_{sk}^2=\sum_{k=1}^{j-1}u_k
\eeaa
for each $2\le j\le q.$ It follows that
\bea\lbl{egg}
\sum_{j=2}^q {\rm tr}[(\bd{\Sigma}_{j-1})_p]=\sum_{j=2}^q \sum_{k=1}^{j-1} u_k=\sum_{k=1}^{q-1}(q-k)u_k.
\eea
We claim that
\be\label{u}\mathbb{E}(u_k)=\frac pn, \quad \mathbb{E}(u_k^2)=\frac{p(p+2)}{n(n+2)}, \quad {\rm Cov}(u_i, u_k)=-\frac{2p(n-p)}{n^2(n-1)(n+2)}
\ee
for any $1\le i\neq k\le q.$ In fact, by {\bf F2}) and  Lemma \ref{what}, it is immediate to see  $\mathbb{E}(u_k)=\frac pn$. Further, by the same argument,
\beaa
\mathbb{E}u_k^2&=&p \mathbb{E}(\gamma_{11}^4) + p(p-1) \mathbb{E}(\gamma_{11}^2\gamma_{12}^2)\\
& = & \frac{3p}{n(n+2)}+ \frac{p(p-1)}{n(n+2)}=\frac{p(p+2)}{n(n+2)}.
\eeaa
Now we turn to prove the third conclusion from (\ref{u}). For any $i\neq k,$ by {\bf F2}) again,
\beaa
{\rm Cov}(u_i,u_k)&=& \mathbb{E}(u_iu_k)-\frac{p^2}{n^2}\\
&=&\sum_{s=1}^p\mathbb{E}(\gamma_{si}^2\gamma_{sk}^2)+\sum_{1\le s\neq t\le p}\mathbb{E}(\gamma_{si}^2 \gamma_{tk}^2)-\frac{p^2}{n^2}\\
&=& p\mathbb{E}(\gamma_{11}^2\gamma_{12}^2)+p(p-1)\mathbb{E}(\gamma_{11}^2\gamma_{22}^2)-\frac{p^2}{n^2}\\
&=& \frac{p}{n(n+2)}+p(p-1)\frac{n+1}{n(n-1)(n+2)}-\frac{p^2}{n^2}\\
&=& -\frac{2p(n-p)}{n^2(n-1)(n+2)},
\eeaa
where we use Lemma \ref{what} for the fourth equality. So claim (\ref{u}) follows.

Now, let us go back to the formula in (\ref{egg}). By (\ref{u}),
\beaa
\lbl{e} \sum_{j=2}^q \mathbb{E}{\rm tr}[(\bd{\Sigma}_{j-1})_p]=\sum_{k=1}^{q-1}(q-k)\mathbb{E} u_k=\frac{pq(q-1)}{2n}.
\eeaa
The first identity from (\ref{monkey}) is concluded. Now we work on the second one.
It is readily seen from the first two conclusions of (\ref{u}) that
${\rm Var}(u_k)=\frac{2p(n-p)}{n^2(n+2)}.$ By \eqref{u} again,
\beaa
{\rm Var}\Big(\sum_{j=2}^q {\rm tr} [(\bd{\Sigma}_{j-1})_p]\Big)
&=&\sum_{k=1}^{q-1}(q-k)^2 {\rm Var}(u_k)+\sum_{1\le i\neq k\le q-1}(q-i)(q-k){\rm Cov}(u_i, u_k)\\
&=&A\sum_{k=1}^{q-1}(q-k)^2-B\sum_{1\le i\neq k\le q-1}(q-i)(q-k)\\
& = & (A+B)\sum_{r=1}^{q-1}r^2 -B\Big(\sum_{r=1}^{q-1}r\Big)^2
\eeaa
by setting $r=q-k$ and $s=q-i$, respectively, where
\beaa
A=\frac{2p(n-p)}{n^2(n+2)}\ \ \ \mbox{and}\ \ \ B=\frac{2p(n-p)}{n^2(n-1)(n+2)}.
\eeaa
From an elementary calculation, we get
$$\aligned
{\rm Var}\Big(\sum_{j=2}^q {\rm tr} [(\bd{\Sigma}_{j-1})_p]\Big)
&=\frac{pq(n-p)(q-1)}{3n(n-1)(n+2)}\Big(2q-1-\frac{3q(q-1)}{2n}\Big)\lbl{come_on}\\
&= O\Big(\frac{pq^3}{n^2}\Big)
 = O\Big(\frac{1}{pq}\Big)\to 0
\endaligned $$
as $ n\to\infty$, where we use the fact $\frac{3q(q-1)}{2n}=O(1)$ under the assumption $pq^2/n\to\sigma$ and $q\to\infty.$ This gives the second conclusion from (\ref{monkey}).

Now we prove 2).  By {\bf F2})  and Lemma \ref{Togo},
\beaa
& & \mathbb{E}{\rm tr}[(\bd{\Sigma}_{j-1})_p^2]\\
&=& p(j-1)\mathbb{E}(\gamma_{11}^4) + [(j-1)p(p-1) + p(j-1)(j-2)]\mathbb{E}(\gamma_{11}^2\gamma_{12}^2)\\
& & \ \ \ \ \ \ \ \ \ \ \ \ \ \ \ \ \ \ \ + \, p(p-1)(j-1)(j-2)\mathbb{E}(\gamma_{11}\gamma_{12}\gamma_{21}\gamma_{22}).
\eeaa
Write
\beaa
& & (j-1)p(p-1) + p(j-1)(j-2)=(j-1)p(p-2)+p(j-1)^2;\\
& & p(p-1)(j-1)(j-2)=p(p-1)\big[(j-1)^2-(j-1)\big].
\eeaa
Then, by computing $\sum_{j=1}^{q}(j-1)$ and $\sum_{j=1}^{q}(j-1)^2$, we obtain
\beaa
& & \sum_{j=2}^{q}\mathbb{E}\,{\rm tr}[(\bd{\Sigma}_{j-1})_p^2]\\
&=&
\frac{1}{2}pq(q-1)\mathbb{E}(\gamma_{11}^4) + \frac{1}{6}pq(q-1)\big[3p+2q-7\big]\mathbb{E}(\gamma_{11}^2\gamma_{12}^2)\\
& & \ \ \ \ \ \ \ \ \ \ \ \ \ \ \ \ \ \ \ \ \ + \,
\frac{1}{3}pq(p-1)(q-1)(q-2)\mathbb{E}(\gamma_{11}\gamma_{12}\gamma_{21}\gamma_{22}).
\eeaa
From  Lemma \ref{what}, it is trivial to get
\beaa
& & \sum_{j=2}^{q}\mathbb{E}\,{\rm tr}[(\bd{\Sigma}_{j-1})_p^2]
=
\frac{pq(q-1)(p+2)}{2n(n+2)}+\frac{pq(q-1)(q-2)(n-p)}{3n(n-1)(n+2)}.
\eeaa
The proof is completed.  \hfill$\blacksquare$

\medskip

Our target is a submatrix of an Haar-orthogonal matrix. Based on the argument in the proof of Lemma \ref{trofSigma}, an estimate of the submatrix is provided now.

\bprop\label{lala} Let $\bd{\Gamma}_{p\times q}$ and $\bd{Y}_{p\times q}$ be
 as in  (\ref{Hiha}). Then
 \beaa
 \mathbb{E}\|\sqrt{n}\bd{\Gamma}_{p\times q}-\bd{Y}_{p\times q}\|^2_{\rm HS} \leq \frac{24pq^2}{n}
\eeaa
for any $n\ge 2$ and $1\leq p, q\leq n.$
\nprop

\noindent\textbf{Proof}.  Review the notation from (\ref{sea}), identity
$\bd{w}_j=(\bd{I}-\bd{\Sigma}_{j-1})\bd{y}_j$ and $\bm{\gamma}_j={\bd w}_j/\|{\bf w}_j\|$.
We first write
\bea\lbl{aha}
\sqrt{n}\bm{\gamma}_j-\bd{y}_j=\sqrt{n}{\bm \gamma}_j-{\bf w}_j
-\bd{\Sigma}_{j-1}\bd{y}_j=(\sqrt{n}-\|{\bf w}_j\|){\bm \gamma}_j-\bd{\Sigma}_{j-1}\bd{y}_j
\eea
for $1\leq j \leq n$, where $\bd{\Sigma}_{0}=\bd{0}$. Define $\bd{M}=\bd{M}_{p\times n}=(\bd{I}_p, \bd{0})$ for $1\leq p \leq n-1$ and $\bd{M}_{n\times n}=\bd{I}_n$, where $\bd{I}_p$ is the $p\times p$ identity matrix and $\bd{0}$ is the $p\times (n-p)$ matrix whose entries are all equal to zero. Evidently, $\bd{M}(\sqrt{n}\bm{\gamma}_j-\bd{y}_j)$ is the upper $p$-dimensional vector of  $\sqrt{n}\bm{\gamma}_j-\bd{y}_j$. Hence, by (\ref{aha}),
\be\lbl{expnorm}\aligned
\|\sqrt{n}\bd{\Gamma}_{p\times q}-\bd{Y}_{p\times q}\|^2_{\rm HS}&=\sum_{j=1}^q\sum_{i=1}^p(\sqrt{n}\gamma_{ij}-y_{ij})^2 \\
&= \sum_{j=1}^q\|\bd{M}(\sqrt{n}\bm{\gamma}_j-\bd{y}_j)\|^2 \\
&=\sum_{j=1}^q\big\|(\sqrt{n}-\|{\bf w}_j\|)\bd{M}{\bm \gamma}_j-
\bd{M}\bd{\Sigma}_{j-1}\bd{y}_j\big\|^2\\
& \leq 2\sum_{j=1}^q(\sqrt{n}-\|{\bd w}_j\|)^2\|\bd{M}{\bm \gamma}_j\|^2 + 2
\sum_{j=1}^q\|\bd{M}\bd{\Sigma}_{j-1}\bd{y}_j\big\|^2\\
\endaligned
\ee
by the triangle inequality and the formula $(a+b)^2\leq 2a^2+ 2b^2$ for any $a,b\in \mathbb{R}.$ Define
\beaa
A_j=\big(\sqrt{n}-\|\bd{w}_j\|\big)^2;\ \ B_j=\|\bd{M}\bm{\gamma}_j\|^2;\ \
C_j=\|\bd{M}\bd{\Sigma}_{j-1}\bd{y}_j\big\|^2
\eeaa
for $1\leq j \leq n$ with $C_0=0$. Then,
\bea\lbl{simple}
\mathbb{E}\|\sqrt{n}\bd{\Gamma}_{p\times q}-\bd{Y}_{p\times q}\|^2_{\rm HS} \leq 2\sum_{j=1}^q \mathbb{E}(A_jB_j) + 2\sum_{j=1}^q \mathbb{E}C_j.
\eea
We next bound $A_j$, $B_j$ and $C_j$, respectively, in  terms of their moments.

{\it The estimate of $A_j$}. Trivially,
$$\aligned
A_j^2=\Big[\frac{\|\bd{w}_j\|^2-n}{(\|\bd{w}_j\|+\sqrt{n})}\Big]^4
\leq\frac{1}{n^2}\cdot(\|\bd{w}_j\|^2-n)^4. \endaligned $$
By Lemma 3.1, $\|\bd{w}_j\|^2\sim\chi^2(n-j+1).$ Set $c_j:=n-j+1$ and ${\bf z}_j=\|{\bf w}_j\|^2-c_j.$ By Lemma \ref{party_love} and the binomial formula, we have
$$\aligned
n^2\mathbb{E}(A_j^2) &\le \mathbb{E}\big[(\bd{z}_j-j+1)^4\big]\\
&=\mathbb{E}(\bd{z}_j^4)-4(j-1)\mathbb{E}(\bd{z}_j^3)+6(j-1)^2\mathbb{E}(\bd{z}_j^2)-4(j-1)^3\mathbb{E}\bd{z}_j+(j-1)^4\\
&=12c_j(c_j+4)-32(j-1)c_j+12c_j(j-1)^2+(j-1)^4\\
&\le 12(c_j+2)^2+12(c_j+2)(j-1)^2+3(j-1)^4\\
&=12\big(c_j+2+\frac{(j-1)^2}{2}\big)^2.
\endaligned $$
This immediately implies that
$$(\mathbb{E}A_j^2)^{1/2}\le\frac{2\sqrt{3}}{n}\Big[n-j+3+\frac{(j-1)^2}{2}\Big].$$

{\it The estimate of $B_j$}.  Recall (\ref{sea}). The vector $\bm{\gamma}_j=\frac{\bd{w}_j}{\|\bd{w}_j\|}$ has the same distribution as $\bm{\gamma}_1=\frac{\bd{y}_1}{\|\bd{y}_1\|}$.
Note $\bd{y}_1=(y_{11}, \cdots, y_{n1})'$, hence  $\|\bd{M}\bd{\gamma}_1\|^2=U_1+\cdots + U_p$ where $U_i=y_{i1}^2/(y_{11}^2 +\cdots y_{n1}^2)$ for $1\leq i \leq n$. By Lemma \ref{Jiang2009},
\beaa
\mathbb{E}\big(\|\bd{M}\bd{\gamma}_1\|^4\big) & = & \mathbb{E}\big[(U_1+\cdots + U_p)^2\big] \\
& = & p\mathbb{E}(U_1^2)+ p(p-1)\mathbb{E}(U_1U_2)\\
& = &\frac{p(p+2)}{n(n+2)} .
\eeaa
Therefore $$\mathbb{E}B_j^2=\mathbb{E}\big(\|\bd{M}\bd{\gamma}_1\|^4\big)=\frac{p(p+2)}{n(n+2)}\le \frac{3p^2}{n^2}.$$
In particular, the two estimates above  conclude that
$$ |\mathbb{E}(A_jB_j)|\leq
(\mathbb{E}A_j^2)^{1/2}(\mathbb{E}B_j^2)^{1/2}\leq \frac{6p}{n^2}\Big(n-j+3+\frac{(j-1)^2}{2}\Big),\\
$$
which guarantees
\be\lbl{ABsum}\aligned  2\sum_{j=1}^q\mathbb{E}(A_jB_j)&\le \frac{12p}{n^2}\sum_{j=1}^q\Big(n-j+3+\frac{(j-1)^2}{2}\Big) \\
&=\frac{12pq}{n}\Big(1+\frac{2q^2-9q+31}{12n}\Big).\\
\endaligned
\ee

{\it The estimate of $C_j$.}  Now, conditioning on $\bd{y}_1, \cdots, \bd{y}_{j-1}$, we get from (\ref{normal_fact}) that
\be\lbl{Cespress}
C_j=\|\bd{M}\bd{\Sigma}_{j-1}\bd{y}_j\big\|^2
= \bd{y}_j'\bd{\Sigma}_{j-1}
\begin{pmatrix}
\bd{I}_p & \bd{0}\\
\bd{0} & \bd{0}
\end{pmatrix}
\bd{\Sigma}_{j-1}\bd{y}_j  \overset{d}{=}  \sum_{k=1}^n\lambda_k\xi_k^2, \ee
where $\xi_1, \cdots, \xi_n$ are i.i.d. $N(0,1)$-distributed random variables and $\lambda_1, \cdots, \lambda_n$ are the eigenvalues of
\be\lbl{defA}
\bd{A}_{j-1}:=\bd{\Sigma}_{j-1}
\begin{pmatrix}
\bd{I}_p & \bd{0}\\
\bd{0} & \bd{0}
\end{pmatrix}
\bd{\Sigma}_{j-1}.
\ee
In particular,
\beaa
\mathbb{E}C_j = \mathbb{E}\,\mbox{tr}(\bd{A}_{j-1}).
\eeaa
By the fact $\mbox{tr}(\bd{A}\bd{B})=\mbox{tr}(\bd{B}\bd{A})$ for any matrix $\bd{A}, \bd{B}$ and the fact that both $\bd{\Sigma}_{j-1}$ and $\begin{pmatrix}
\bd{I}_p & \bd{0}\\
\bd{0} & \bd{0}
\end{pmatrix}$ are idempotent, we see
\be\lbl{traceA}\aligned
\mbox{tr}(\bd{A}_{j-1})&=\mbox{tr}\, \Big[\bd{\Sigma}_{j-1}
\begin{pmatrix}
\bd{I}_p & \bd{0}\\
\bd{0} & \bd{0}
\end{pmatrix}
\begin{pmatrix}
\bd{I}_p & \bd{0}\\
\bd{0} & \bd{0}
\end{pmatrix}
\bd{\Sigma}_{j-1}\Big]\\
& = \mbox{tr}\, \Big[
\begin{pmatrix}
\bd{I}_p & \bd{0}\\
\bd{0} & \bd{0}
\end{pmatrix}
\bd{\Sigma}_{j-1}\begin{pmatrix}
\bd{I}_p & \bd{0}\\
\bd{0} & \bd{0}
\end{pmatrix}\Big]\\
&={\rm tr}\big((\bd{\Sigma}_{j-1})_p\big),
\endaligned
\ee
where $(\bd{\Sigma}_{j-1})_p$ is as in the statement of Lemma \ref{Togo}. Hence by \eqref{monkey}, we have
\be\lbl{Cesti}
\sum_{j=2}^q\mathbb{E}C_j=\sum_{j=2}^q\mathbb{E}\,\mbox{tr}(\bd{A}_{j-1})
=\sum_{j=2}^q\mathbb{E}\,\mbox{tr}((\bd{\Sigma}_{j-1})_p)=\frac{pq(q-1)}{2n}.\ee
Therefore plugging \eqref{ABsum} and \eqref{Cesti} into \eqref{simple}, we know
\beaa
\mathbb{E}\|\sqrt{n}\bd{\Gamma}_{p\times q}-\bd{Y}_{p\times q}\|^2_{\rm HS}
&\leq & 2\sum_{j=1}^q\mathbb{E}(A_jB_j) + 2\sum_{j=1}^q\mathbb{E}C_j \nonumber\\
&\le &\frac{12pq}{n}\Big(1+\frac{2q^2-9q+31}{12n}\Big)+\frac{pq(q-1)}{n}\\
&= &\frac{pq^2}{n}\cdot\Big(1+\frac{2q^2-9q+31+11n}{nq}\Big).\\
\eeaa
Define $f(q)=2q+\frac{31+11n}{q}-9$ for $q>0.$ Since $f{''}(q)>0$ for all $q>0$, we know $f(q)$ is a convex function. Therefore, $\max_{1\leq q \leq n}f(q)= f(1)\vee f(n).$
Trivially, $$f(1)-f(n)=24+11n-\Big(2n+\frac{31}{n}+2\Big)=22+9n-\frac{31}{n}\ge 22+9-31=0$$ for all $n\geq 1$.  We then have
$\max_{1\leq q \leq n}f(q)=24+11n$
for any $n\geq 2.$
Thus,
$$
\mathbb{E}\|\sqrt{n}\bd{\Gamma}_{p\times q}-\bd{Y}_{p\times q}\|^2_{\rm HS}
\leq\frac{pq^2}{n}\Big(1+\frac{24+11n}{n}\Big)=
\frac{12pq^2}{n}\Big(1+\frac2n\Big)\le\frac{24pq^2}{n}
$$
for any $n\geq 2.$ The proof is completed.  \hfill$\blacksquare$

\medskip
\medskip
Similar to  Lemma \ref{High_teacher}, the next result will serve as the framework of the proof of  Theorem \ref{ttt}. The spirit of the proof is close to that of Lemma \ref{High_teacher}. We therefore omit it. For a sequence of numbers $\{a_n;\, n\geq 1\}$, we write $\lim_{n\to\infty}a_n\in (0, \infty)$ if $\lim_{n\to\infty}a_n=a$ exists and $a\in (0, \infty)$.
\begin{lem}\lbl{Low_beauty}
For each $n\geq 1$, let $f_n(p,q):\{1, 2,\cdots, n\}^2\to [0, \infty)$ satisfy that  $f_n(p, q)$ is non-decreasing in  $p\in \{1, 2,\cdots, n\}$ and $q\in \{1, 2,\cdots, n\}$, respectively. Suppose
\bea\lbl{red_red}
\liminf_{n\to\infty}f_n(p_n, q_n)>0
\eea
for any sequence $\{(p_n, q_n)\in \{1, 2,\cdots, n\}^2\}_{n=1}^{\infty}$  if any of the next two conditions holds:

(i) $q_n\equiv 1$ and  $\lim_{n\to\infty}p_n/n\in (0,1)$;


(ii) $\lim_{n\to\infty}q_n=\infty$ and  $\lim_{n\to\infty}(p_nq_n^2)/n\in (0, \infty)$.

\noindent Then (\ref{red_red}) holds for any sequence $\{(p_n, q_n)\}_{n=1}^{\infty}$ with $1\leq p_n, q_n\leq n$ for each $n\geq 1$ and  $\lim_{n\to\infty}(p_nq_n^2)/n\in (0, \infty)$.
\end{lem}

\subsection{The Proof of Theorem \ref{ttt}}\lbl{big_mouse}

After many pieces of understanding, we are now ready to prove the second main result in this paper.

\medskip

\noindent\textbf{Proof of Theorem \ref{ttt}}. The first part follows immediately from Proposition \ref{lala}. The second part is given next.

We first prove that
\bea\lbl{God_creat_math}
\|\sqrt{n}\bd{\Gamma}_{p\times q}-\bd{Y}_{p\times q}\|_{\rm HS} \overset{p}{\to} \Big(\frac{\sigma}{2}\Big)^{1/2}
\eea
for any $1\leq p_n, q_n \leq n$ satisfying  $q_n\to \infty$ and $\frac{pq^2}{n}\to  \sigma>0$. We claim this  implies that
\bea\lbl{sun_cars}
\liminf_{n\to\infty}P(\|\sqrt{n}\bd{\Gamma}_{p\times q}-\bd{Y}_{p\times q}\|_{\rm HS}\geq \epsilon)>0
\eea
for any $\epsilon \in (0, \sqrt{\sigma/2})$ and  any $1\leq p_n, q_n \leq n$ with  $\frac{pq^2}{n}\to  \sigma>0$. In fact, for given $\epsilon \in (0, \sqrt{\sigma/2})$, set
\beaa
f_n(p,q)=P(\|\sqrt{n}\bd{\Gamma}_{p\times q}-\bd{Y}_{p\times q}\|_{\rm HS}\geq \epsilon)
\eeaa
 for all $1\leq p, q\leq n$. Here we slightly abuse some notation: $\bd{\Gamma}_{p\times q}$ and $\bd{Y}_{p\times q}$ are $p\times q$ matrices with $p$ and $q$ being arbitrary instead of  fixed sizes $p_n$ and $q_n.$ By (\ref{Hiha}), it is obvious that $f_n(p,q)$ is non-decreasing in $p$ and $q$, respectively, for any $n\geq 1$. Assume (\ref{God_creat_math}) holds, then $\liminf_{n\to\infty}f_n(p_n, q_n)=1$ for any $1\leq p_n, q_n \leq n$ under condition $\lim_{n\to\infty}q_n= \infty$ and $\lim_{n\to\infty}\frac{pq^2}{n}=\sigma \in (0, \infty).$ By Lemma \ref{piano_black},  $\liminf_{n\to\infty}f_n(p_n, 1)=P(|N(0,1)|\geq \epsilon\sqrt{2/c})$ for any $1\leq p_n \leq n$ with $\lim_{n\to\infty}p_n/n=c\in (0,1)$. Then we obtain (\ref{sun_cars}) from Lemma \ref{Low_beauty}.

Now we start to prove (\ref{God_creat_math}).  Let us continue to use the notation in the proof of Proposition \ref{lala}. Review $\bd{M}=\bd{M}_{p\times n}=(\bd{I}_p, \bd{0})$ for $1\leq p \leq n-1$ and $\bd{M}_{n\times n}=\bd{I}_n$, where $\bd{I}_p$ is the $p\times p$ identity matrix and $\bd{0}$ is the $p\times (n-p)$ matrix whose entries are all equal to zero. By \eqref{expnorm},
\beaa
\|\sqrt{n}\bd{\Gamma}_{p\times q}-\bd{Y}_{p\times q}\|_{\rm HS}^2=
 \sum_{j=1}^q\|(\sqrt{n}-\|\bd{w}_j\|)\bd{M}\bm{\gamma}_j-
\bd{M}\bd{\Sigma}_{j-1}\bd{y}_j\big\|^2.
\eeaa
Review
$$
A_j=(\sqrt{n}-\|\bd{w}_j\|)^2;\ \ B_j=\|\bd{M}\bm{\gamma}_j\|^2;\ \ C_j=\|\bd{M}\bd{\Sigma}_{j-1}\bd{y}_j\big\|^2.
$$
For vectors $\bd{u}, \bd{v}\in \mathbb{R}^p$, we know $\|\bd{u}+\bd{v}\|^2=\|\bd{u}\|^2+\|\bd{v}\|^2 + 2\langle\bd{u}, \bd{v}\rangle$. By the Cauchy-Schwartz inequality, $|\langle\bd{u}, \bd{v}\rangle| \leq \|\bd{u}\|\cdot \|\bd{v}\|.$ So we can write
$$
\|\sqrt{n}\bd{\Gamma}_{p\times q}-\bd{Y}_{p\times q}\|_{\rm HS}^2=\sum_{j=1}^qA_jB_j + \sum_{j=1}^qC_j + \sum_{j=1}^q\epsilon_j
$$
where $|\epsilon_j|\leq 2\sqrt{A_jB_jC_j}$ for $1\leq j \leq q.$ From (\ref{ABsum}), we see that
$$
\mathbb{E}\sum_{j=1}^qA_jB_j
\le \frac{6pq}{n}\Big(1+\frac{2q^2-9q+31}{12n}\Big)
=O\Big(\frac{1}{q}+\frac1{pq}\Big)\to 0
$$
since $q=q_n\to\infty$ and $\frac{pq^2}{n}\to  \sigma>0$ as $n\to\infty.$ In particular,
\bea\lbl{bad_child}
\sum_{j=1}^qA_jB_j \overset{p}{\to} 0
\eea
as $n\to\infty$. We claim that it suffices to show
\bea\lbl{last}
\sum_{j=1}^qC_j \overset{p}{\to} \frac{\sigma}{2}
\eea
as $n\to\infty$. In fact, once \eqref{last} holds, we have by the Cauchy-Schwartz inequality and \eqref{bad_child}
$$
(\sum_{j=1}^q|\epsilon_j|)^2 \le4(\sum_{j=1}^q\sqrt{A_jB_jC_j} )^2\le 4\Big(\sum_{j=1}^qA_jB_j\Big)\sum_{j=1}^q C_j\stackrel{p}{\to} 0
$$
as $n\to\infty.$ Then $\sum_{j=1}^q\epsilon_j\stackrel{p}{\to}0 $ as $n\to\infty.$ Now we prove \eqref{last}.

Recall the notation $(\bd{\Sigma}_{k})_{p}$ stands for the upper-left $p\times p$ submatrix of $\bd{\Sigma}_{k}$. Let $\mathcal{F}_{k}$ be
 the sigma-algebra generated by $\bd{y}_1, \bd{y}_2, \cdots, \bd{y}_k.$ We first claim
\bea\lbl{conclusion2}
X_j:=C_j-{\rm tr }[(\bd{\Sigma}_{j-1})_{p}],\ \ j=2,3,\cdots,q,
\eea
forms  a martingale difference with respect to $\mathcal{F}_2,\mathcal{F}_3, \cdots, \mathcal{F}_{q}.$ In fact, as in \eqref{Cespress}, we write
\bea\lbl{no_need}
C_j=\bd{y}_j'\bd{A}_{j-1}\bd{y}_j
\eea
for any $2\le j\le q,$ where the symmetric matrix
$
\bd{A}_{j-1}$
is defined in \eqref{defA} and is independent of $\bd{y}_j$. Let $\mu_1, \cdots, \mu_n$ be the eigenvalues of $\bd{A}_{j-1}$. By (\ref{normal_fact}), \eqref{traceA} and independence,
 \beaa\lbl{no_useless}
 \mathbb{E}(\bd{y}_j' \bd{A}_{j-1} \bd{y}_j|\mathcal{F}_{j-1})=\mathbb{E}\sum_{i=1}^n\mu_i\xi_i^2=\mbox{tr}\,(\bd{A}_{j-1})
 =\mbox{tr}\,[(\bd{\Sigma}_{j-1})_p],
 \eeaa
 where $\xi_1, \cdots, \xi_n$ are i.i.d. standard normals.
This confirms (\ref{conclusion2}).

Obviously,  $$B_k:=\sum_{j=1}^k X_j,\ \ k=2,\cdots, q$$ is a martingale relative to $\{\mathcal{F}_k;\, k=2,\cdots, q\}$. Therefore,
$$\sum_{j=2}^qC_j=B_q+\sum_{j=2}^q{\rm tr}[(\bd{\Sigma}_{j-1})_p].$$
By Lemma \ref{trofSigma}, when $pq^2/n\to \sigma ,$ $$\sum_{j=2}^q{\rm tr}[(\bd{\Sigma}_{j-1})_p]\stackrel{p}{\longrightarrow}\frac{\sigma}2$$ as $n\to\infty.$ To get (\ref{last}), it is enough to show
\bea\lbl{inspiration}
{\rm Var}(B_q)\to 0
\eea
as $n\to\infty.$
Since $(X_i)_{1\le i\le q}$ is a martingale difference, it entails that
$$\mathbb{E}(X_iX_j)=\mathbb{E}[X_i\mathbb{E}(X_j|\mathcal{F}_{j-1})]=0$$ for any $2\le i<j\le q.$  Also, recall the conditional variance has the formula
\beaa
{\rm Var}(X_j)&=& \mathbb{E}\, {\rm Var}\,(X_j|\mathcal{F}_{j-1}) + {\rm Var}\,[\mathbb{E} (X_j|\mathcal{F}_{j-1})]\\
& = & \mathbb{E}\, {\rm Var}\,(X_j|\mathcal{F}_{j-1})
\eeaa
since $X_j$ is a martingale difference,
where ${\rm Var}\,(X_j|\mathcal{F}_{j-1})=
\mathbb{E}\big[(X_j-\mathbb{E}(X_j|\mathcal{F}_{j-1}))^2|\mathcal{F}_{j-1}\big]$; see, for example, \cite{Casella}.
Therefore, by (\ref{no_need}) and then (\ref{variance2})
\beaa
 {\rm Var}(B_q)=\sum_{j=2}^q {\rm Var}(X_j)&=&\sum_{j=2}^q\mathbb{E}\, {\rm Var}(\bd{y}_j'\bd{A}_{j-1} \bd{y}_j|\mathcal{F}_{j-1})\\
&=& 2\sum_{j=2}^q \mathbb{E}\,{\rm tr}(\bd{A}_{j-1}^2).
\eeaa
Repeatedly using the facts
\beaa
\bd{\Sigma}_{j-1}^2=\bd{\Sigma}_{j-1},\ \
\begin{pmatrix}\bd{I}_{p} & \bd{0} \\
\bd{0} & \bd{0} \end{pmatrix}^2=\begin{pmatrix}\bd{I}_{p} & \bd{0} \\
\bd{0} & \bd{0} \end{pmatrix}, \ \ \ \mbox{tr}\, (\bd{U}\bd{V})=\mbox{tr}\, (\bd{V}\bd{U})
\eeaa
for any $n\times n$ matrices $\bd{U}$ and $\bd{V}$,
it is not difficult to see ${\rm tr}(\bd{A}_{j-1}^2)={\rm tr}[(\bd{\Sigma}_{j-1})_p^2]$. From
Lemma \ref{trofSigma},
\beaa
\sum_{j=2}^q \mathbb{E}\,{\rm tr}(\bd{A}_{j-1}^2)&= &\frac{pq(q-1)(p+2)}{2n(n+2)}+\frac{pq(q-1)(q-2)(n-p)}{3n(n-1)(n+2)}\\
& \leq & C\cdot \Big(\frac{p^2q^2}{n^2} + \frac{pq^3}{n^2}\Big)=O(\frac1{q^2}+\frac1{pq})\to 0
\eeaa
as $n\to\infty$ by the assumption $q\to\infty$ and $\frac{pq^2}{n}\to\sigma>0.$  We gets (\ref{inspiration}). The proof is completed. \hfill$\blacksquare$

\section{Appendix}\lbl{Good_Tue}

In this section we will prove Lemmas \ref{great}, \ref{awesome_pen} and \ref{var-e}. We start with Lemma \ref{great}, which computes the mean values of monomials of the matrix elements
from an Haar-orthogonal matrix.

\begin{figure}
\begin{center}
\includegraphics[height=3.5in,width=6.4in]{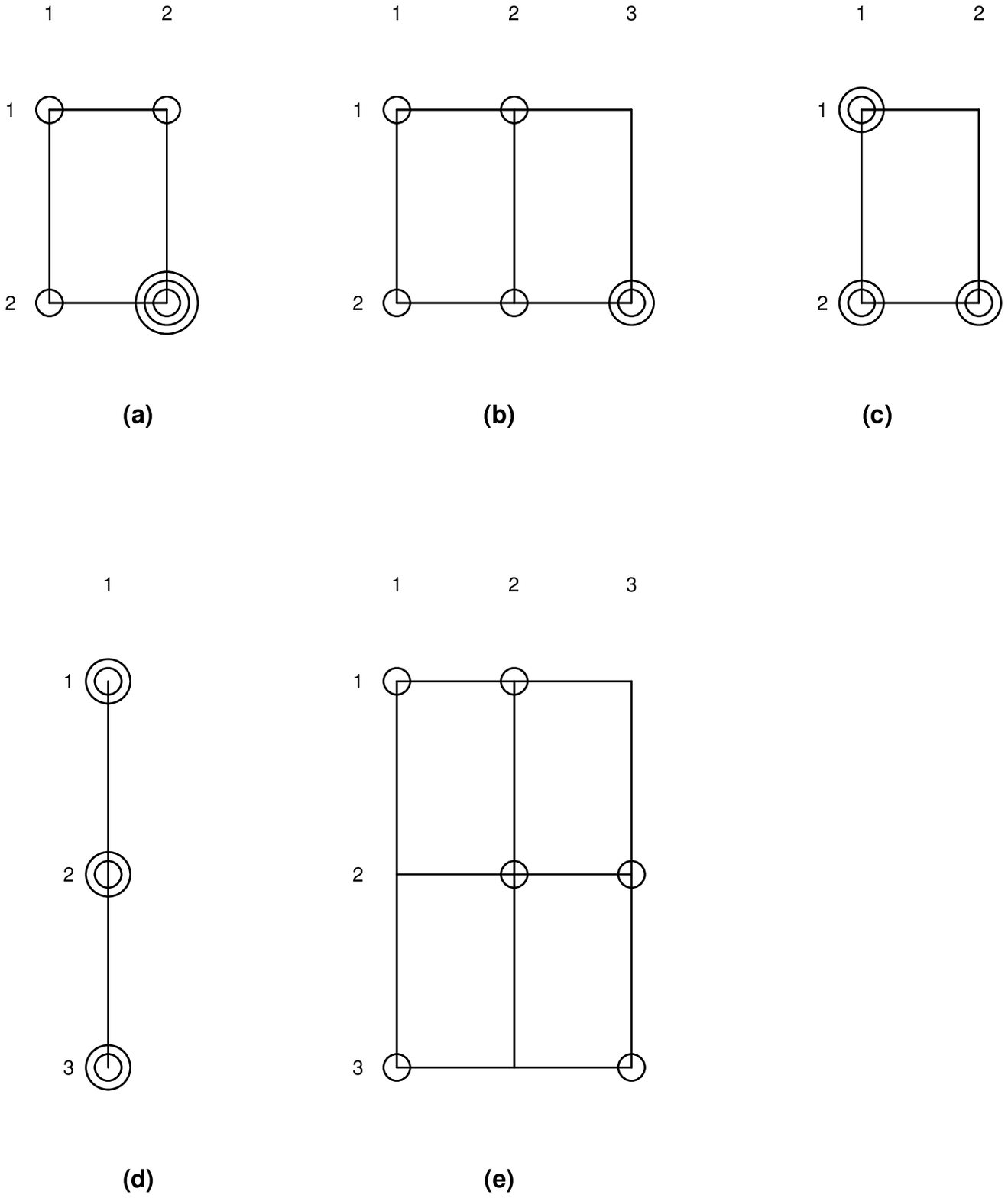}
\caption{}
\end{center}
\label{fig200}
\end{figure}

To make the monomials more intuitive, we make Figure \ref{fig200}. For each plot inside the graph,
the number of circles  appearing in a corner means the power of the corresponding matrix entry appearing in the monomial. For example, plot (d)
stands for the monomial $\gamma_{11}\gamma_{12}\gamma_{21}\gamma_{22}^3$;
plot (e) represents $\gamma_{11}\gamma_{12}\gamma_{22}\gamma_{23}\gamma_{31}\gamma_{33}.$

\medskip

\noindent\textbf{Proof of Lemma \ref{great}}. Our argument below are based on the unit length of
each row/column, the orthogonality of any two rows/columns and that all row/column random
vectors are exchangeable. We will first prove conclusions (a) and (c), and then prove the rest of them.

(a) Recall {\bf F1}).
Take $a_1=a_2=a_3=1$ and $a_3=\cdots =a_n=0$ in Lemma \ref{Jiang2009}, we get the conclusion.

(c) Since $$\gamma_{11}^2\gamma_{21}^2(\sum_{j=1}^n\gamma_{1j}^2-1)=0,$$ take expectations to get
\beaa
\mathbb{E}(\gamma_{11}^4\gamma_{21}^2) + (n-1)\mathbb{E}(\gamma_{11}^2\gamma_{21}^2\gamma_{12}^2) =\mathbb{E}(\gamma_{11}^2\gamma_{21}^2).
\eeaa
By using Lemma \ref{Jiang2009}, we see
\bea\lbl{rabbit}
\mathbb{E}(\gamma_{11}^2\gamma_{21}^2)=\frac{1}{n(n+2)};\ \ \ \mathbb{E}(\gamma_{11}^4\gamma_{21}^2)=\frac{3}{n(n+2)(n+4)}.
\eea
The conclusion (c) follows.

(b), (d) \& (e). Since the first and the second columns of $\f {\Gamma}_n$ are mutually orthogonal, we know
\bea\lbl{tri-1}\aligned
0&=\mathbb{E}\Big(\sum_{i=1}^n \gamma_{i1}\gamma_{i2}\sum_{j=1}^n \gamma_{j2}\gamma_{j3} \sum_{k=1}^n \gamma_{k1}\gamma_{k3}\Big)\\
&=n\mathbb{E}(\gamma_{11}^2\gamma_{12}^2\gamma_{13}^2)+
3n(n-1)\mathbb{E}(\gamma_{11}^2\gamma_{12}\gamma_{13}\gamma_{22}\gamma_{23})\\
&+n(n-1)(n-2)\mathbb{E}(\gamma_{11}\gamma_{12}\gamma_{22}\gamma_{23}\gamma_{31}\gamma_{33}).
\endaligned \eea
Similarly we have
\bea\lbl{tri-2}\aligned
0&=\mathbb{E}\Big(\sum_{i=1}^n \gamma_{i1}^2\gamma_{i2}\gamma_{i3}\sum_{j=1}^n \gamma_{j2} \gamma_{j3}\Big)\\
&=n\mathbb{E}(\gamma_{11}^2\gamma_{12}^2\gamma_{13}^2)+
n(n-1)\mathbb{E}(\gamma_{11}^2\gamma_{12}\gamma_{13}\gamma_{22}\gamma_{23})\\
\endaligned \eea
and
\bea\lbl{tri-3}\aligned
0&=\mathbb{E}\Big(\sum_{i=1}^n \gamma_{i1}^3\gamma_{i2}\sum_{j=1}^n \gamma_{j1} \gamma_{j2}\Big)\\
&=n\mathbb{E}(\gamma_{11}^4\gamma_{12}^2)+n(n-1)\mathbb{E}(\gamma_{11}^3\gamma_{12}\gamma_{21}\gamma_{22}).\\
\endaligned \eea
Combining the expressions \eqref{tri-1}, \eqref{tri-2} and \eqref{tri-3} together with  conclusion (a) and (\ref{rabbit}),   we arrive at
\beaa\lbl{tri-t}\aligned
&\mathbb{E}(\gamma_{11}^2\gamma_{12}\gamma_{13}\gamma_{22}\gamma_{23})=
-\frac{1}{(n-1)n(n+2)(n+4)}; \\
& \mathbb{E}( \gamma_{11}^3\gamma_{12}\gamma_{22}\gamma_{21})=-\frac{3}{(n-1)n(n+2)(n+4)};\\
& \mathbb{E}(\gamma_{11}\gamma_{12}\gamma_{22}\gamma_{23}\gamma_{31}\gamma_{33})
=\frac{2}{(n-2)(n-1)n(n+2)(n+4)}.
\endaligned
\eeaa
By swapping rows and columns and using the invariance, we get
\beaa
& & \mathbb{E}(\gamma_{11}^2\gamma_{12}\gamma_{13}\gamma_{22}\gamma_{23})=
\mathbb{E}(\gamma_{11}\gamma_{12}\gamma_{21}\gamma_{22}\gamma_{23}^2);\\
& & \mathbb{E}(\gamma_{11}\gamma_{12}\gamma_{21}\gamma_{22}^3)=\mathbb{E}( \gamma_{11}^3\gamma_{12}\gamma_{22}\gamma_{21}).
\eeaa
The proof is completed. \hfill$\blacksquare$

\medskip

We will derive the central limit theorem appearing in Lemma \ref{awesome_pen} next. Two preliminary calculations are needed.

\blem\lbl{doublegreat}
Let $\bd{y}=(\xi_1, \cdots, \xi_p)'\sim N_p(\bd{0}, \bd{I}_p)$. Let $\bd{a}=(\alpha_1, \cdots, \alpha_p)'\in \mathbb{R}^p$
and $\bd{b}=(\beta_1, \cdots, \beta_p)'\in \mathbb{R}^p$ with $\|\bd{a}\|=\|\bd{b}\|=1$.  Then,
$$\mathbb{E}\big[(\bd{a}' \bd{y})^2(\bd{b}'\bd{y})^2\big]=2(\bd{a}'\bd{b})^2+1.$$
\nlem
\noindent\textbf{Proof}. Rewrite
$$\bd{a}'\bd{y}\bd{b}'\bd{y}=\bd{y}' \bd{a} \bd{b}' \bd{y}=\bd{y}'\bd{b}\bd{a}'\bd{y}=\frac12 \bd{y}'(\bd{a}\bd{b}'+\bd{b}\bd{a}')\bd{y}.$$  Define
$\bd{A}=\frac12(\bd{a}\bd{b}'+\bd{b}\bd{a}').$ Then, $\bd{a}'\bd{y}\bd{b}'\bd{y}=\bd{y}'
\bd{A}\bd{y}$. By (\ref{normal_fact}) and (\ref{variance2}),
\beaa
\mathbb{E}(\bd{y}'\bd{A}\bd{y})=\,\mbox{tr}(\bd{A})\ \
\mbox{and}\ \ \mbox{Var}(\bd{y}'\bd{A}\bd{y})=2\, \mbox{tr} (\bd{A}^2).
\eeaa
Note that ${\rm tr}(\bd{A})=\bd{a}' \bd{b}.$  Therefore,
\beaa
\mathbb{E}\big[(\bd{a}' \bd{y})^2(\bd{b}'\bd{y})^2\big]
&= &2{\rm tr}(\bd{A}^2)+{\rm tr}^2(\bd{A})\\
&=&\frac12{\rm tr}\big(\bd{a}\bd{b}'\bd{a}\bd{b}'+\bd{b}\bd{a}'\bd{b}\bd{a}'+
\bd{a}\bd{b}'\bd{b}\bd{a}'+\bd{b}\bd{a}'\bd{a}\bd{b}'\big)+(\bd{a}'\bd{b})^2\\
&=& 2(\bd{a}'\bd{b})^2+1
\eeaa
by the assumption $\|\bd{a}\|=\|\bd{b}\|=1$.
\hfill$\blacksquare$

\blem\lbl{soup_good} Let $\bd{u}=(\xi_1, \cdots, \xi_p)'$ and $\bd{v}=(\eta_1, \cdots, \eta_p)'$ be independent random vectors with distribution $N_p(\bd{0}, \bd{I}_p).$ Set $\bd{w}=(\bd{u}'\bd{v})^2-\|\bd{u}\|^2$. Then
\bea
& & \mathbb{E}[\bd{w} | \bd{v}]
=\|\bd{v}\|^2-p; \lbl{Clinton}\\
& & \mathbb{E}[\bd{w}^2 | \bd{u}]
=2\|\bd{u}\|^4; \lbl{Trump}\\
& \mathbb{E}[\bd{w}^2| \bd{v}]&=3\|\bd{v}\|^4+(p^2+2p)-2(p+2)\|\bd{v}\|^2. \lbl{none}
\eea
\nlem
\noindent\textbf{Proof}. The assertion (\ref{Clinton}) follows from independence directly. Further,
\beaa
\mathbb{E}\big[\bd{w}^2 | \bd{u}\big]&=&
\mathbb{E}\big[(\bd{u}'\bd{v})^4-2\|\bd{u}\|^2 (\bd{u}'\bd{v})^2+\|\bd{u}\|^4 | \bd{u}\big]\\
& = & 3\|\bd{u}\|^4-2\|\bd{u}\|^4+\|\bd{u}\|^4=2 \|\bd{u}\|^4.
\eeaa
We then obtain (\ref{Trump}). Finally, since $\|\bd{u}\|^2 \sim \chi^2(p)$,
we have
 \beaa
  \mathbb{E}[\bd{w}^2| \bd{v}]&=&\mathbb{E}\big[\big((\bd{u}'\bd{v})^4+\|\bd{u}\|^4-2\|\bd{u}\|^2
  (\bd{u}'\bd{v})^2\big)|\bd{v}\big]\\
&=& 3\|\bd{v}\|^4+(p^2+2p)-2\, \mathbb{E}\Big[\|\bd{u}\|^2
\big(\sum_{l=1}^p\xi_l\eta_l\big)^2\big| \bd{v}\Big].
\eeaa
Expanding the last sum, we see from independence that
\beaa
\mathbb{E}\Big[\|\bd{u}\|^2
\big(\sum_{l=1}^p\xi_l\eta_l\big)^2\big| \bd{v}\Big] &=& \sum_{l=1}^p\eta_l^2\mathbb{E}\big[\|\bd{u}\|^2\xi_1^2\big] + \sum_{1\leq k\ne l \leq p}\eta_k\eta_l\mathbb{E}\big[\|\bd{u}\|^2\xi_1\xi_2\big]\\
& = & \Big(\sum_{l=1}^p\eta_l^2\Big)\cdot \frac{1}{p}\Big(\mathbb{E}\sum_{l=1}^p\|\bd{u}\|^2\xi_l^2\Big)\\
& = & \|\bd{v}\|^2\cdot \frac{1}{p}\mathbb{E}\big(\|\bd{u}\|^4\big)=(p+2)\|\bd{v}\|^2
\eeaa
by the fact $\mathbb{E}\big[\|\bd{u}\|^2\xi_1\xi_2\big]=0$ due to the symmetry of normals random variables. The above two identities imply (\ref{none}). \hfill$\blacksquare$

\medskip

\begin{figure}
\title{Plot of CLT}
\begin{center}
\begin{tabular}{ccc}
\includegraphics[height=1.5in,width=1.7in]{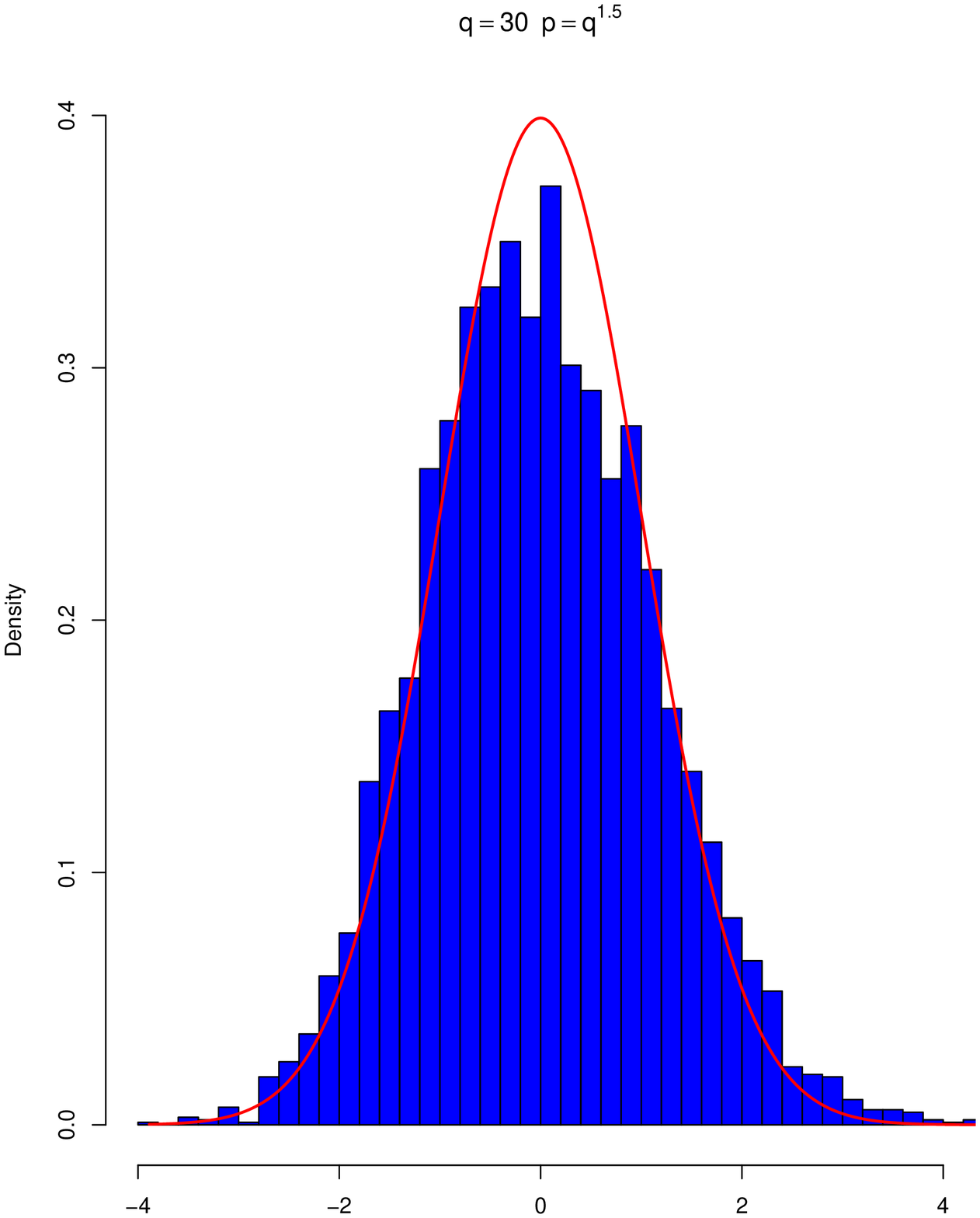}&\includegraphics[height=1.5in,width=1.7in]{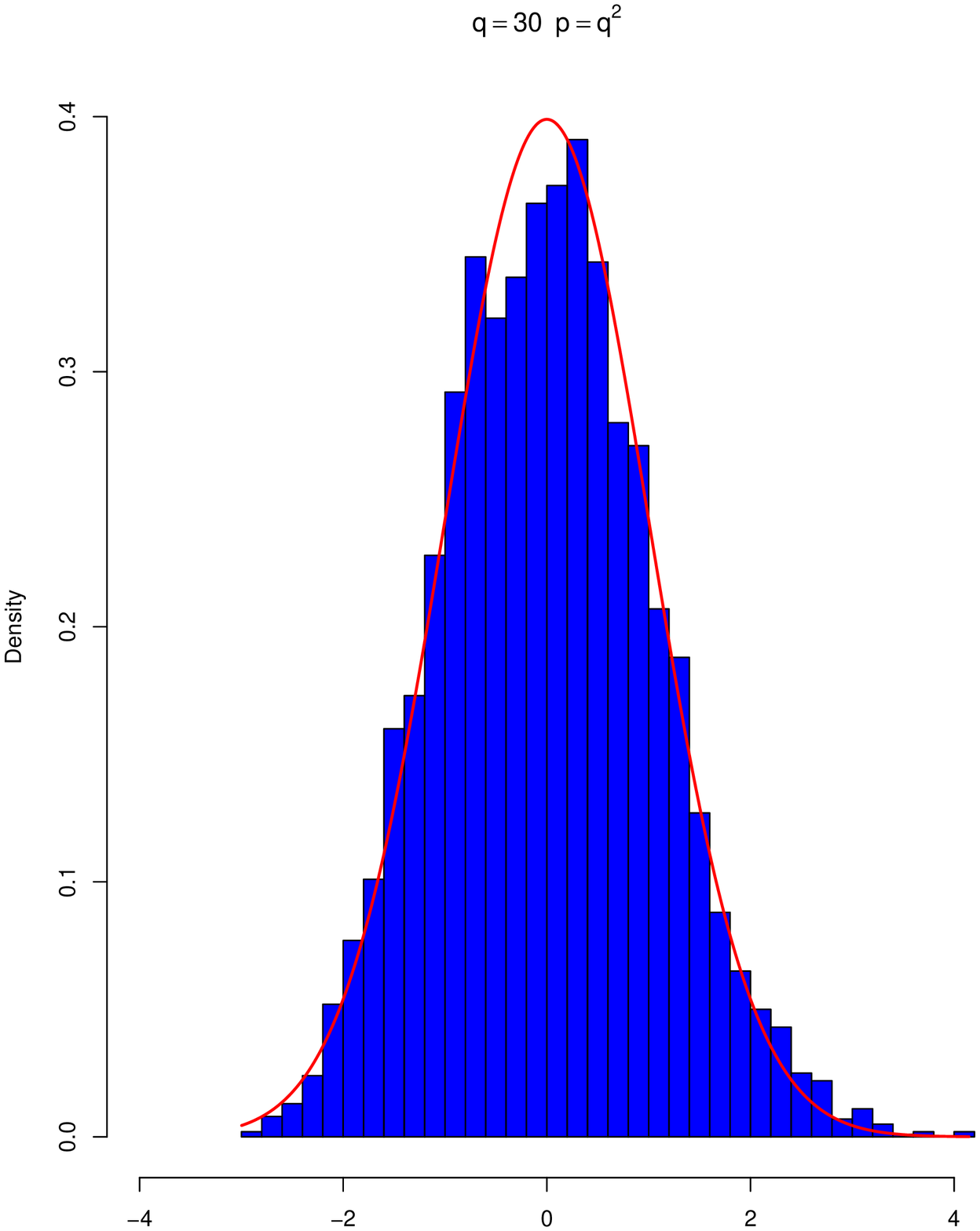} & \includegraphics[height=1.5in,width=1.7in]{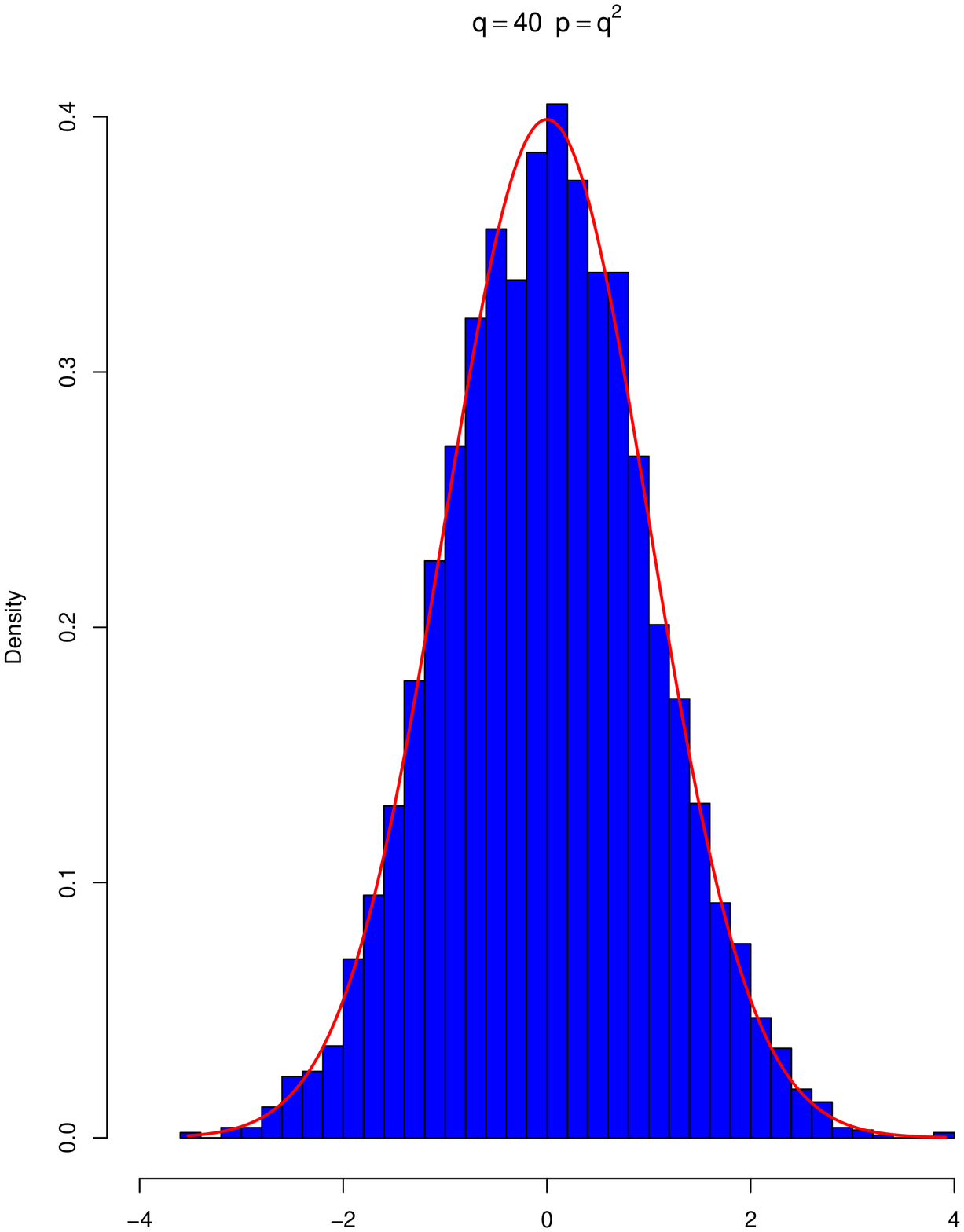} \\
\includegraphics[height=1.5in,width=1.7in]{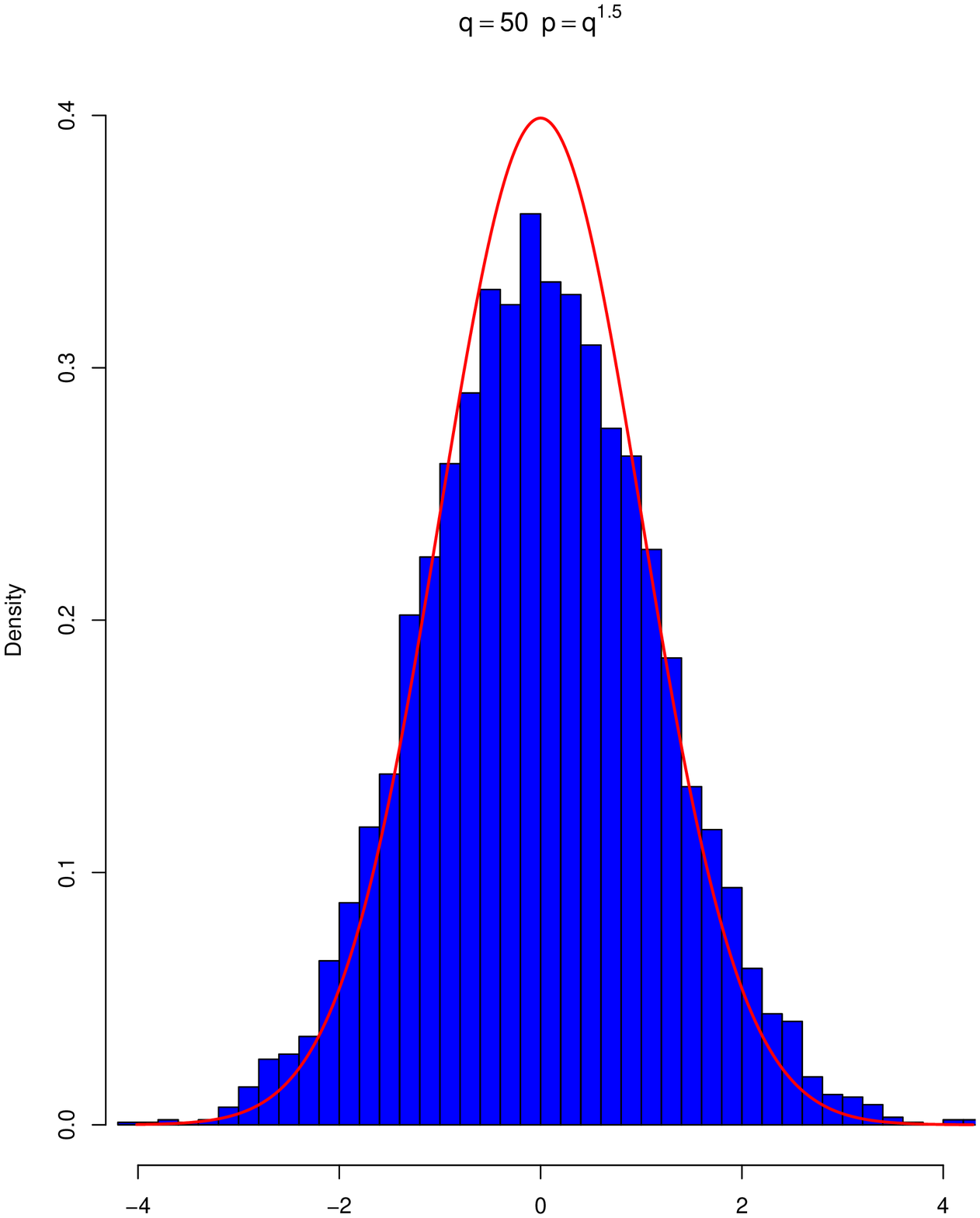} & \includegraphics[height=1.5in,width=1.7in]{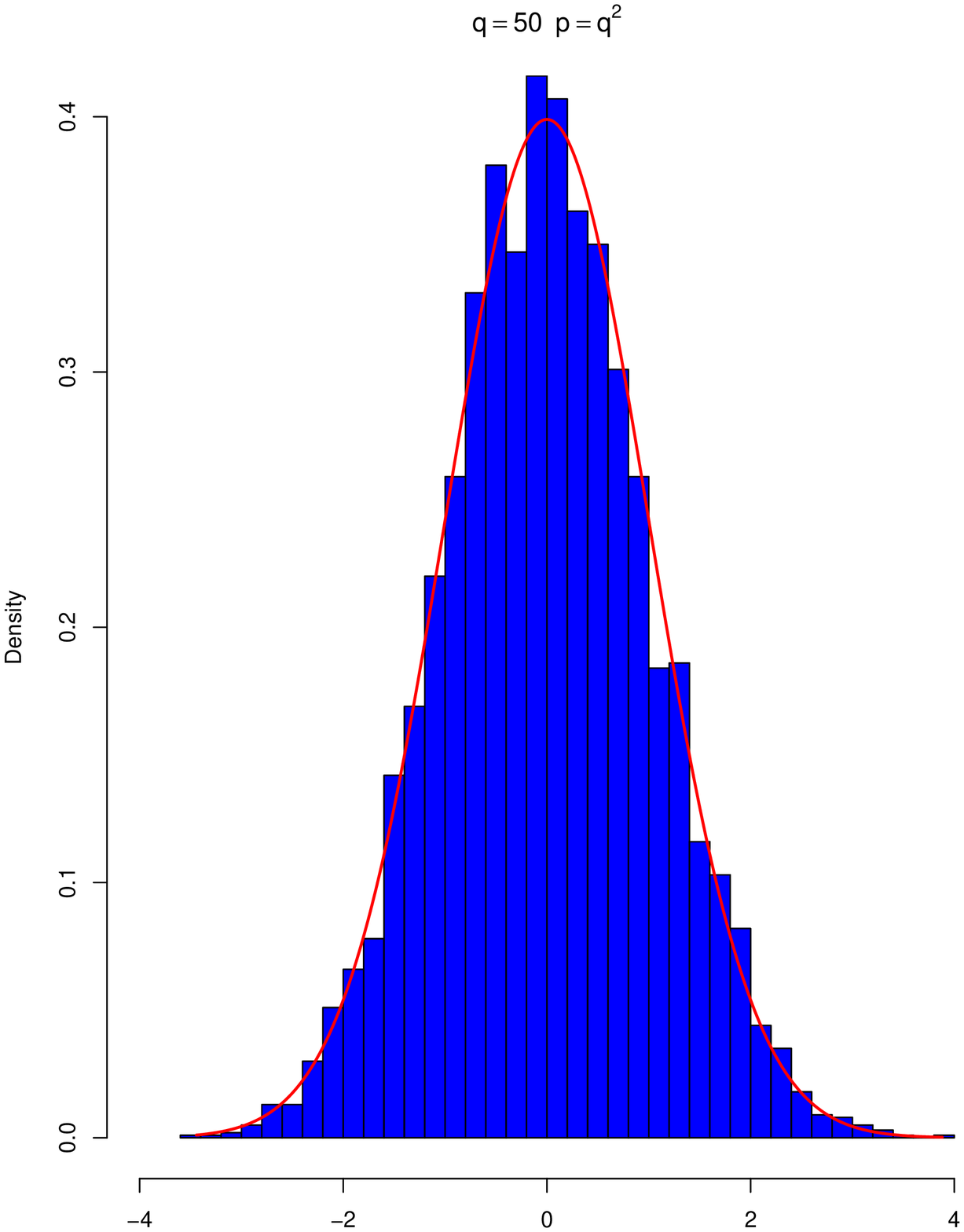} & \includegraphics[height=1.5in,width=1.7in]{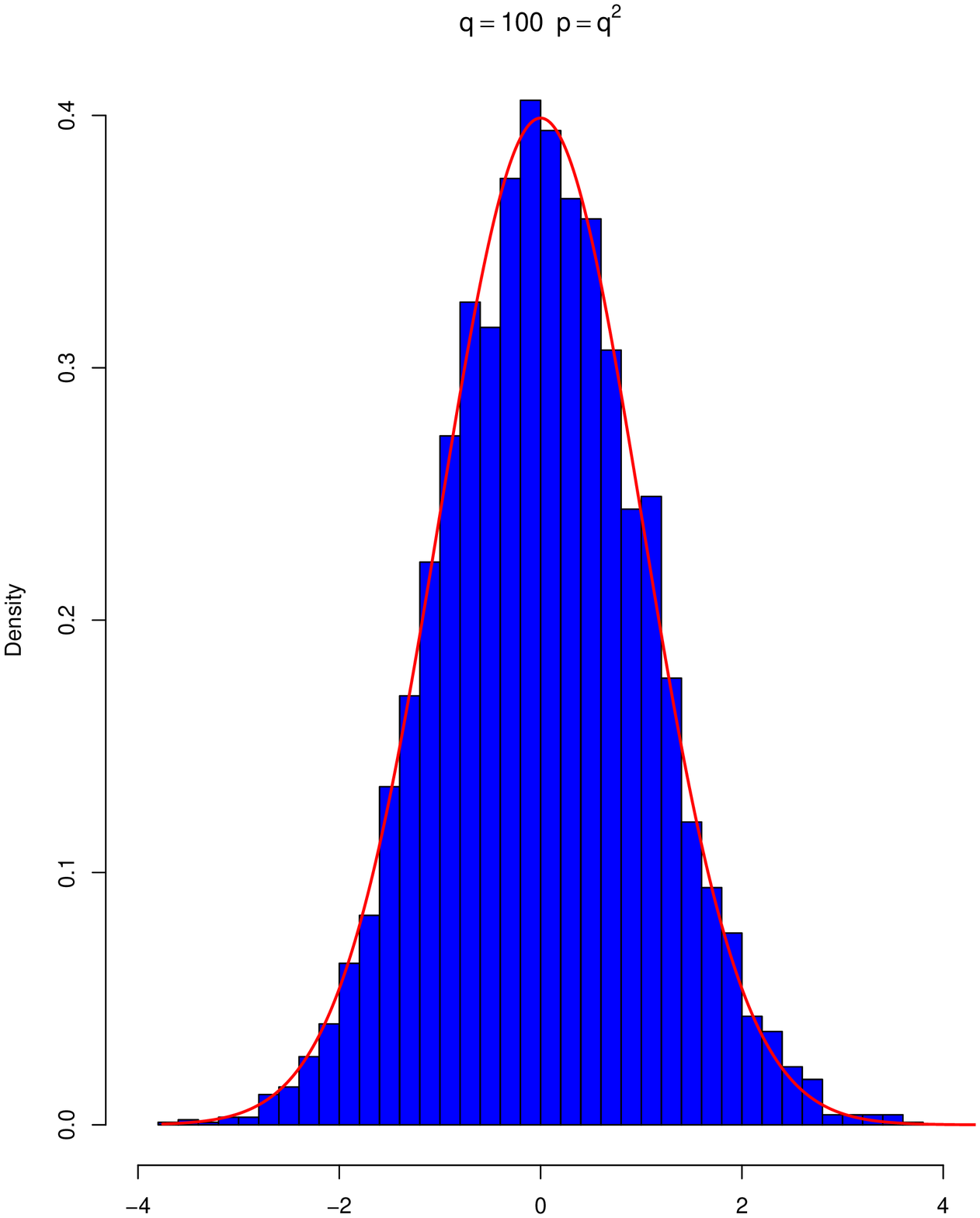}
\end{tabular}
\caption{\sl The plots simulate the density of $\frac{1}{2pq}\sum_{1\leq i\ne j\leq q}\big[(\bd{g}_i'\bd{g}_j)^2-p\big]$ in  Lemma ~\ref{awesome_pen} for $(p,q)= (165, 30), (900, 30), (1600, 40)$ in the top row, and $(p,q)= (355, 50), (2500, 50), (10000, 100)$ in the bottom row, respectively. They match the density curve of $N(0,1)$ better as both $p$ and $q$ are larger.}
\label{fig_Xinmei}
\end{center}
\end{figure}

Now we prove the second main result in this section.

\medskip
\noindent\textbf{Proof of Lemma \ref{awesome_pen}}. Since $\frac{q}{p}\to 0$, we assume that, without loss of generality, $q< p$ for all $n\geq 3$. Let $$T_n=\frac{1}{pq}\sum_{1\leq i\ne j\leq q}\big[(\bd{g}_i'\bd{g}_j)^2-p\big],$$
which can be rewritten by
 $$T_n=\frac{2}{pq}\sum_{j=2}^q \sum_{i=1}^{j-1}\big[(\bd{g}_i'\bd{g}_j)^2-p\big].$$
 Define $C_0=0$ and
 $$C_j=\sum_{i=1}^{j-1}\big[(\bd{g}_i'\bd{g}_j)^2-p\big]$$
 for $2\leq j \leq q$.
 It is easy to see
\bea\lbl{putin_three}
B_j:=\mathbb{E}[C_j|\mathcal{F}_{j-1}]=\sum_{i=1}^{j-1}
\mathbb{E}[\big((\bd{g}_i'\bd{g}_j)^2-p\big)|\mathcal{F}_{j-1}]=\sum_{i=1}^{j-1}(\|\bd{g}_i\|^2-p)
\eea
 where the sigma algebra $\mathcal{F}_k=\sigma\big(\bd{g}_1, \bd{g}_2, \cdots, \bd{g}_k\big)$ for all
 $1\le k\le q.$
Thus
\bea\lbl{siren_green}
X_j:=C_j-\mathbb{E}[C_j|\mathcal{F}_{j-1}]=
\sum_{i=1}^{j-1}\big[ (\bd{g}_i' \bd{g}_j)^2-\|\bd{g}_i\|^2\big], \quad 2\le j\le q,
\eea
form a martingale difference with respect to the $\sigma$-algebra $(\mathcal{F}_{j})_{2\le j\le q}.$
Therefore $T_n$ can be further written by
\be\lbl{subofT}T_n=\frac{2}{pq}\sum_{j=2}^q X_j+\frac{2}{pq}\sum_{j=2}^q B_j.\ee
By using (\ref{putin_three}) and changing sums, one gets
$$\frac{2}{pq} \sum_{j=2}^qB_j=\frac{2}{pq}\sum_{i=1}^{q-1}(\|\bd{g}_i\|^2-p)(q-i).$$
Since $\|\bd{g}_i\|^2\sim \chi^2(p)$ for each $1\leq i \leq q$, we know ${\rm Var}(\|\bd{g}_i\|^2-p)={\rm Var}(\|\bd{g}_i\|^2)=2p$. Hence,
\beaa
 {\rm Var}\Big(\frac{1}{pq} \sum_{j=2}^qB_j\Big) & = & \frac{1}{p^2q^2}
 \sum_{i=1}^{q-1}(q-i)^2 {\rm Var}(\|\bd{g}_i\|^2-p)\\
&\leq & \frac{2p}{p^2q^2}\cdot q^3=\frac{2q}{p}\to 0.
\eeaa
This together with the fact $\mathbb{E}\sum_{j=2}^qB_j=0$ indicates that
$$\frac{1}{pq} \sum_{j=2}^q B_j\stackrel{p}{\to} 0.$$
By (\ref{subofT}), to prove the theorem, it suffices to prove
$$\frac{1}{pq} \sum_{j=2}^q X_j \to N(0, 1)$$ weakly as $n\to\infty.$
By the Lindeberg-Feller central limit theorem for martingale differences
(see, for example, p. 414 from \cite{Durrett}), it is enough to verify that
 \be\lbl{squareX}  W_n:=\frac{1}{p^2q^2} \sum_{j=2}^q\mathbb{E}[X_j^2|\mathcal{F}_{j-1}]\stackrel{p}{\to} 1 \ee
 and
\be\lbl{fourX} \frac{1}{p^4q^4} \sum_{j=2}^q\mathbb{E}\big(X_j^4\big)\to 0\ee
as $n\to\infty.$ To prove \eqref{squareX}, it suffices to show
\bea\lbl{Weight_loss}
\mathbb{E}(W_n)\rightarrow 1
\eea
and
\bea\lbl{Weight_loss2}
 {\rm Var}(W_n)\rightarrow 0
\eea
as $n\to\infty.$ In the rest of the proof, due to their lengths we will show the above three assertions in the order of  (\ref{Weight_loss}),
(\ref{fourX}) and (\ref{Weight_loss2}), respectively. The proof will be
finished then.

{\it The proof of (\ref{Weight_loss})}. For simplicity, given $2\le j\le q,$ define
\beaa\lbl{fire_bad}
\bd{w}_i:=(\bd{g}_i'\bd{g}_j)^2-\|\bd{g}_i\|^2
\eeaa
for any $i=1, \cdots, j-1.$ Recall that $\bd{a}_i:=\frac{\bd{g}_i}{\|\bd{g}_i\|}$ and $\|\bd{g}_i\|$ are independent, we have a useful fact that
\bea\lbl{night_sun}
\bd{w}_i =\|\bd{g}_i\|^2[(\bd{a}_i'\bd{g}_j)^2-1]=\|\bd{g}_i\|^2(\chi-1)
\eea
given $\bd{g}_i$, where $\chi$ is a random variable with distribution $\chi^2(1)$ and is independent of $\|\bd{g}_i\|$. It is easy to see
$X_j=\sum_{i=1}^{j-1} \bd{w}_i$ from (\ref{siren_green}).
By Lemma \ref{soup_good}, for any $2\le j\le q,$ we see
\beaa\lbl{x2expe}
\mathbb{E}\big(X^2_j\big)
&=&\sum_{i=1}^{j-1}\mathbb{E}(\bd{w}_i^2)+\sum_{1\le i\neq k\le j-1}\mathbb{E}(\bd{w}_i\bd{w}_k)\\
&=&\sum_{i=1}^{j-1}\mathbb{E}\big[\mathbb{E}(\bd{w}_i^2 |\bd{g}_{i})\big]
+\sum_{1\le i\neq k<j}\mathbb{E}\big(\mathbb{E}[\bd{w}_i \bd{w}_k | \bd{g}_{j}]\big)\\
&=&2\sum_{1\le i<j}\mathbb{E}\|\bd{g}_i\|^4+\sum_{1\le i\neq k<j}\mathbb{E}(\|\bd{g}_j\|^2-p)^2\\
&=&2(p^2+2p)(j-1)+2p(j-1)(j-2),
\eeaa
where in the third equality we use the fact that  $\bd{w}_i$ and  $\bd{w}_k$ are conditionally independent given $\bd{g}_j$ for any $j \ne k$ and $\mathbb{E}[\bd{w}_i | \bd{g}_{j}]=\|\bd{g}_j\|^2-p.$
Thereby,
\beaa
\mathbb{E}(W_n)=\frac{1}{p^2q^2}\sum_{j=2}^q
\mathbb{E}(X_j^2)=\frac{(p^2+2p)q(q-1)+O(pq^3)}{p^2q^2}\to 1
\eeaa
as $n\to\infty$. This justifies (\ref{Weight_loss}). In particular,
\bea\lbl{where_moon}
\sum_{j=2}^q \mathbb{E}(X_j^2)\leq Cp^2q^2
\eea
for all $n\geq 4$, which will be used later.

{\it The proof of (\ref{fourX})}. Fix $j$ with $2\le j\le q$.
Observe that $(\bd{w}_i)_{1\le i\le j-1}$ form again a martingale difference with respect to
the sigma algebra $(\mathcal{F}_{i})_{1\le i\le j-1}.$ The Burkholder inequality
(see, for example, \cite{Shiryayev}) says that, for any
$s>1,$
$$ \mathbb{E}\Big(\sum_{i=1}^{j-1} \bd{w}_i\Big)^s\le C\cdot\mathbb{E}\Big(\sum_{i=1}^{j-1}\bd{w}_i^2\Big)^{s/2}, $$
where $C$ is a universal constant depending on $s$ only. By taking $s=4$, we see from $X_j=\sum_{i=1}^{j-1} \bd{w}_i$ that
\bea
\mathbb{E}(X^4_j)
& \le & C
\sum_{i=1}^{j-1}\mathbb{E}(\bd{w}_i^4)+C\sum_{1\le i\neq k\le j-1}\mathbb{E}(\bd{w}_i^2\bd{w}_l^2)\nonumber\\
&=& C\sum_{i=1}^{j-1}\mathbb{E}[\mathbb{E}(\bd{w}_i^4|\bd{g}_i)]+C\sum_{1\le i\neq k\le j-1}
\mathbb{E} \big[\mathbb{E}(\bd{w}_i^2| \bd{g}_j)\cdot \mathbb{E}(\bd{w}_k^2| \bd{g}_j)\big]
\lbl{kill_cloud}
\eea
by using the conditional independence. From (\ref{night_sun}),
\beaa
\mathbb{E}\big[\mathbb{E}(\bd{w}_i^4|\bd{g}_i)\big]
&\leq & \mathbb{E}(\|\bd{g}_i\|^8)\cdot\mathbb{E} \big[(\chi-1)^4\big]\\
& =& \mathbb{E}\big(\chi^2(p)^4)\cdot\mathbb{E} \big[(\chi-1)^4\big]\\
& \leq & Cp^4
\eeaa
by Lemma \ref{party_love}. Now, from (\ref{Trump}), the second sum in (\ref{kill_cloud}) is bounded by
\beaa
& & \mathbb{E} \big[3\|\bd{g}_j\|^4+(p^2+2p)-2(p+2)\|\bd{g}_j\|^2\big]^2\\
& \leq & 3\big[\mathbb{E}(\chi^2(p)^4)+(p^2+2p)^2 + 4(p+2)^2\mathbb{E}(\chi^2(p)^2)\big]\\
& \leq & Cp^4
\eeaa
by the Minkowski inequality and Lemma \ref{party_love} again. The above two estimates together with (\ref{kill_cloud}) imply
\beaa
\mathbb{E}(X^4_j) \leq C\cdot\big[(j-1)p^4 + (j-1)(j-2)p^4]
\eeaa
for $2\le j\le q$, where $C$ is free of $n$ and $j$. Consequently,

$$\frac{1}{p^4q^4}\sum_{j=2}^q \mathbb{E}\big(X_j^4\big)\le
\frac{C}{p^4q^4}\sum_{j=2}^qj^2p^4
\leq \frac{C p^4q^3}{p^4q^4}\to 0$$
as $n\to\infty.$ This concludes (\ref{fourX}).

{\it The proof of (\ref{Weight_loss2})}.
 We need to prove that
\beaa
\frac{1}{p^4q^4}{\rm Var}\bigg[ \sum_{j=2}^q\mathbb{E}\big(X_j^2|\mathcal{F}_{j-1}\big)\bigg]
\to 0
\eeaa
as $n\to\infty.$
Let us first compute $\mathbb{E}[X_j^2|\mathcal{F}_{j-1}]$.
Set $\alpha_i=\frac{\bd{g}_i}{\|\bd{g}_i\|}$ for $1\le i\le q.$ Then,
\beaa
\mathbb{E}\big(X_j^2|\mathcal{F}_{j-1}\big)
&=&\mathbb{E}
\Big[\big(\sum_{i=1}^{j-1} \bd{w}_i\big)^2|\mathcal{F}_{j-1}\Big]\\
&=&\sum_{i=1}^{j-1}\mathbb{E}\big(\bd{w}_i^2|\mathcal{F}_{j-1}\big)+
I(j\geq 3)\cdot\sum_{1\le i\neq k\le j-1}\mathbb{E}\big(\bd{w}_i\bd{w}_k|\mathcal{F}_{j-1}\big) \\
&=&\sum_{i=1}^{j-1}\|\bd{g}_i\|^4\cdot \mathbb{E}[(\chi-1)^2]+I(j\geq 3)\cdot
\sum_{1\le i\neq k\le j-1}\|\bd{g}_i\|^2\|\bd{g}_k\|^2\cdot\\
&&\quad\quad\quad\quad\quad\quad\quad
\mathbb{E}\big[(\alpha_i' \bd{g}_j)^2(\alpha_k' \bd{g}_j)^2-(\alpha_i' \bd{g}_j)^2-
(\alpha_k' \bd{g}_j)^2+1|\mathcal{F}_{j-1}\big],
\eeaa
where $I(j\geq 3)$ is the indicator function of the set $\{j\geq 3\}.$
Given $\mathcal{F}_{j-1}$, evidently $\alpha_i' \bd{g}_j\sim N(0,1)$ for $1\leq i \leq j-1.$
Therefore, from Lemma \ref{doublegreat} we have
\beaa
\mathbb{E}\big(X_j^2|\mathcal{F}_{j-1}\big)&=&2\sum_{i=1}^{j-1}\|\bd{g}_i\|^4+2I(j\geq 3)\cdot\sum_{1\le i\neq k\le j-1}
\|\bd{g}_i\|^2\|\bd{g}_k\|^2(\alpha_i'\alpha_k)^2\\
&=&2\sum_{i=1}^{j-1}\|\bd{g}_i\|^4+2I(j\geq 3)\cdot\sum_{1\le i\neq k\le j-1}(\bd{g}_i'\bd{g}_k)^2.
\eeaa
By changing the order of sums, it is not difficult to verify that
$$ \sum_{j=2}^q\mathbb{E}\big(X_j^2|\mathcal{F}_{j-1}\big)
=2\sum_{i=1}^{q-1}\|\bd{g}_i\|^4(q-i)+
4\sum_{i=1}^{q-2}\sum_{k=i+1}^{q-1}(\bd{g}_i'\bd{g}_k)^2(q-k).$$
Therefore,
\bea & & {\rm Var}\Big[\sum_{j=2}^q\mathbb{E}\big(X_j^2|\mathcal{F}_{j-1}\big)\Big]\nonumber\\
&\le & 8\cdot{\rm Var}\Big[\sum_{i=1}^{q-1}\|\bd{g}_i\|^4(q-i)\Big]+32\cdot
 {\rm Var}\bigg[\sum_{i=1}^{q-2}\sum_{j=i+1}^{q-1}(\bd{g}_i'\bd{g}_j)^2(q-j) \bigg]\nonumber\\
&=& 8\cdot\sum_{i=1}^{q-1}(q-i)^2 {\rm Var}(\|\bd{g}_i\|^4)+32\cdot {\rm Var}
\bigg[\sum_{j=2}^{q-1}(q-j)\sum_{i=1}^{j-1}(\bd{g}_i'\bd{g}_j)^2\bigg].\lbl{I+II}
\eea
On the one hand, by Lemma \ref{party_love} we know \be\lbl{I-I}
\sum_{i=1}^{q-1}(q-i)^2 {\rm Var}(\|\bd{g}_i\|^4)\leq q^3\cdot{\rm Var}\big[\chi^2(p)^2\big]
= 8pq^3(p+2)(p+3).
\ee
Moreover, for $2\le j\le q$ fixed, recall the notation
$$\bd{w}_i:=(\bd{g}_i'\bd{g}_j)^2-\|\bd{g}_i\|^2, \quad X_j=\sum_{i=1}^{j-1}\bd{w}_i.$$
Then  $$\sum_{i=1}^{j-1} (\bd{g}_i'\bd{g}_j)^2=
\sum_{i=1}^{j-1}\big(\bd{w}_i+\|\bd{g}_i\|^2\big)=X_j+\sum_{i=1}^{j-1}\|\bd{g}_i\|^2,$$
which implies
\beaa
\sum_{j=2}^{q-1}(q-j)\sum_{i=1}^{j-1}(\bd{g}_i'\bd{g}_j)^2
&=& \sum_{j=2}^{q-1}(q-j) X_j+\sum_{j=2}^{q-1}(q-j)\sum_{i=1}^{j-1}\|\bd{g}_i\|^2\\
&=& \sum_{j=2}^{q-1}(q-j) X_j+\sum_{i=1}^{q-2}\|\bd{g}_i\|^2\sum_{j=i+1}^{q-1}(q-j)\\
&=& \sum_{j=2}^{q-1}(q-j) X_j+\frac{1}{2}\sum_{i=1}^{q-2}\|\bd{g}_i\|^2(q-i-1)(q-i).
\eeaa
 Consequently we have
\bea
& & {\rm Var}\bigg[\sum_{j=2}^{q-1}(q-j)\sum_{i=1}^{j-1}(\bd{g}_i'\bd{g}_j)^2\bigg]\nonumber\\
&\le & 2\cdot {\rm Var}\bigg[\sum_{j=2}^{q-1}(q-j) X_j\bigg]+\frac{1}{2}\cdot{\rm Var}
\bigg[\sum_{i=1}^{q-2}\|\bd{g}_i\|^2(q-i-1)(q-i)\bigg] \nonumber\\
&=&2\sum_{j=2}^{q-1} (q-j)^2{\rm Var}(X_j)+
\frac12\sum_{i=1}^{q-2}(q-i-1)^2(q-i)^2{\rm Var}\big(\|\bd{g}_1\|^2\big)\lbl{mouse_day}\\
&\le & 2 q^2\sum_{j=2}^{q}\mathbb{E}(X_j^2)+ q^5\,{\rm Var}\big(\|\bd{g}_1\|^2\big)
\lbl{California}\\
&=& C\big(p^2q^4+pq^5),\lbl{miracle_2}
\eea
where we use the fact $(X_j)_{2\le j\le q}$ is a martingale with respect
to $(\mathcal{F}_j)_{2\le j\le q}$ in (\ref{mouse_day}); the trivial bounds $(q-j)^2\leq q^2$
and $(q-i-1)^2(q-i)^2\leq q^4$ are applied in (\ref{California});
the last step is obtained by (\ref{where_moon}) and
the identity ${\rm Var}\big(\|\bd{g}_1\|^2\big)=2p.$

Plugging \eqref{I-I} and \eqref{miracle_2} into \eqref{I+II},
we have from the fact $q\leq p$ that
\beaa
\frac{1}{p^4q^4}{\rm Var}\bigg( \sum_{j=2}^q\mathbb{E}[X_j^2|\mathcal{F}_{j-1}]\bigg)
\leq \frac{Cp^3q^3}{p^4q^4}=\frac{C}{pq}\to 0
\eeaa
as $n\to\infty$. This finishes the verification of (\ref{Weight_loss2}).
The proof is completed.
\hfill$\blacksquare$

\medskip

Now we prove  Lemma \ref{var-e}.

\noindent{\bf Proof of Lemma \ref{var-e}}. Write ${\bf X}=({\bf g}_1, {\bf g}_2, \cdots, {\bf g}_q)$ where ${\bf g}_i\sim N_p(0, {\bf I}_p).$ A repeatedly used fact is that $\|{\bf g}_i\|^2\sim \chi^2(p)$ for each $i$. Using this fact, independence,  Lemma \ref{party_love} and (\ref{night_sun}),  we have
$$\aligned &\mathbb{E}({\bf g}_1'{\bf g}_2)^2=p, \quad \mathbb{E}({\bf g}_1'{\bf g}_2)^4=3\mathbb{E}(\|{\bf g}_1\|^4)=3p(p+2);\\
& \mathbb{E}[\|{\bf g}_1\|^2\cdot ({\bf g}^{\prime}_1 {\bf g}_2)^2]=p(p+2);\\
& \mathbb{E}\big[({\bf g}_1'{\bf g}_2)^2({\bf g}_1'{\bf g}_3)^2\big]=\mathbb{E}(\|{\bf g}_1\|^4)=p(p+2);\\
&\mathbb{E}\big[\|{\bf g}_1\|^4({\bf g}_1'{\bf g}_2)^2\big]=\mathbb{E}(\|{\bf g}_1\|^6)=p(p+2)(p+4).\endaligned $$
Easily, ${\rm tr}({\bf X}^{\prime}{\bf X})=\sum_{i=1}^q\|{\bf g}_i\|^2$ and
\be\lbl{trace-one}{\rm tr}[({\bf X'X})^2]=\sum_{1\le i\neq j\le q}({\bf g}_i'{\bf g}_j)^2+\sum_{k=1}^q\|{\bf g}_k\|^4.\ee
(i)
By using the set of formulas right before (\ref{trace-one}), we obtain
\be\lbl{power2}\mathbb{E}{\rm tr}[({\bf X'X})^2]=q(q-1)p+qp(p+2)=pq(p+q+1).\ee
From \eqref{trace-one} we see
\be\lbl{trace-two}
 {\rm tr}^2[({\bf X'X})^2]=\big(\sum_{k=1}^q\|{\bf g}_k\|^4\big)^2+\big[\sum_{1\le i\neq j\le q}({\bf g}_i'{\bf g}_j)^2\big]^2
+2\sum_{k=1}^q\|{\bf g}_k\|^4\sum_{1\le i\neq j\le q}({\bf g}_i'{\bf g}_j)^2.
\ee
It is easy to check
\bea\lbl{thin_of}
\Big(\sum_{k=1}^q\|{\bf g}_k\|^4\Big)^2=\sum_{k=1}^q\|{\bf g}_k\|^8+\sum_{1\le i\neq j\le q}\|{\bf g}_i\|^4\|{\bf g}_j\|^4
\eea
and
\bea
\Big[\sum_{1\le i\neq j\le q}({\bf g}_i'{\bf g}_j)^2\Big]^2 &=&\sum_{1\le i\neq j\neq k\neq l\le q} ({\bf g}_i'{\bf g}_j)^2({\bf g}_k'{\bf g}_l)^2+2\sum_{1\le i\neq j\le q}({\bf g}_i'{\bf g}_j)^4\nonumber\\
&&+4\sum_{1\le i\neq j\neq k\le q}({\bf g}_i'{\bf g}_j)^2({\bf g}_i'{\bf g}_k)^2\lbl{head_foot}
\eea
and
\bea\lbl{dudu}
2\sum_{k=1}^q\|{\bf g}_k\|^4\sum_{1\le i\neq j\le q}({\bf g}_i'{\bf g}_j)^2=2\sum_{1\le i\neq j\neq k\le q}\|{\bf g}_k\|^4({\bf g}_i'{\bf g}_j)^2+4\sum_{1\le i\neq j\le q}\|{\bf g}_i\|^4({\bf g}_i'{\bf g}_j)^2.
\eea
By the formulas right before (\ref{trace-one}), Lemma \ref{party_love}, (\ref{thin_of}), (\ref{head_foot}) and (\ref{dudu}), respectively,   we have that
$$
\mathbb{E}\Big(\sum_{k=1}^q\|{\bf g}_k\|^4\Big)^2=qp(p+2)(p+4)(p+6)+q(q-1)p^2(p+2)^2$$
and
\beaa
\mathbb{E}\Big[\sum_{1\le i\neq j\le q}({\bf g}_i'{\bf g}_j)^2\Big]^2&=&q(q-1)(q-2)(q-3)p^2+6p(p+2)q(q-1)\\
&&+4p(p+2)q(q-1)(q-2)
\eeaa
and
\beaa
2\mathbb{E}\Big[\sum_{k=1}^q\|{\bf g}_k\|^4\sum_{1\le i\neq j\le q}({\bf g}_i'{\bf g}_j)^2\Big]&=& 2q(q-1)(q-2)p^2(p+2)\\
 &&+ 4q(q-1)p(p+2)(p+4).
\eeaa
Putting all these three expressions back into \eqref{trace-two}, one gets
\beaa
 \mathbb{E}{\rm tr}^2[({\bf X'X})^2]&=&qp(p+2)(p+4)(p+6)+q(q-1)p^2(p+2)^2\\
&&+q(q-1)(q-2)(q-3)p^2+6p(p+2)q(q-1)\\
&&+4p(p+2)q(q-1)(q-2)\\
 &&+2q(q-1)(q-2)p^2(p+2)+4q(q-1)p(p+2)(p+4).
\eeaa
This together with \eqref{power2} implies
$$
{\rm Var}\big({\rm tr}[({\bf X'X})^2]\big)
=4p^2q^2+8pq(p+q)^2+20pq(p+q+1).
$$

(ii) Recall (\ref{trace-one}). By independence,
\beaa
 {\rm Cov}\big({\rm tr}({\bf X}^{\prime}{\bf X}), {\rm tr}[({\bf X}^{\prime}{\bf X})^2]\big)&=&{\rm Cov}\Big(\sum_{i=1}^q\|{\bf g}_i\|^2, \sum_{i=1}^q\|{\bf g}_i\|^4+\sum_{1\le i\neq j\le q}({\bf g}^{\prime}_i {\bf g}_j)^2 \Big)\\
&=&\sum_{i=1}^q{\rm Cov}\big(\|{\bf g}_i\|^2, \|{\bf g}_i\|^4\big)+2\sum_{1\le i\neq j\le q}{\rm Cov}(\|{\bf g}_i\|^2,  ({\bf g}^{\prime}_i {\bf g}_j)^2)\\
&=&\sum_{i=1}^q\big[\mathbb{E}(\|{\bf g}_i\|^6)-\mathbb{E}(\|{\bf g}_i\|^4)\cdot\mathbb{E}(\|{\bf g}_i\|^2)\big]\\
&&+2\sum_{1\le i\neq j\le q}\big(\mathbb{E}[\|{\bf g}_i\|^2\cdot ({\bf g}^{\prime}_i {\bf g}_j)^2]-\mathbb{E}(\|{\bf g}_i\|^2)\cdot\mathbb{E}({\bf g}^{\prime}_i {\bf g}_j)^2\big).
\eeaa
From the formulas right before (\ref{trace-one}) again, we have
$$\aligned {\rm Cov}\big({\rm tr}[{\bf X}^{\prime}{\bf X}], {\rm tr}[({\bf X}^{\prime}{\bf X})^2]\big)&=\sum_{i=1}^q[p(p+2)(p+4)-p^2(p+2)]+2\sum_{1\le i\neq j\le q}[p(p+2)-p^2]\\
&=4pq(p+2)+4pq(q-1)\\
&=4pq(p+q+1).\endaligned $$
The proof is completed now.
\hfill$\blacksquare$

\section{This part is for referees only}\lbl{only_for_referee}

\noindent\textbf{Proof of Lemma \ref{Low_beauty}}.
Suppose the conclusion is not true, that is, $\liminf_{n\to\infty}f_n(p_n, q_n)=0$ for some  sequence $\{(p_n, q_n)\}_{n=1}^{\infty}$ with $1\leq p_n, q_n\leq n$ for each $n\geq 1$ and  $\lim_{n\to\infty}(p_nq_n^2)/n=\alpha\in (0, \infty)$.
Then there exists a subsequence $\{n_k; \, k\geq 1\}$ satisfying  $1\leq p_{n_k}, q_{n_k}\leq n_k$ for all $k\geq 1$ and $\lim_{k\to\infty}(p_{n_k}q_{n_k}^2)/n_k=\alpha\in (0, \infty)$ such that
\bea\lbl{shawn_23}
\lim_{k\to\infty}f_{n_k}(p_{n_k}, q_{n_k})=0.
\eea
There are two possibilities: $\liminf_{k\to\infty}q_{n_k}<\infty$ and $\liminf_{k\to\infty}q_{n_k}=\infty$. Let us discuss the two cases separately.

{\it (a)}. Assume $\liminf_{k\to\infty}q_{n_k}<\infty$. Then there exists a further subsequence $\{n_{k_j}\}_{j=1}^{\infty}$ such that $q_{n_{k_j}}\equiv m\geq 1$. For convenience of notation, write $\bar{n}_j=n_{k_j}$ for all $j\geq 1.$ The condition $\lim_{n\to\infty}\frac{p_nq_n^2}{n}=\alpha$ implies that $\lim_{j\to\infty}\frac{p_{\bar{n}_j}}{\bar{n}_j}=\frac{\alpha}{m^2}\in (0, 1].$ By (\ref{shawn_23}) and the monotonocity,
\bea\lbl{privacy_it_51}
\lim_{j\to\infty}f_{\bar{n}_j}(p_{\bar{n}_j}, 1)=\lim_{j\to\infty}f_{\bar{n}_j}(p_{\bar{n}_j}, q_{\bar{n}_j})=0.
\eea
Define $\tilde{p}_{\bar{n}_j}=[p_{\bar{n}_j}/2]+1$ for all $j\geq 1.$ Then, $\lim_{j\to\infty}\frac{\tilde{p}_{\bar{n}_j}}{\bar{n}_j}=c:=\frac{\alpha}{2m^2}\in (0, \frac{1}{2}].$ Construct a new sequence such that
\beaa\lbl{green_beer}
\tilde{p}_r=
\begin{cases}
\tilde{p}_{\bar{n}_j}, & \text{if $r=\bar{n}_j$ for some $j\geq 1$};\\
[cr] \vee 1, & \text{if not}
\end{cases}
\eeaa
and $\tilde{q}_r=1$ for $r=1,2, \cdots.$ Obviously, $\tilde{p}_{\bar{n}_j}\leq p_{\bar{n}_j}$ for each $j\geq 1.$ It is easy to check $1\leq  \tilde{p}_r, \tilde{q}_r\leq r$ for all $r\geq 1$ and $\lim_{r\to\infty}\tilde{p}_r/r=c\in (0,1/2)$. So $\{(\tilde{p}_r, \tilde{q}_r);\, r\geq 1\}$ satisfies condition (i), and hence  $\liminf_{r\to\infty}f_r(\tilde{p}_r, \tilde{q}_r)>0$ by (\ref{red_red}). This contradicts (\ref{privacy_it_51}) since $f_r(\tilde{p}_r, \tilde{q}_r)=f_{\bar{n}_j}(\tilde{p}_{\bar{n}_j}, 1)\leq f_{\bar{n}_j}(p_{\bar{n}_j}, 1)$ if  $r=\bar{n}_j$ for some $j\geq 1$ by monotonocity.

{\it (b).} Assume $\liminf_{k\to\infty}q_{n_k}=\infty$. Then  $\lim_{k\to\infty}q_{n_k}=\infty$.
%
%
Define
\beaa
\tilde{p}_r=
\begin{cases}
p_{n_k}, & \text{if $r=n_k$ for some $k\geq 1$};\\
[r^{1/3}], & \text{if not}
\end{cases}
\eeaa
and
\beaa
\tilde{q}_r=
\begin{cases}
q_{n_k}, & \text{if $r=n_k$ for some $k\geq 1$};\\
\big([\sqrt{\alpha} r^{1/3}] + 1\big)\wedge r, & \text{if not}.
\end{cases}
\eeaa
 Trivially, $1\leq \tilde{p}_r, \tilde{q}_r\leq  r$ for all $r\geq 1$, $\lim_{r\to\infty}\tilde{q}_r=\infty$ and $\lim_{r\to\infty}\tilde{p}_r\tilde{q}_r^2/r=\alpha$.
By (ii),
\beaa
\liminf_{k\to\infty}f_{n_k}(p_{n_k}, q_{n_k})\geq \liminf_{r\to\infty}f_r(\tilde{p}_r, \tilde{q}_r)>0
\eeaa
since $\tilde{p}_r=p_{n_k}$ and $\tilde{q}_r=q_{n_k}$ if $r=n_k$.
This contradicts (\ref{shawn_23}).

In summary, each of the cases that $\liminf_{k\to\infty}q_{n_k}<\infty$ and that $\liminf_{k\to\infty}q_{n_k}=\infty$ results with a contradiction. Therefore, we obtain our desired conclusion.
\hfill$\blacksquare$\\

\noindent\textbf{Acknowledgement}. We thank Professor Xinmei Shen for very helpful  communications. In particular we thank her for producing Figure \ref{fig_Xinmei} for us.

\end{document}